\documentclass[preprint,11pt]{elsarticle}
\usepackage[toc,page]{appendix}
\usepackage{amsmath,amsthm,amsfonts,mathtools,bm,hyperref}

\usepackage{multirow}
\usepackage{tikz,pgfplots}
\usepgfplotslibrary{groupplots}

\usetikzlibrary{patterns}
\definecolor{BurntOrange}{rgb}{0.8, 0.33, 0.0}
\usepackage{subcaption}
\pgfplotsset{compat=newest}
\usepackage{float,indentfirst}
\usepackage{pgfplots}
\usepackage{booktabs}
\usepackage[margin=1in]{geometry}

\theoremstyle{plain}

\usepackage[capitalise]{cleveref}
\usepackage{appendix}
\usepackage{multicol}
\usepackage{nomencl}
\usepackage{ifthen}

\usepackage{times}
\usepackage{fancyhdr,graphicx,amsmath,amssymb}
\usepackage[ruled,vlined]{algorithm2e}

\usepackage{todonotes}
\usepackage{url}

\newcommand{\bt}[0]{\bm{\theta}}
\newcommand{\bx}{\mathbf{x}}
\newcommand{\bn}[0]{\mathbf{n}}
\newcommand{\bu}{\mathbf{u}}
\newcommand{\bv}[0]{\mathbf{v}}

\newcommand{\bsigma}[0]{\bm{\sigma}}
\newcommand{\bepsilon}[0]{\bm{\epsilon}}

\makenomenclature
\usepackage{etoolbox}
\renewcommand\nomgroup[1]{%
  \item[\bfseries
  \ifstrequal{#1}{A}{Symbols}{%
  \ifstrequal{#1}{B}{Superscripts}{%
  \ifstrequal{#1}{C}{Subscripts}{
  \ifstrequal{#1}{D}{Operators}{}}}}%
]}

\makeatletter
\def\ps@pprintTitle{%
   \let\@oddhead\@empty
   \let\@evenhead\@empty
   \def\@oddfoot{\reset@font\hfil\thepage\hfil}
   \let\@evenfoot\@oddfoot
}
\makeatother

\begin{document}

\begin{abstract}

We propose a novel approach to model viscoelasticity materials using neural networks, which capture rate-dependent and nonlinear constitutive relations. However, inputs and outputs of the neural networks are not directly observable, and therefore common training techniques with input-output pairs for the neural networks are inapplicable. To that end, we develop a novel computational approach to both calibrate parametric and learn neural-network-based constitutive relations of viscoelasticity materials from indirect displacement data in the context of multi-physics interactions. We show that limited displacement data hold sufficient information to quantify the viscoelasticity behavior. We formulate the inverse computation---modeling viscoelasticity properties from observed displacement data---as a PDE-constrained optimization problem and minimize the error functional using a gradient-based optimization method. The gradients are computed by a combination of automatic differentiation and physics constrained learning. The effectiveness of our method is demonstrated through numerous benchmark problems in geomechanics and porous media transport. 
\end{abstract}


\begin{frontmatter}

\title{Inverse Modeling of Viscoelasticity Materials using Physics Constrained Learning}
\author[rvt1]{Kailai~Xu}
\ead{kailaix@stanford.edu}

\author[rvt3]{Alexandre M. Tartakovsky}
\ead{alexandre.tartakovsky@pnnl.gov}

\author[rvt3]{Jeff Burghardt}
\ead{jeffrey.burghardt@pnnl.gov}

\author[rvt1,rvt2]{Eric~Darve}
\ead{darve@stanford.edu}

\address[rvt1]{Institute for Computational and Mathematical Engineering,
               Stanford University, Stanford, CA, 94305}
\address[rvt2]{Mechanical Engineering, Stanford University, Stanford, CA, 94305}
\address[rvt3]{Pacific Northwest National Laboratory, 902 Battelle Blvd, Richland, WA 99354}

\end{frontmatter}

\section{Introduction}

Viscoelastic material exhibits both viscous fluid and elastic solid characteristics \cite{sorvari2009modelling}. Examples of viscoelasticity materials include skins, memory foam mattress, wood, polymers and concrete \cite{xu2016inverse}.  Viscoelasticity plays a crucial role because all materials exhibit viscoelastic behavior to some extent. Therefore, it is very important to model viscoelasticity properties.
 The viscoelasticity makes the behavior of material rate-dependent; that is, the material's response to deformation or force may change over time. This rate-dependent behavior makes viscosity modeling unique and challenging compared to its counterparts, such as elasticity.

Inverse modeling of viscoelasticity uses observed data to infer the viscoelastic constitutive relation of a material. Typically, the strain and stress tensors in the constitutive relations are not directly observable. Therefore, we need to consider the conservation laws---usually described by partial differential equations (PDEs)---and conduct inverse computation using indirect data, such as displacement data from sensors located in a material. The observed data can also be a result of multi-physics interactions, such as those in coupled geomechanics and multi-phase flow models. 

In general, there are two types of hypothetical constitutive relations:  
\begin{enumerate}
	\item Parametric models. In parametric models, the forms of the constitutive relations are known but the physical parameters in the models are unspecified. This setting leads to a parameter inverse problem.
	\item Nonparametric models. In these models, the forms of the constitutive relations are unknown or only partially known (for example, a physical parameter in a parametric model may be a \textit{function} of state variables). This setting leads to a function inverse problem.
\end{enumerate}
 In the first case, the inverse modeling problem can be formulated as a PDE-constrained optimization problem \cite{tarantola2005inverse,nakamura2015inverse,biegler2003large,rees2010optimal}, where the objective function (also known as the \textit{loss function}) measures the discrepancy between the predicted and observed data, the conservations laws are the constraints, and the unknown physical parameters are the free optimization variables. However, this optimization formulation is inapplicable to the function inverse problem because the function space is infinite dimensional. We propose deep neural networks as surrogate models for the unknown functions because deep neural networks are promising to approximate high dimensional and nonsmooth functions \cite{weinan2017deep,han2018solving,han2017overcoming,xu2019neural}, which are quite challenging for traditional basis functions, such as the radial basis functions and piecewise linear functions.  Neural network approaches have been demonstrated effective in constitutive modeling in much recent work \cite{tartakovsky2018learning,2019arXiv190512530H,raissi2019physics}.

Inverse modeling for viscoelasticity has been studied extensively in science and engineering. Some applicable examples: \cite{catheline2004measurement} uses  transient elastography to measure the viscoelastic properties of a homogeneous soft solid from experimental data;
\cite{janno1997inverse} proposes a Fourier's method of eigenfunction expansion to identify memory kernels in linear viscoelasticity;
\cite{kim2005characterization} develops a three dimensional finite element model to simulate the forces at the indenter and uses the Levenenberg-Marquardt method to update new parameters in the viscoelastic and hyperelastic materials;
\cite{brigham2007inverse} uses measured acoustic pressure amplitudes in the fluid to quantify both elastic and viscoelastic material behaviors. Particularly, \cite{brigham2007inverse} combined a surrogate model optimization strategy and a classical random search technique for finding global extrema;
\cite{barkanov2009characterisation} proposes a new inverse technique for characterizing the nonlinear mechanical properties of viscoelastic core layers in sandwich panels by planning the experiments and  response surface techniques in order to minimize an error functional;
\cite{pagnacco2007inverse} determines elasticity and viscoelasticity parameters of an isotropic material by a hybrid numerical/experimental approach using full-field measurements. However, existing methods generally focus on a specific parametric viscoelasticity model and do not consider the coactions of multi-physics, partially due to difficulties in extracting data for fitting parametric viscoelasticity models and lack of an expressive model suitable for fitting a vast of viscoelasticity models. 

\begin{figure}[htpb]
\centering
  \includegraphics[width=0.6\textwidth]{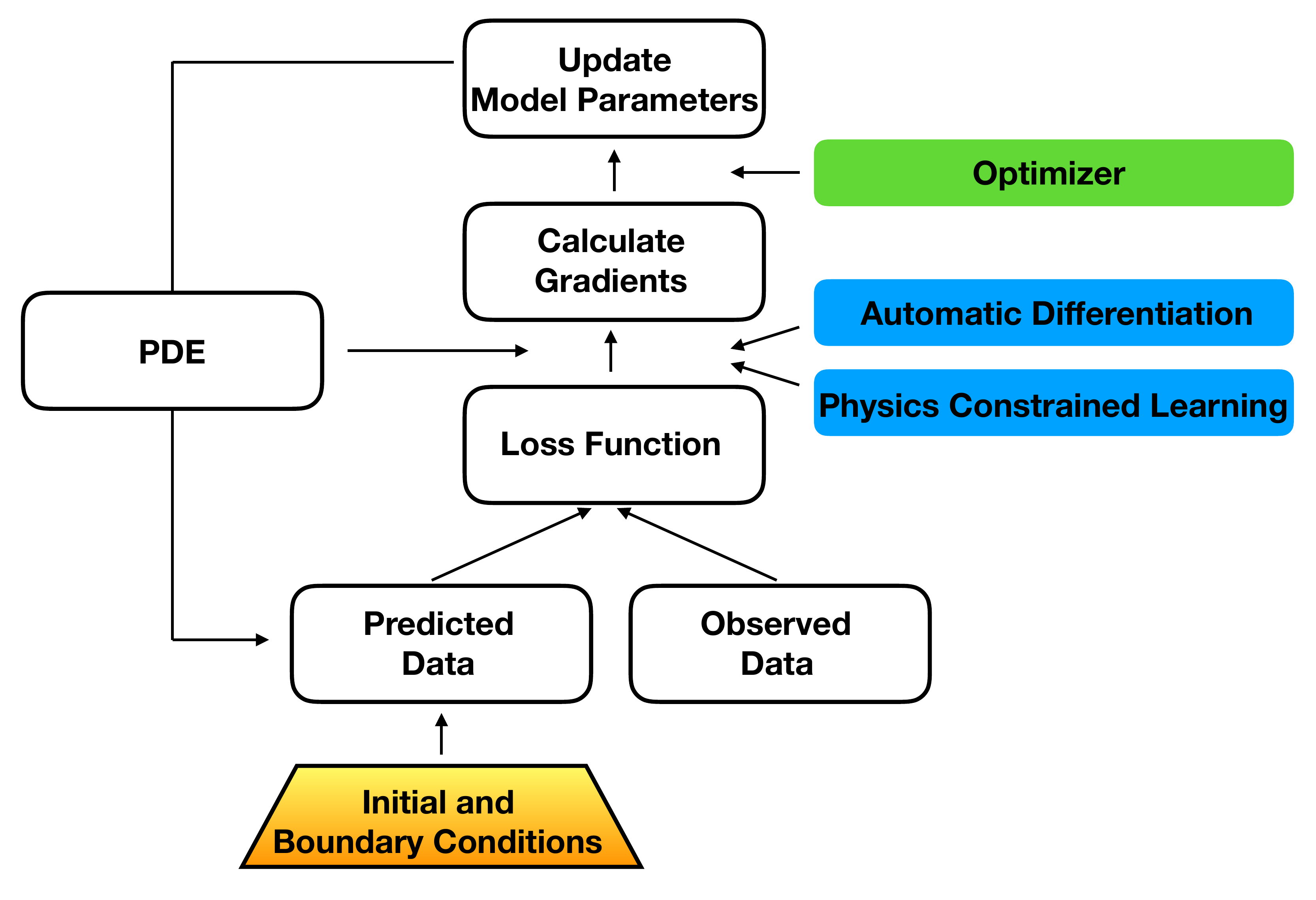}
  \caption{Flow chart for inverse modeling procedure. The gradients of the loss function are computed using AD or physics constrained learning. These gradients are used in an off-the-shelf optimizer for updating the parameters, i.e., physical parameters or neural network weights and biases. The outer loop repeats until the discrepancy between the predicted and observed data are sufficiently close.}
  \label{fig:routine}
\end{figure}

The goal of this research is to develop a novel inverse modeling approach for calibrating parametric viscoelastic models and learning neural-network-based constitutive relations using reverse-mode automatic differentiation (AD) techniques \cite{paszke2017automatic, baydin2017automatic} in the context of multi-physics interactions \cite{zimmerman2006multiphysics}. AD is used to compute the gradients of the objective function in the PDE-constrained optimization problem. 
With the advent of deep learning, plenty of highly efficient, robust and scalable AD software have been developed in the context of machine learning and scientific computing \cite{abadi2016tensorflow,paszke2017automatic,hogan2017adept,revels2016forward,phipps2009sacado}. However, the plain AD technique can only calculate the derivatives of explicit functions whose derivatives can be analytically computed, or compositions of explicit functions. For implicit functions, such as the solution to a nonlinear equation (finding the solution may require an iterative algorithm), the applicability of AD is not obvious\footnote{It is possible to include the iterative process into the automatic differentiation; however, this approach can be costly in both computation and storage. Additionally, the accuracy is compromised \cite{ablin2020super}.}. 
Fortunately, it is shown that reverse-mode AD and adjoint-state methods are mathematically equivalent \cite{li2019time}; therefore, we can supplement the plain AD with the so-called \textit{physics constrained learning}, proposed in \cite{xu2020physics}, to enable automatic differentiation through implicit functions. The inverse modeling procedure that integrates AD and physics constrained learning is shown in \Cref{fig:routine}. The physics constrained learning is essential to inverse modeling of viscoelastic models because solving the PDE systems usually require an implicit solver for stability and efficiency. 

We conducted extensive numerical results to demonstrate the effectiveness of our method. We divided the numerical examples into two categories according to the governing equations: the geomechanics equation, and the coupled geomechanics and single-/multi-phase flow equations. In both categories, we considered parametric viscoelastic models and neural-network-based models. The major contribution of this work is to develop and implement a unified framework for inverse viscoelasticity computation in the context of different physical models. The code for this work is available at 
\begin{center}
	\url{https://github.com/kailaix/ADCME.jl}
\end{center}

\section{Mathematical Models}

\subsection{Conservation Laws and Constitutive Relations}

In order to give a full description of viscous constitutive models in coupled geomechanics and multi-phase flows, we start with the simpler linear elasticity model. 
In the linear elasticity model, the governing equations are described by two components: the conservation law and the constitutive relation. Along with the boundary and initial conditions, the governing equations describe the deformation of linear elastic materials. 

The conservation law is formulated by the balance of linear momentum: 
\begin{equation}\label{equ:balance}
	\mathrm{div}\bsigma + \rho \mathbf{g} = 0
\end{equation}
where $\bsigma$ is the stress tensor, $\rho$ is the mass density of the material, $\mathbf{g}$ is the gravity vector, and $\mathrm{div}$ is the divergence operator. 

The constitutive relation describes the mapping between the strain and stress tensors:
\begin{equation}\label{equ:cons}
	\bsigma(\bu) = 2\mu\bepsilon(\bu) + \lambda \nabla\cdot \bu \mathbf{I}\qquad \bepsilon(\bu)  =\frac{1}{2}\left(\nabla \bu + (\nabla \bu)^T \right)
\end{equation}  
where $\mathbf{I}$ is the identity tensor, and $\lambda, \mu>0$ are the Lam\'e coefficients, which are related to the Young's modulus $E>0$ and the Poisson's ratio $\nu\in(-1,0.5)$
\begin{equation}\label{equ:cons2}
	\lambda = \frac{\nu}{(1-2\nu)(1+\nu)}E, \quad \mu = \frac{1}{2(1+\nu)}E
\end{equation}
Here, the strain $\bepsilon$ is the Cauchy strain
$$\bepsilon = \frac{1}{2}(\nabla \mathbf{u} + (\nabla \mathbf{u})^T)$$
As for the boundary conditions, we consider two types of boundaries
\begin{itemize}
	\item Traction: for tractions $\mathbf{t}$ imposed on the portion of the surface of the body $\Gamma_N$
	$$\bsigma\bn = \mathbf{t}$$
	Here $\bn$ is a surface normal.	Particularly, when $\mathbf{t}=\mathbf{0}$, we say the boundary condition is \textit{traction-free}. 
	\item Dirichlet boundary conditions: for displacements $\bu_0$ imposed on the portion of the surface of the body $\Gamma_u$
	$$\bu = \bu_0$$
\end{itemize}
\Cref{equ:balance,equ:cons,equ:cons2}, together with the boundary conditions, form the governing equations of linear elasticity materials. In the next subsection, we extend the constitutive relation to viscoelasticity. 

\subsection{Viscoelasticity}

For viscoelastic materials, the strain-stress relationship has the form
\begin{equation}\label{equ:general}
	\bsigma = \bsigma(\bepsilon, \dot\bepsilon)
\end{equation}
Here the dot denotes the time derivative and $\dot\bepsilon$ is called the strain rate. Viscoelastic materials---or in general viscous materials---are rate-dependent, i.e., the strain-stress relation is dependent on the rate at which the strain is developed. 

There are various empirical models of viscoelasticty. Springs and dashpots constitute building blocks for many simplified viscoelasticity models \cite{ozkaya1999mechanical}. For example, the Kelvin-Voight model, the Maxwell model, and the standard linear model are obtained based on connecting springs and dashpots with different configurations \cite{christensen2012theory}. In this work, we generate observation data using the Maxwell model, but our methodologies are applicable to other models as well. The Maxwell model has the following form \cite{banks2011brief}
\begin{equation}\label{equ:maxwell}
	\dot \sigma_{ij} + \frac{\mu}{\eta} \left( \sigma_{ij} - \frac{\sigma_{kk}}{3}\delta_{ij} \right) = 2\mu \dot \epsilon_{ij} + \lambda \dot\epsilon_{kk}\delta_{ij}
\end{equation}

When $\mu, \eta$, and $\lambda$ are independent of $\bepsilon$ and $\bsigma$, \Cref{equ:maxwell} is a first order linear partial differential equation. In numerical simulations, we deploy an implicit scheme to \Cref{equ:maxwell} and obtain the following the semi-discrete equation
\begin{equation}
	 \frac{\sigma_{ij}^{n+1} - \sigma_{ij}^n}{\Delta t} + \frac{\mu}{\eta} \left( \sigma_{ij}^{n+1} - \frac{\sigma_{kk}^{n+1}}{3}\delta_{ij} \right) = 2\mu  \frac{\epsilon_{ij}^{n+1}-\epsilon_{ij}^n}{\Delta t} + \lambda \frac{\epsilon_{kk}^{n+1}-\epsilon_{kk}^n}{\Delta t} \delta_{ij}
\end{equation}
Here $\Delta t$ is the time step size and $n$ is the time index.

For viscoelasticity models and other viscous models (the stress has the general form \Cref{equ:general}), we usually need to calibrate the model parameters empirically in practice. Unfortunately, we usually cannot measure the stress directly \footnote{There have been a few cases where changes in the stress field over time have been measured (somewhat directly), but this is rarely possible.}, but we can measure the displacement data $\bu$ at certain locations with sensors. Additionally, the deformation of the material is usually a result of multi-physics interactions. Therefore, we focus on modeling the viscous relations from limited subsurface displacement data, which are dependent on the coupled system of geomechanics and single-phase flow. 

\subsection{Coupled Geomechanics and Single-Phase Flow}

We first consider single phase porous media flow. The governing equations with unspecified constitutive relations for the single phase flow are \cite{burghardt2017geomechanical}
\begin{equation}\label{equ:single}
	\begin{aligned}
		\mathrm{div}(\bsigma'-bp \mathbf{I}) + \rho \mathbf{g} &= 0\\
		\frac{1}{M}\frac{\partial p}{\partial t} + b\frac{\partial \epsilon_v}{\partial t} + \mathrm{div}\bv & = f
	\end{aligned}
\end{equation}
where $\bsigma'$ is the Biot effective stress, $b$ is the Biot coefficient, $\rho$ is the mass density, $\epsilon_v = \mathrm{tr} \bepsilon$, $f$ is a volumetric source term, and the Darcy's velocity $\bv$ is 
$$\mathbf{v} = -\frac{1}{B_f}\frac{k}{\mu}(\nabla p - \rho \mathbf{g})$$
Here $k$ is the absolute permeability tensor, $\mu$ is the fluid viscosity and $B_f$ is the formation volume factor of the fluid. The mechanical equation and porous media transport (fluid) equation are coupled through $p$ and $\epsilon_v$. 

As a simple example, for linear elastic materials, the constitutive relation is linear, and the system \Cref{equ:single} is reduced to the poroelasticity equation  (assume $\mathbf{g}=\mathbf{0}$) \cite{kolesov2014splitting}
\begin{equation}
	\begin{aligned}
    \mathrm{div}\bsigma'(\bu) - b \nabla p &= 0\\
    \frac{1}{M} \frac{\partial p}{\partial t} + b\frac{\partial \epsilon_v(\bu)}{\partial t} - \nabla\cdot\left(\frac{k}{B_f\mu}\nabla p\right) &= f(x,t)\\
    \bsigma'(\bu) &= 2\mu\bepsilon(\bu) + \lambda \nabla\cdot \bu \mathbf{I}\qquad \bepsilon(\bu)  =\frac{1}{2}\left(\nabla \bu + (\nabla \bu)^T \right)
    \end{aligned}
\end{equation}
with a boundary condition
$$\begin{aligned}
    \bsigma' \bn = 0,\quad x\in \Gamma_{N}^u, \qquad u=0, \quad x\in \Gamma_D^u\\
    -\frac{k}{B_f\mu}\frac{\partial p}{\partial \bn} = \mathbf{0},\quad x\in \Gamma_{N}^p, \qquad p=g, \quad x\in \Gamma_D^p
    \end{aligned}$$

\subsection{Coupled Geomechanics and Two-Phase Flow}\label{sect:twophaseflow}

 We consider the coupling of geomechanics and two-phase flow. Examples of this model includes injecting water into a reservoir to produce oil, or injecting supercritical CO${}_2$ into saline aquifers. The description of the governing equation includes five components. 
 
 The first component is the equation for the solid. We assume that the solid density is constant, and therefore, the mass balance equation of the deforming porous medium is 
\begin{equation}\label{equ:phi}
   \frac{\partial}{\partial t} (1-\phi) + \nabla\cdot(1-\phi)\mathbf{v}_s = 0 \Leftrightarrow \frac{\partial \phi}{\partial t} + \nabla \cdot (\mathbf{v}_s \phi) = \nabla \cdot \mathbf{v}_s
\end{equation}
where $\phi$ is the porosity, and $\bv_s$ is the solid velocity, given by 
$$\mathbf{v}_s = \frac{d\mathbf{u}}{dt}$$
Using \Cref{equ:phi} we obtain the relation between the porosity and the volumetric strain $\epsilon_{v} = \nabla \cdot \mathbf{u}$:
\begin{equation}
	\phi = 1-(1-\phi_0)\exp(-\epsilon_{v})
\end{equation}

The second component describes the equations for the fluid. The mass balance equations of multiphase multicomponent fluid are given by
\begin{equation}
	\frac{\partial }{{\partial t}}(\phi {S_i}{\rho _i}) + \nabla  \cdot ({\rho _i}{\mathbf{v}_{is}}) = {\rho _i}{q_i}, \quad i = 1,2
\end{equation}
Here $q_1$ and $q_2$ are the source functions, $\rho_1$ and $\rho_2$ are the densities, and $S_1$ and $S_2$ are the saturations, which satisfy
$$S_1+S_2=1, \quad S_1\in [0,1], \quad S_2\in [0,1]$$ 

Darcy's velocity of phase $i$ can be represented as \cite{wan2003stabilized}
\begin{equation}\label{equ:vs}
	{\mathbf{v}_{is}} =  - \frac{{K{k_{ri}(S_i)}}}{{{\tilde{\mu}_i}}}(\nabla {P_i} - g{\rho _i}\nabla y), \quad i=1,2
\end{equation}
Here, $K$ is the permeability tensor, but in our case we assume it is a space varying scalar value. $k_{ri}(S_i)$ is a function of $S_i$, and typically the higher the saturation, the easier the corresponding phase is to flow. $\tilde \mu_i$ is the viscosity. $y$ is the depth coordinate, i.e., $\nabla y = [0\;1]'$, $\rho_i$ is the mass density of phase $i$, $\phi$ is the porosity, $q_i$ is the source rate of phase $i$, $P_i$ is the pressure of phase $i$ and $g$ is the gravity constant. We assume the movement of the solid is slow in this study, therefore  Darcy's law is still valid without acceleration terms. $\mathbf{v}_{is}$ is the relative velocity of the phase $i$ with respect to $\mathbf{v}_s$ (also called the interstitial velocity). 

The third component is the equation for the balance of linear momentum. This equation couples the stress tensor in geomechanics and the saturation as well as the pressure in the fluid (assuming the Biot coefficient is 1)
\begin{equation}\label{equ:mechanics}
	\nabla \cdot {\sigma}' -b \nabla \left( S_1P_1 + S_2P_2 \right) + \mathbf{f} = 0
\end{equation}
Here $\sigma'$ should be understood as the effective stress \cite{merxhani2016introduction}, which allows us to treat a multiphase porous medium as a mechanically equivalent single-phase continuum. The external force $\mathbf{f}$ includes the effect of the gravity force. 

The fourth component is the constitutive relation. The constitutive relation connects $\sigma'$ and the displacement $\mathbf{u}$. For example, the linear elastic relation is expressed as 
\begin{equation}\label{equ:linear2}
  \bsigma' = \lambda \mathbf{I}\nabla \cdot \mathbf{u} + 2\mu \bepsilon
\end{equation}

Instead of assuming a linear elasticity model for the geomechanics, we can also model the subsurface solid material by a viscoelasticity model
\begin{equation}\label{equ:maxwell2}
	\dot \sigma_{ij}' + \frac{\mu}{\eta} \left( \sigma_{ij}' - \frac{\sigma_{kk}'}{3}\delta_{ij} \right) = 2\mu \dot \epsilon_{ij} + \lambda \dot\epsilon_{kk}\delta_{ij}
\end{equation}

\section{Inverse Modeling Methodology}

In the last section, we introduced several geomechanics models and coupled geomechanics and single-/multi-phase flow models. Particularly, we introduced the rate-dependent viscous constitutive relations. In this section, we consider how to calibrate the model parameters or train neural networks based constitutive relations from limited displacement data. 

\subsection{Physics Constrained Learning}

The basic idea of calibrating model parameters or neural network using observation data is to formulate the inverse modeling problem as a PDE-constrained optimization problem. In the optimization problem, the objective function measures the discrepancy between observation and prediction. The physical constraints, which are described by PDEs, serve as the optimization constraints. Solving the PDEs (numerically) produces the prediction for the observation data. The free optimization variables are the model parameters or the neural network weights and biases. The optimal values for these parameters are found by minimizing the objective function. In this work, we focus on using a gradient descent method to minimize the objective function; more specifically, throughout the paper we use the L-BFGS-B \cite{byrd1995limited} optimization method.  For L-BFGS-B, we use the line search routine in \cite{more1994line}, which attempts to enforce the Wolfe conditions \cite{byrd1995limited} by a sequence of polynomial interpolations. Note L-BFGS-B is applicable in our case since the data sets are typically small, otherwise we should resort to first order optimization methods or gradient-free methods \cite{luenberger1984linear} to satisfy the limited memory. 

Mathematically, the PDE-constrained optimization problem has the following representation 
\begin{gather}\label{equ:hhh}
    \min_{\bt}\; L_h(u_h) \\
    \mathrm{s.t.}\;\; F_h(\bt, u_h) = 0 \notag
\end{gather}
where $\bt$ is the unknown physical parameters, $F_h$ is the  physical constraints (discretized PDEs and boundary conditions), $u_h$ is the numerical solution to the PDE system $F_h(\bt, u_h) = 0$, and $L_h$ is the loss function that measures the discrepancy between predicted and observed data. In this work, if we observe values of $\bu(\bx)$, i.e., $\{u_i\}_{i\in \mathcal{I}_{\mathrm{obs}}}$, at location $\{\bx_i\}$, we can formulate the loss function with least squares
\begin{equation*}
    L_h(u_h) =  \sum_{i\in \mathcal{I}_{\mathrm{obs}}}\left( u_h(\bx_i) - u_i \right)^2
\end{equation*}

One way to solve \Cref{equ:hhh} is by treating the constraints as a penalty term and solving an unconstrained optimization problem 
\begin{equation}\label{equ:penalty}
    \min_{\bt, u_h}\;\tilde L_{h,\lambda}(\bt, u_h) :=  L_h(u_h) + \lambda\|F_h(\bt, u_h)\|_2^2
\end{equation}
where $\lambda\in (0,\infty)$ is the penalty parameter. 

In terms of implementation, the unconstrained optimization problem does not require solving the PDE, which can be expensive. The gradients $\frac{\partial L_{h,\lambda}(\bt, u_h)}{\partial \bt}$ and $\frac{\partial L_{h,\lambda}(\bt, u_h)}{\partial u_h}$ are computed via automatic differentiation, which is available in many scientific computing and machine learning packages. Upon calculating the gradients, the objective function and the gradients are provided to an optimization algorithm, which seeks the optimal values of the parameter $\bt$. 

Despite a simple and elegant implementation of the penalty method, this method does not eventually satisfy the physical constraints, and thus can lead to a spurious solution that deviates from the true physics. Additionally, \cite{xu2020physics} shows that for stiff physical problems, the penalty method is much more ill-conditioned than the PDE itself for the PDE-constrained optimization problem. The intuition is that the penalty method treats $u_h$ as an independent variable. The set of added independent variables can be very large in time-dependent problems because solutions at each time step are included in the set. Therefore, we apply another method, physics constrained learning, proposed in \cite{xu2020physics}, to find the optimal parameter $\bt$. 

The basic idea of physics constrained learning is to first solve the PDE system and obtain a solution $u_h=G_h(\bt)$, s.t.,
$$F_h(\bt, G_h(\bt))=0$$
The new loss function becomes
\begin{equation}\label{equ:new_loss}
    \tilde L_h(\bt)  = L_h(u_h) = L_h(G_h(\bt))
\end{equation}

The challenge here is how to compute the gradient $\frac{\tilde L_h(\bt)}{\partial\bt}$ in the context of numerical PDE solvers. The reverse-mode automatic differentiation is suitable for explicit functions, such as functions with analytical derivatives, e.g., $\cos$, $\sin$, or compositions of these analytically differentiable functions. However, for a function that defines its output implicitly, which usually requires an iterative algorithm to numerically find the solution, the applicability of reverse-mode automatic differentiation is limited. For example, consider an implicitly defined function $f: x\mapsto y$, where $y$ satisfies $x^3-y^3-y=0$. It is not obvious how to compute $\frac{\partial f}{\partial x}$ using automatic differentiation. The idea of physics constrained learning is to apply the implicit function theorem \cite{krantz2012implicit} to extracts the gradients. For the example above, let $y=y(x)$ depends on $x$, and take the derivative with respect to $x$ on both sides, then we have
$$3x^2 - 3y(x)^2y'(x)-y'(x)=0\Rightarrow y'(x) = \frac{3x^2}{1+3y(x)^2}$$

Applying the same implicit function theorem to \Cref{equ:new_loss}, we obtain \cite{xu2020physics}
\begin{equation}\label{equ:s2}
    \frac{{\partial {{\tilde L}_h}(\bt )}}{{\partial \bt }} 
    = - \frac{{\partial {L_h}({u_h})}}{{\partial {u_h}}} \;
    \Big( {\frac{{\partial {F_h(\bt, u_h)}}}{{\partial {u_h}}}\Big|_{u_h = {G_h}(\bt )}} \Big)^{ - 1} \;
    \frac{{\partial {F_h(\bt, u_h)}}}{{\partial \bt }}\Big|_{u_h = {G_h}(\bt )}
\end{equation}

The key observation in physics constrained learning is that we can leverage the reverse-mode automatic differentiation to compute \Cref{equ:s2} efficiently. The procedure is shown in Algorithm \ref{algo:pcl}. The physics constrained learning can be integrated with the original AD framework and thus cooperates with other AD components to calculate the gradients. 
\begin{algorithm}[htpb]
\KwIn{$u_h$, $\bt$, $\frac{\partial L_h(u_h)}{\partial u_h}$}
\KwOut{$\frac{\tilde L_h(\bt)}{\partial \bt}$}

\nl Compute ${\frac{{\partial {F_h(\bt, u_h)}}}{{\partial {u_h}}}\Big|_{u_h = {G_h}(\bt )}}$, which may already be available from forward computation. 

\nl Solve for $\bx$ the linear system 
$$\Big( {\frac{{\partial {F_h(\bt, u_h)}}}{{\partial {u_h}}}\Big|_{u_h = {G_h}(\bt )}} \Big)^T \bx = \left(\frac{{\partial {L_h}({u_h})}}{{\partial {u_h}}}\right)^T$$

\nl Compute the following expression (a scalar) symbolically (the meaning of ``symbolically'' is clear in the context of automatic differentiation)  
$$g(\bt) = \bx^TF_h(\bt, u_h)$$

\nl Treating $\bx$ independent of $\bt$, compute the gradient using automatic differentiation 
$$\frac{\tilde L_h(\bt)}{\partial \bt} = \frac{\partial g(\bt)}{\partial \bt}$$

\caption{{\bf Physics constrained learning for computing the gradients of implicit functions} \label{algo:pcl}}

\end{algorithm}

\subsection{Neural-Network-based Viscous Constitutive Relations}

In this work we mainly solve two types of inverse modeling problems: 
\begin{enumerate}
	\item In the first type of problems, we calibrate the model parameters in the constitutive relation, assuming the form of the relation is \textit{known}. For example, we want to estimate $\mu$, $\eta$, and $\lambda$ in \Cref{equ:maxwell} from the horizontal displacement data on the top surface of the domain. 
	\item In the second type of problems, the form of the constitutive relation is \textit{unknown}, or only \textit{partially known} (for example, $\eta$ is an unknown function of strain or stress in the numerical example in \Cref{sect:3}, although the viscoelastic model is \Cref{equ:maxwell}). This case is usually known as the function inverse problem, where the unknown is a function instead of a parameter.
\end{enumerate}

The solution to the first case using physics constrained learning is straightforward, where the unknown physical parameters serve as the free optimization variables in \Cref{equ:penalty}. 

In the second case, we propose replacing the unknown function $y=f(x)$ by a neural network $y=\mathcal{NN}_{\bt}(x)$, where $\bt$ is the weights and biases of the neural network. Alternative functional approximators, such as radial basis functions and linear basis functions, are also possible. However, the unknown function is usually a high dimensional mapping and may be nonsmooth. These properties make the application of traditional basis functions very challenging. Once we ``parametrize'' the unknown function using a neural network, the inverse modeling problem is reduced to the first case, where free optimization variables $\bt$ are the weights and bias of the neural network. The same optimization technique as in the first case can be applied. 

\paragraph{Form of the Neural-Network-based Viscous Constitutive Relations} It is important that we design neural-network-based constitutive relations so that certain physical constraints are satisfied in order to generalize the constitutive relations to  broad contexts and new data. When we have no prior knowledge on the form of viscous constitutive relations, we propose to model the viscosity part using a neural network and approximate the constitutive relation using  the incremental form
\begin{equation}\label{equ:form_nn_incremental}
  \Delta \bsigma^{n} = \mathcal{NN}_{\bt}^* (\bsigma^n, \bepsilon^n) + H\Delta \bepsilon^{n}
\end{equation}
where 
$$\Delta \bsigma^{n} = \bsigma^{n+1} - \bsigma^n\quad \Delta \bepsilon^{n} = \bepsilon^{n+1} - \bepsilon^n$$
Here $\bt$ is the weights and biases of the neural network, $H$ is a symmetric positive definite matrix so that $\frac{\partial \bsigma^{n+1}}{\partial \bepsilon^{n+1}} = H\succ 0$, and 
The symmetric positive definiteness of the tangent stiffness matrix  $H$ ensures that the strain energy is weak convex. The weak convexity is beneficial for stabilizing both training and predictive modeling \cite{xu2020learning}. Equivalently, \Cref{equ:form_nn_incremental} can be written in the following form 
\begin{equation}\label{equ:form_nn}
\bsigma^{n+1} - H\bepsilon^{n+1} = \mathcal{NN}_{\bt} (\bsigma^n, \bepsilon^n):= \mathcal{NN}_{\bt}^* (\bsigma^n, \bepsilon^n) + \bsigma^n - H\bepsilon^n
\end{equation}
In \Cref{equ:form_nn}, the current stress   $\bsigma^{n+1}$ depends on the current strain $\bepsilon^{n+1}$, and strain-stress tensors $(\bepsilon^n, \bsigma^n)$ at last time step. Therefore,  \Cref{equ:form_nn} implicitly encodes the strain rate  for constant time steps. We will only consider constant time steps in this work, but the form \Cref{equ:form_nn} can be generalized to varying time steps by interpolating the strain and stress data. 

\paragraph{Training the Neural Network} Besides the form of the neural-network-based viscous constitutive relations, the training method is also very important for stability and generalization. In the case we have access to all strain and stress data, we can train the neural-network-based on the input-output pair $\left((\bsigma^n, \bepsilon^n, \bepsilon^{n+1}), \bsigma^{n+1}\right)$ by minimizing 
\begin{equation}\label{equ:supervised}
  \min_{\bt}\; \sum_n \left\|\bsigma^{n+1} - \left(\mathcal{NN}_{\bt} (\bsigma^n, \bepsilon^n) + H\bepsilon^{n+1}\right)\right\|^2
\end{equation}
A second approach is to generate a sequence of hypothetical stress tensors using the neural-network-based constitutive relation
\begin{equation}\label{equ:constitutive_relation}
  \tilde\bsigma^{n+1} = \mathcal{NN}_{\bt} (\tilde\bsigma^n, \bepsilon^n) + H\bepsilon^{n+1}, \qquad \tilde\bsigma^1 = \bsigma^1
\end{equation}
Thus $\tilde \bsigma^{n}$ is a function of $\bt$.  We train the neural network by minimizing 
\begin{equation}\label{equ:rnn}
  \min_{\bt}\; \sum_n \|\tilde\bsigma^n(\bt) - \bsigma^n\|^2
\end{equation}
The training approach \Cref{equ:supervised} is more computationally efficient because terms in the summation are decoupled. The second approach \Cref{equ:rnn} mimics the structure of recurrent neural network \cite{lipton2015critical}, where $\tilde \bsigma^n$ serves as the hidden states. This form is expensive to train due to the sequential dependence of $\tilde \bsigma^n$ but retains history as contextual information. As an important finding, we demonstrate using experimental data in \Cref{sect:4a} that  training using the first approach \Cref{equ:supervised} leads to a nonphysical and instable solution, but training using \Cref{equ:rnn} leads to a much more generalizable and robust neural-network-based constitutive relation.

 Moreover, in practice, the stress and strain tensors are not available in the training dataset. The lack of input-output training pairs make the first approach inapplicable. The second approach can still be applied, except that it must be coupled with the conservation laws to extract hypothetical displacement data for computing an error functional (loss function). The coupling makes the implementation more challenging but the generalization and robustness in predicting processes makes it an ideal candidate for constitutive modeling.

\subsection{Assessment}

There are many hyperparameters or settings that we can use to control the optimization process. For example, thee hyperparameters we consider in this work include the stop criterion of optimizers, the architectures of neural networks, and the initial guesses. Eventually the estimated parameters (or neural network weights and biases) converge to certain values according to the stop criterion (e.g., the relative error is smaller than a certain threshold). However, the convergence does not guarantee that that solution is optimal; therefore, we need some approaches to assess or validate our methods. 

In this work, to validate our method, we can generate synthetic data with known parameters and use these data to solve the inverse problem to obtain estimated values for those parameters. We can then compare the estimated values and true values. For convenience, we call this approach the \textit{direct method}. 

In some cases, the model used in the synthetic data do not have the same form as our employed model in the inverse problem. In this case, the direct method is inapplicable. We propose the \textit{reproduction method}, which  substitutes both the calibrated and the true models into the PDE model and compare representative quantities. One caveat is that the representative quantities should reflect the influence of models, otherwise we may obtain the same representative quantities no matter how our calibrated models deviate from true ones. 

The last two approaches serve as tools to verify our method on synthetic data. In reality, the true physical law may not have the same form as our employed model in the inverse problem. For example, we use a Kelvin-Voigt model to characterize the viscosity, while the true constitutive relation may be much more complex. In this case, we use the \textit{train-test method}. The idea is as follows. If we have multiple observation data derived from different inputs, e.g., different boundary conditions or source functions, we can split the dataset into training set $T_t$ and testing set $T_e$. The set $T_t$ is used to estimated the optimal parameters $\bt^*$. The calibrated parameters are provided to the PDE solver to generate the prediction $u_h(\bt^*)$, which are compared with the corresponding data in the testing set $T_e$. The train-test sets method follows the general training and testing philosophy in machine learning and serves as a good indicator for the generality of the learned model. The method requires multiple sets of data and implicitly assumes that the training and testing data sets are from the same distribution \cite{lecun2015deep}.

In the numerical examples below, we will use one of the three methods---the direct method, the reproduction method, or the training and testing sets method---that is appropriate for each scenario to assess the quality of our learned constitutive models.

\section{Numerical Experiments}

The numerical examples consist of two parts. In the first part, we only consider the geomechanics, and in the second part, we consider the system of coupled geomechanics and flow equations. We also consider two learning problems: learning the parameters in a physics based constitutive model and training neural-network-based constitutive model. The former can serve as a verification of the algorithm we developed since we know the true parameters. 

\begin{figure}[htpb]
\centering
  \includegraphics[width=0.3\textwidth]{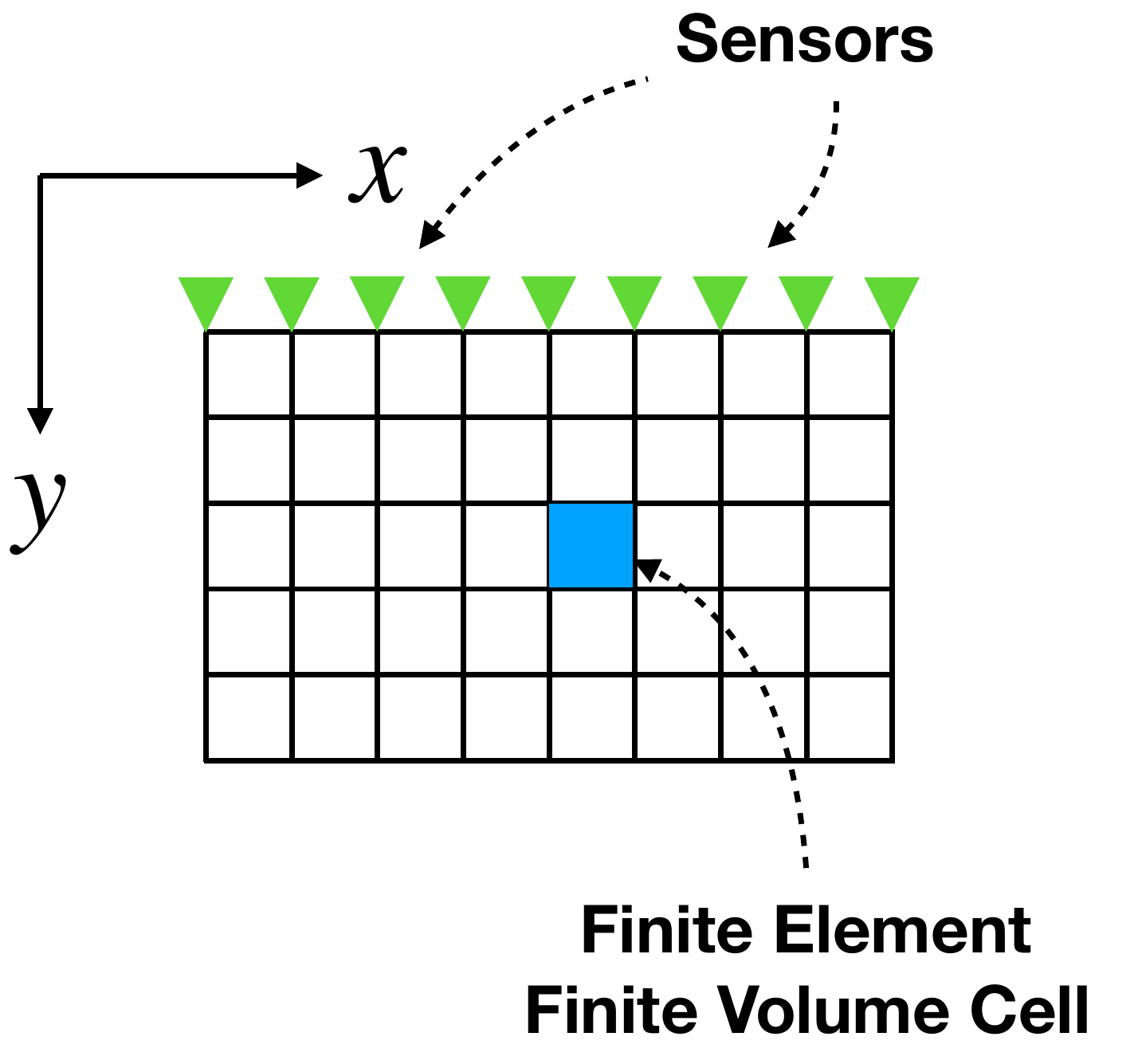}
  \caption{Finite element and finite volume computational domain for the numerical examples. The green triangles are the sensors where the horizontal ($x$ direction) displacement data are measured.}
  \label{fig:cell}
\end{figure}

In the examples below, we assume that the observations are horizontal displacement data on the ground surface, which can be collected by various approaches (e.g., via differential GPS) in practice. For numerical simulation, we use the finite element method for the geomechanics and the finite volume method for the flow equations (see \Cref{fig:cell} for the computational domain). \Cref{tab:numerical_examples} shows a summary of the numerical examples. For visualization of the stress tensor, we use the von Mises stress \cite{lubliner2008plasticity}
$$\sigma_{\mathrm{vm}} = \sqrt{\sigma_{11}^2 - \sigma_{11}\sigma_{22} + \sigma_{22}^2 + 3\sigma_{12}^2}$$
The simulation parameters and codes for reproducing the results in this section are available at 
\begin{center}
	\url{https://github.com/kailaix/PoreFlow.jl}
\end{center}

\begin{table}[htpb]
\centering
\caption{Summary of Numerical Examples}
\label{tab:numerical_examples}
\begin{tabular}{@{}llll@{}}
\toprule
Example & Governing Equation & Unknown Object & Assessment Method \\ \midrule
1 & Geomechanics & Parameter & Direct \\
2 & Coupled Geomechanics and Single-Phase Flow & Parameter & Direct \\
3 & Geomechanics & Function & Reproduction \\
4 & (Experiment) & Function & Train-test \\
5 & Coupled Geomechanics and Single-Phase Flow & Function & Train-test \\
6 & Coupled Geomechanics and Two-Phase Flow & Parameter & Direct \\ \bottomrule
\end{tabular}
\end{table}
\subsection{Learning Space Varying Viscosity Coefficients}\label{sect:1}

\begin{figure}[htpb]
\centering
  \includegraphics[width=0.4\textwidth]{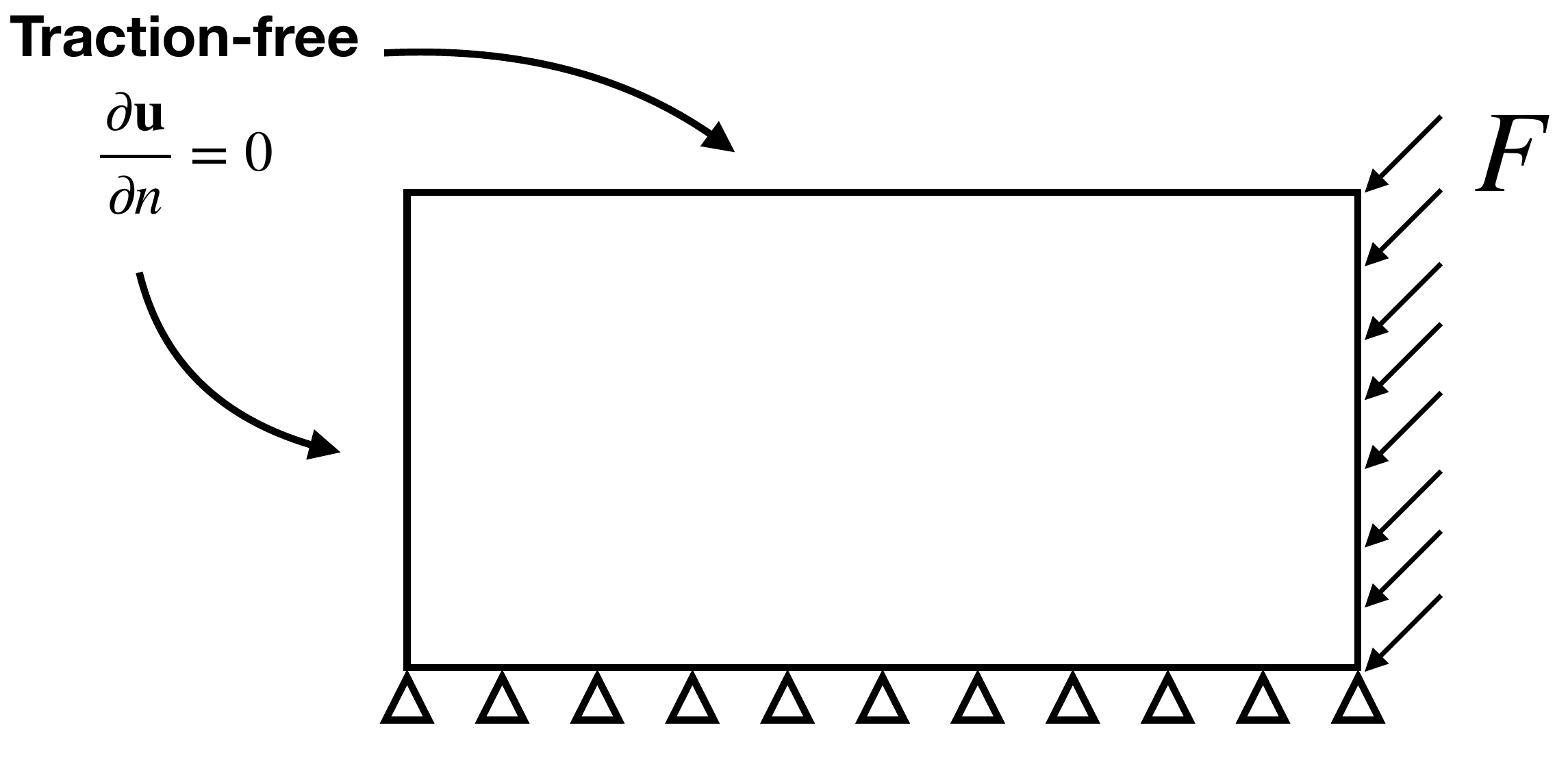}
  \caption{Geometry of inverse modeling problems for \Cref{sect:1,sect:3}.}
  \label{fig:space}
\end{figure}

\Cref{fig:space} shows the geometry of the problem involving only the geomechanics. In this setting, we impose a Dirichlet boundary condition on the bottom side, a traction-free boundary condition on the left and top sides, and external force on the right side. 

In this numerical example, we assume that the region is composed of viscoelastic materials with a viscosity coefficient that varies spatially with depth but is laterally uniform, i.e., $\eta$ is a function of depth $y$. The governing equation is given by 
\begin{equation}\label{equ:1governing}
	\begin{aligned}
			\mathrm{div}\bsigma + \rho \mathbf{g} &= \ddot \bu\\
			\dot \sigma_{ij} + \frac{\mu}{\eta(y)} \left( \sigma_{ij} - \frac{\sigma_{kk}}{3}\delta_{ij} \right) &= 2\mu \dot \epsilon_{ij} + \lambda \dot\epsilon_{kk}\delta_{ij}
	\end{aligned}
\end{equation}
Here $\eta(y)$ depends on the depth $y$. In our algorithm, $\eta(y)$ is discretized on Gauss points, and therefore the optimization problem \Cref{equ:1governing} finds those discrete values by minimizing the discrepancy between observed and predicted displacement data. To solve \Cref{equ:1governing} numerically, we use the $\alpha$-scheme for the simulation of finite element method, which offers accuracy at low-frequency and numerical damping at high frequency \cite{arnold2007convergence,arnold2007convergence,newmark1959method,chung1993time} .

\begin{figure}[htpb]
\centering
  \includegraphics[width=0.33\textwidth]{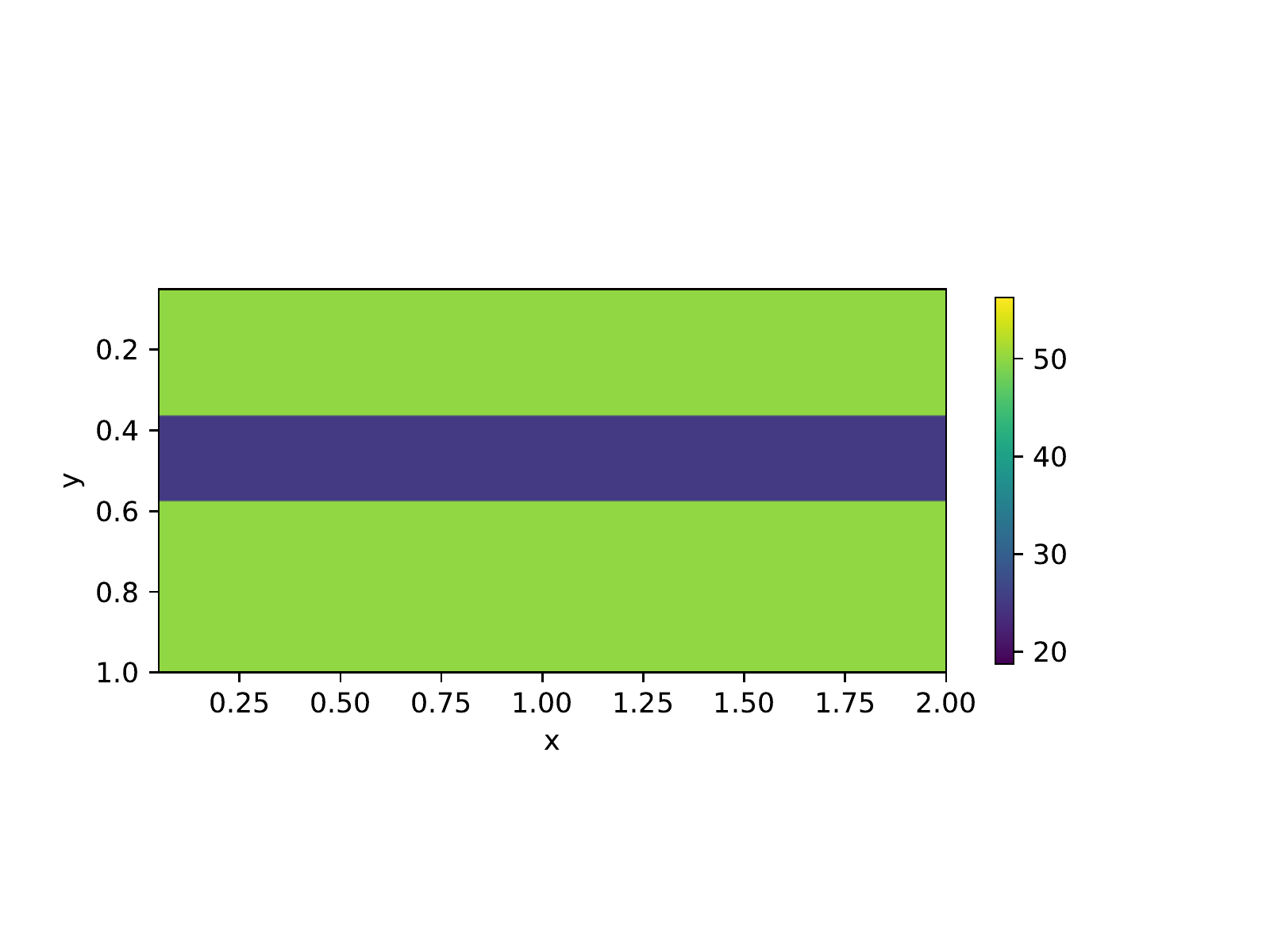}~
  \includegraphics[width=0.33\textwidth]{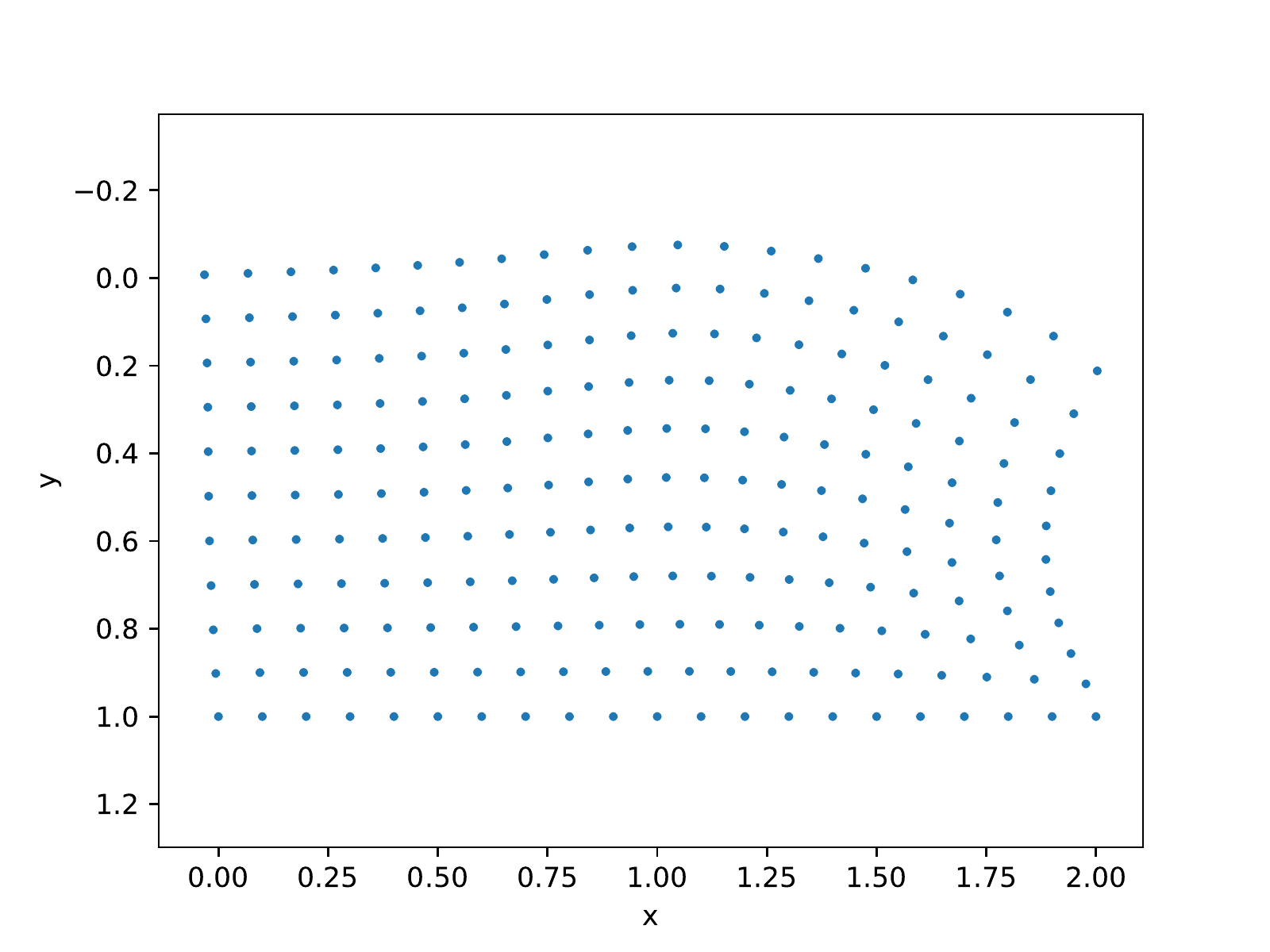}~
  \includegraphics[width=0.33\textwidth]{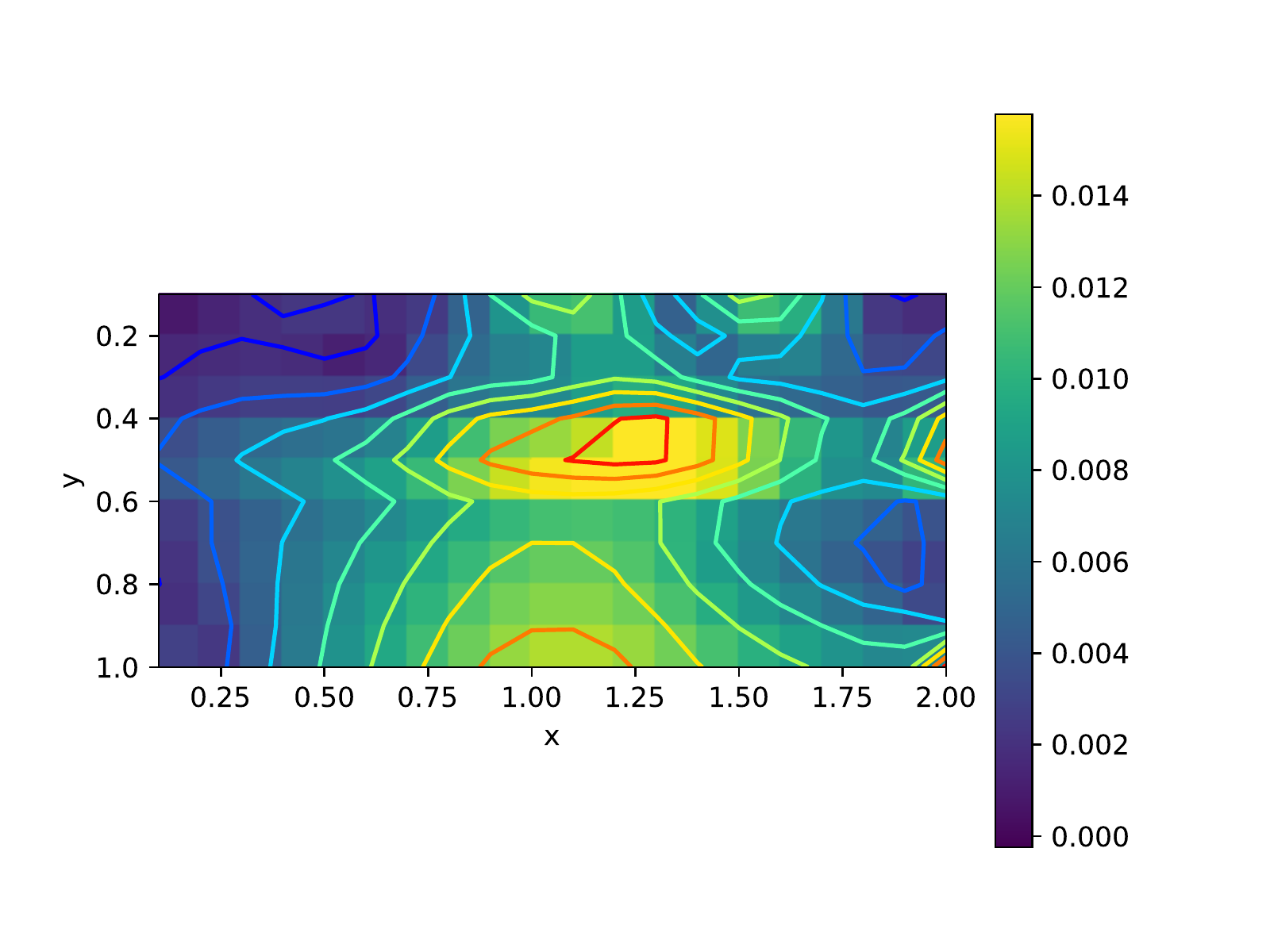}
  \caption{Left: the true model of $\eta(y)^{-1}$. Middle: the displacement at the terminal time. Right: the von Mises stress distribution at the terminal time.}
  \label{fig:spaceus}
\end{figure}

\begin{figure}[htpb]
\centering
  \includegraphics[width=0.245\textwidth]{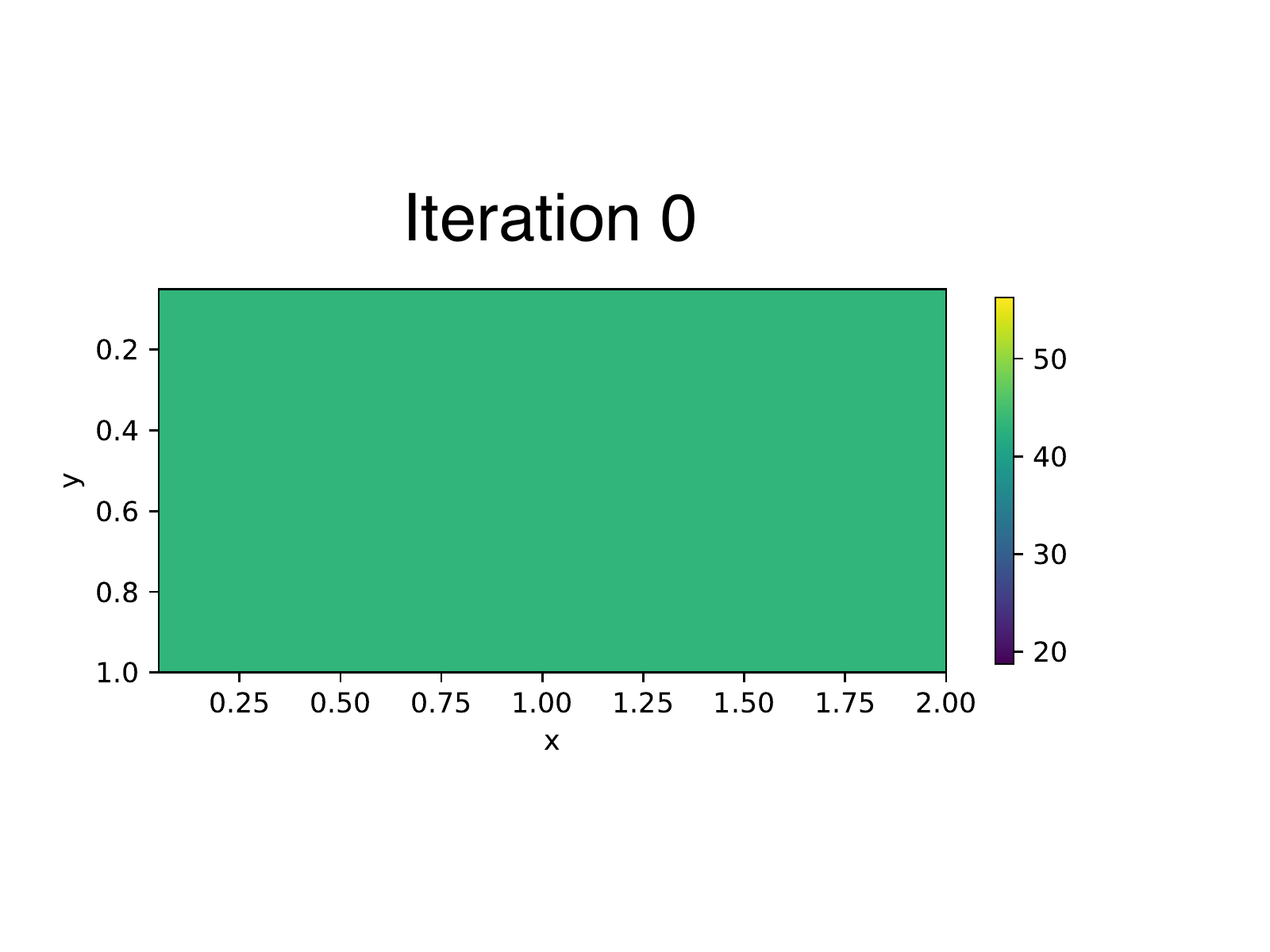}~
  \includegraphics[width=0.245\textwidth]{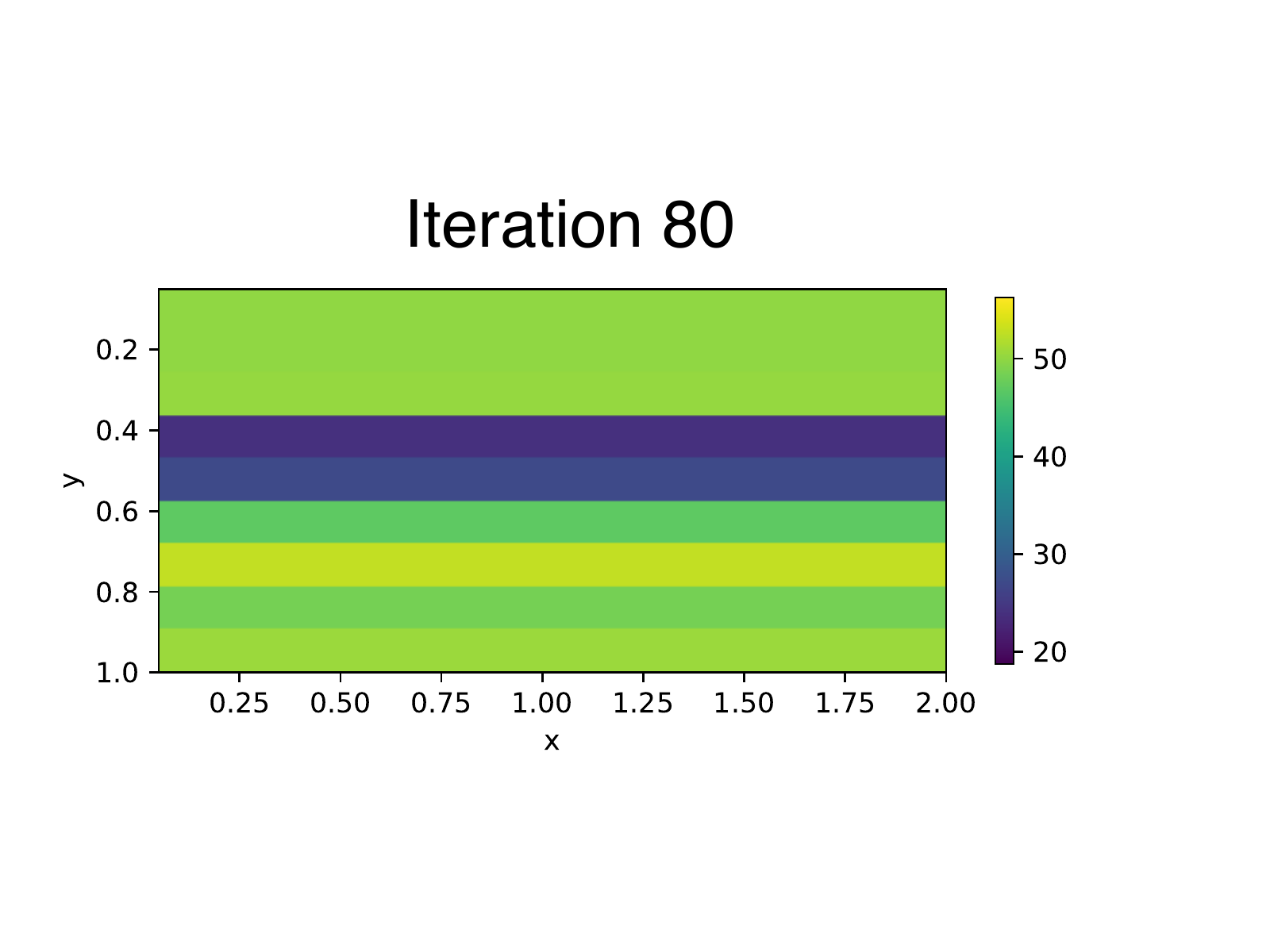}~
\includegraphics[width=0.245\textwidth]{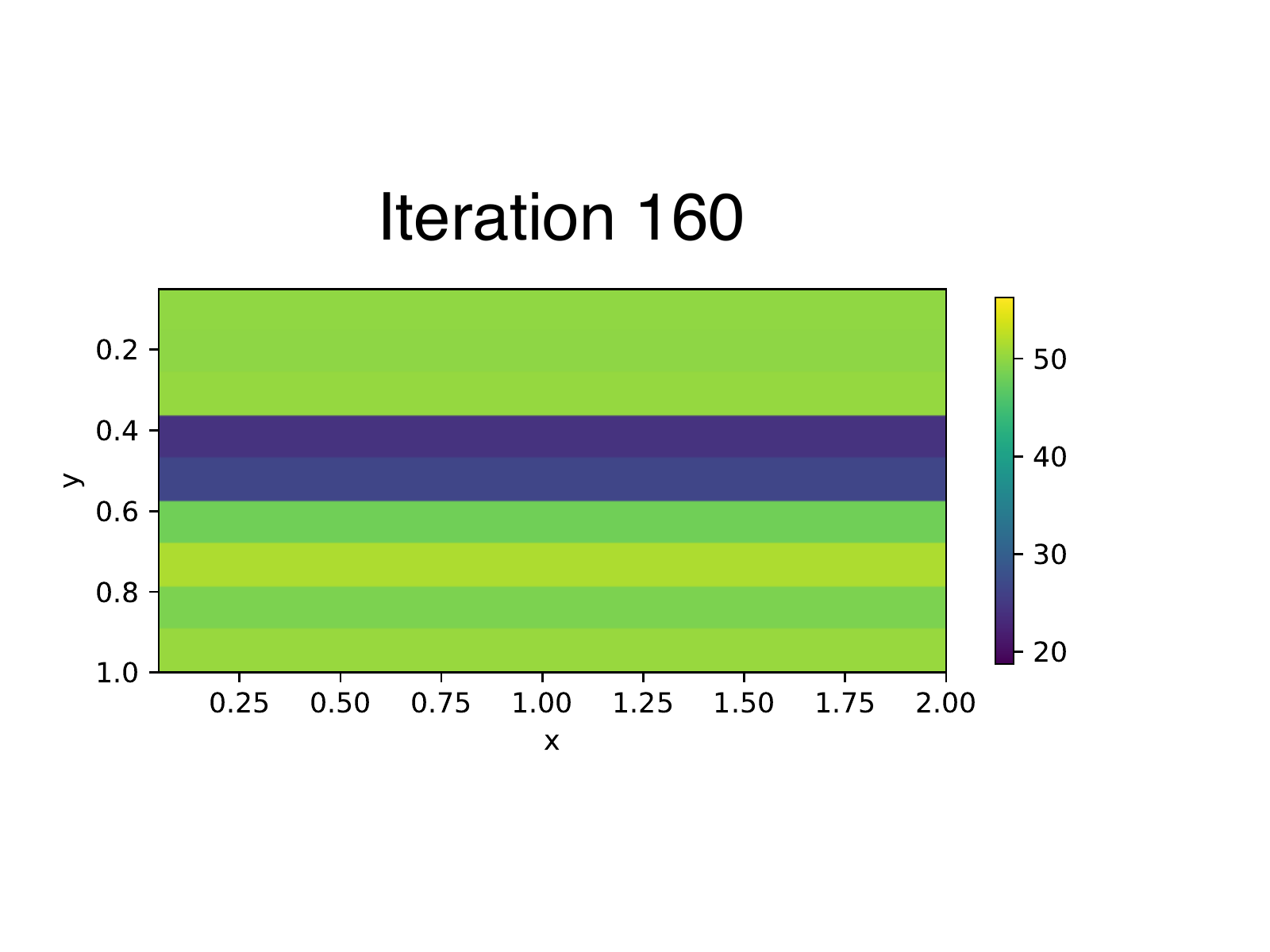}~
\includegraphics[width=0.245\textwidth]{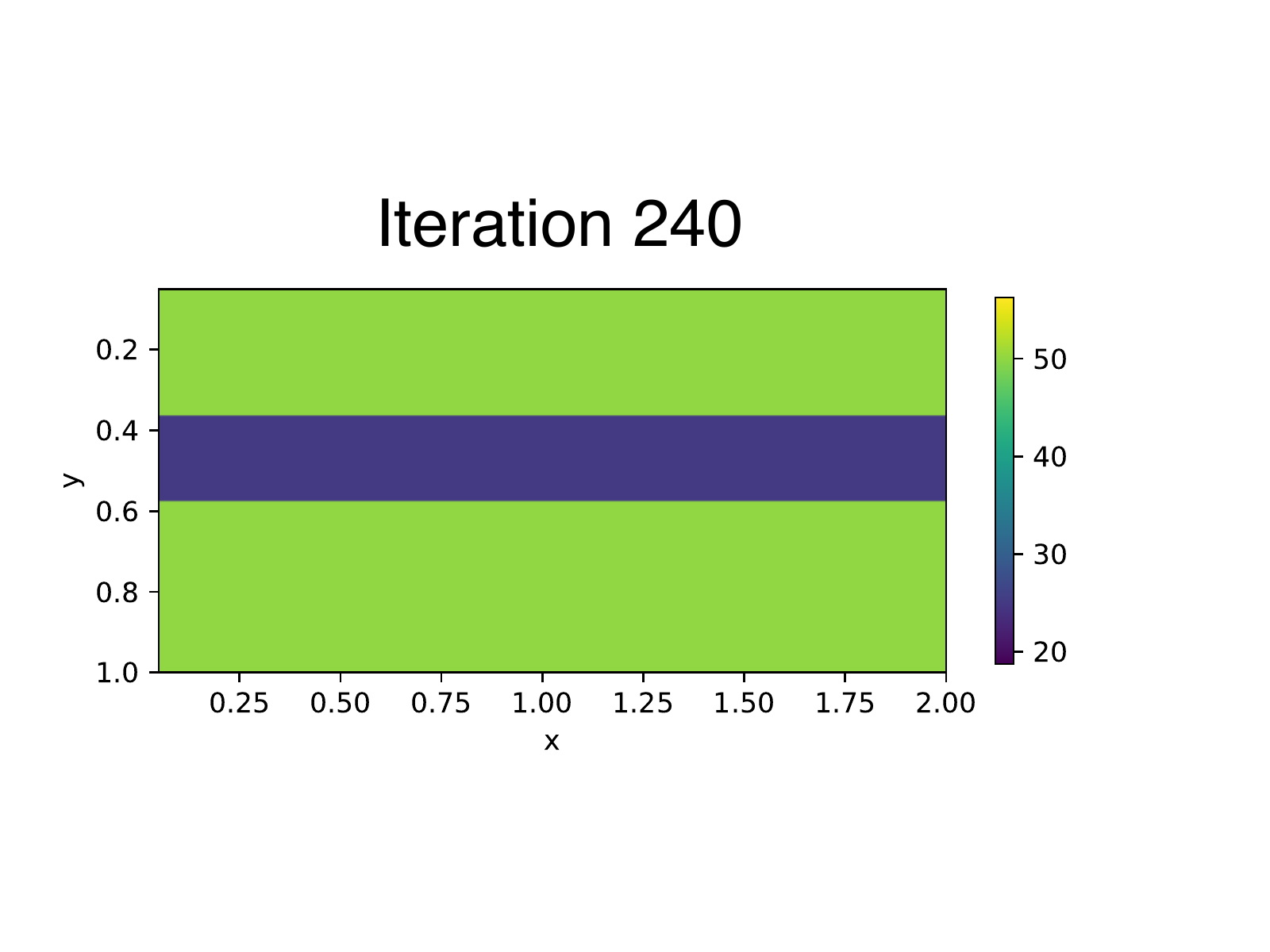}
  \caption{Evolution of learned $\eta(y)^{-1}$ at iteration 0, 80, 160, and 240.}
  \label{fig:spaceiter}
\end{figure}

\Cref{fig:spaceus} shows the true model, the corresponding displacement, and the von Mises stress distribution at the terminal time. \Cref{fig:spaceiter} shows the evolution of estimated viscosity parameter distribution in the optimization process. The results show that we can recover the space varying subsurface viscosity parameters quite accurately with only the horizontal displacement data on the surface.

\subsection{Learning Linear Elasticity Matrix in Poroelasticity Materials}
\label{sect:2}

In this case, we study the estimation of the linear elasticity matrix in the poroelasticity materials. Although the constitutive relation is not related to viscosity, this example serves as a verification of our algorithm in the coupled geomechanics and multi-phase flow. 

\begin{figure}[htpb]
\centering
  \includegraphics[width=0.6\textwidth]{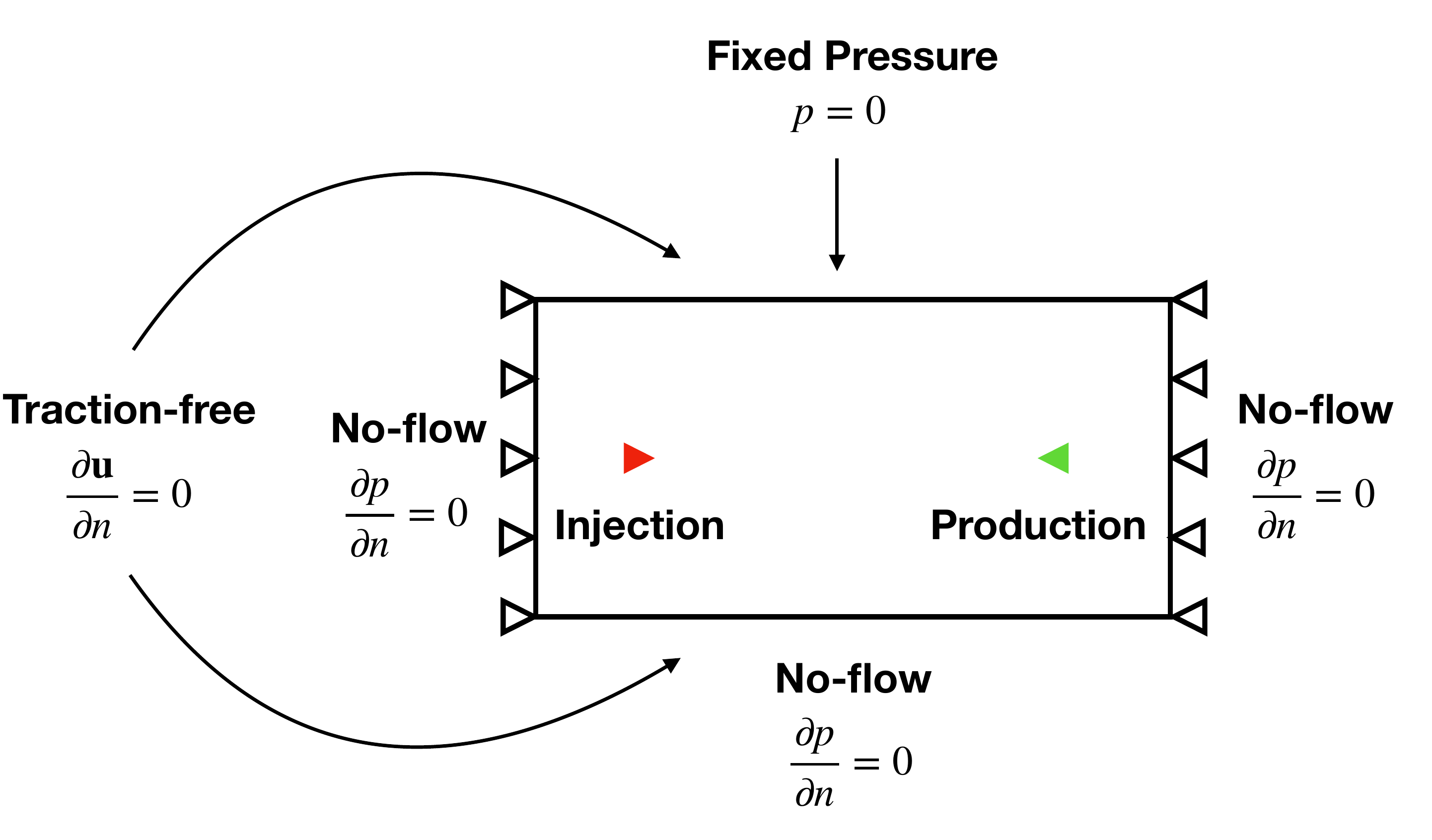}
  \caption{Geometry of inverse modeling problems for \Cref{sect:2,sect:4,sect:5}.}
  \label{fig:ip}
\end{figure}

The geometry of the problem in shown in \Cref{fig:ip}, where we assume no-flow boundary conditions for the pressure on the left, right, and bottom sides; fixed pressure on the top side; traction-free boundary conditions on the top and bottom sides; and zero Dirichlet boundary conditions for the displacement on the left and right sides. The red triangle denotes the injection point, which serves as a positive source, and the green triangle denotes the production point, which serves as a negative source. 

\begin{figure}[htpb]
\centering
  \includegraphics[width=0.3\textwidth]{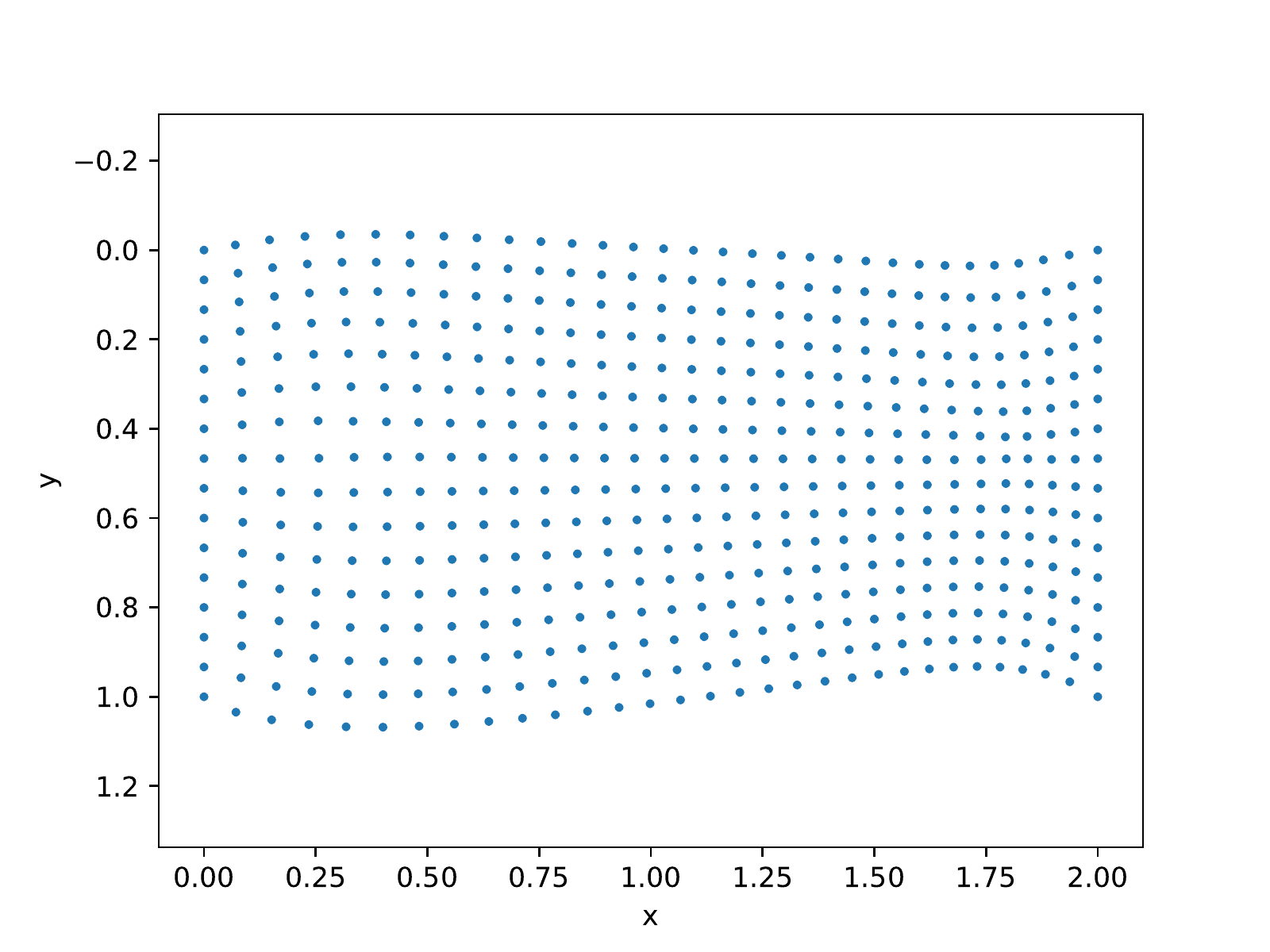}~
  \includegraphics[width=0.33\textwidth]{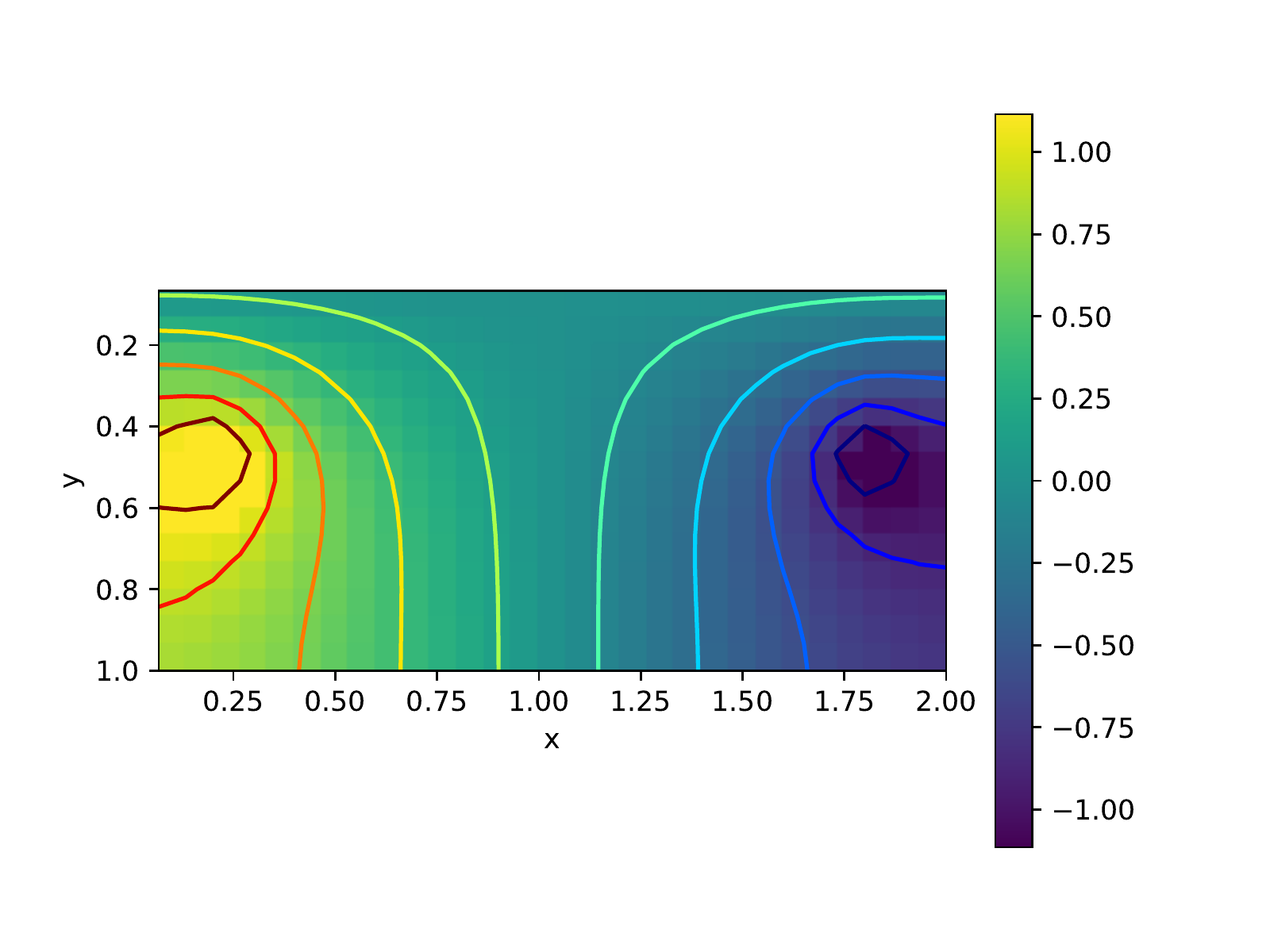}~
\includegraphics[width=0.33\textwidth]{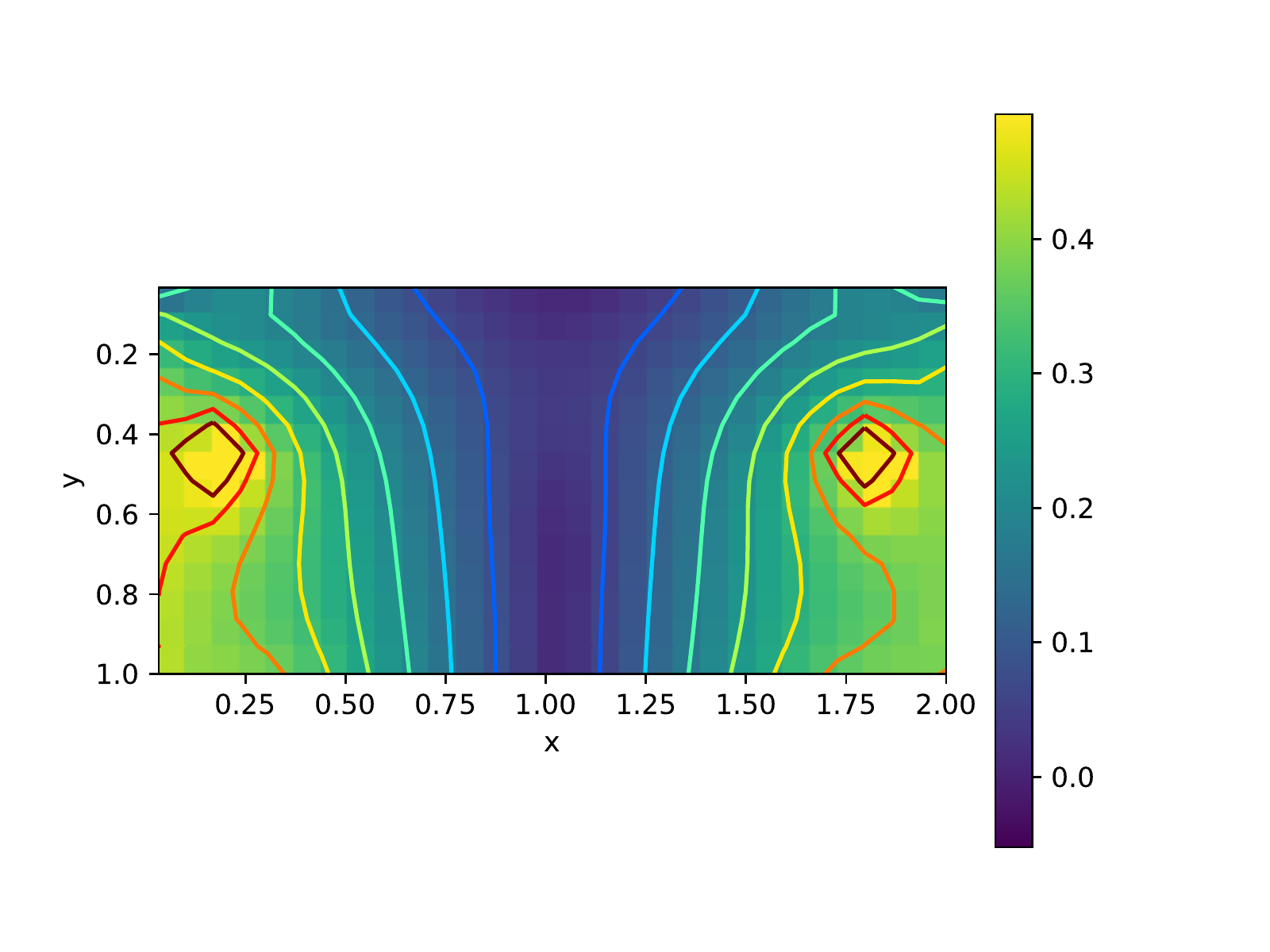}
  \caption{An example of the displacement, pressure and von Mises stress distribution at the terminal time.}
  \label{fig:poro}
\end{figure}

\begin{figure}[htpb]
\centering
	\scalebox{0.6}{
\begin{tikzpicture}

\definecolor{color0}{rgb}{0.12156862745098,0.466666666666667,0.705882352941177}

\begin{axis}[
log basis y={10},
log basis x={10},
tick align=outside,
tick pos=left,
grid=both,
x grid style={white!69.0196078431373!black},
xlabel={\(\displaystyle \sigma_{\mathrm{noise}}\)},
xmin=-0.005, xmax=0.105,
xtick style={color=black},
xtick={0,0.02,0.04,0.06,0.08,0.10},
 xticklabels={0,0.02,0.04,0.06,0.08,0.10},
y grid style={white!69.0196078431373!black},
ylabel={Error},
ymin=4.07600620992903e-07, ymax=0.718877290530786,
ymode=log,
ytick style={color=black}
]
\path [fill=color0, fill opacity=0.5]
(axis cs:0,7.83723763836465e-07)
--(axis cs:0,7.83723763836465e-07)
--(axis cs:0.002,0.00158132130322771)
--(axis cs:0.004,0.00279067352003508)
--(axis cs:0.006,0.00394521102773815)
--(axis cs:0.008,0.000349285027821401)
--(axis cs:0.01,0.00699638757621611)
--(axis cs:0.03,0.00461509316602693)
--(axis cs:0.05,0.0110712847413186)
--(axis cs:0.06,0.0670814259961408)
--(axis cs:0.07,0.0375887667314422)
--(axis cs:0.075,0.040654352874947)
--(axis cs:0.08,0.0560566935916344)
--(axis cs:0.09,0.0331367543693384)
--(axis cs:0.1,0.0260155909185255)
--(axis cs:0.1,0.211137432133811)
--(axis cs:0.1,0.211137432133811)
--(axis cs:0.09,0.373875137591445)
--(axis cs:0.08,0.150975590514042)
--(axis cs:0.075,0.217861989738841)
--(axis cs:0.07,0.0982544842870446)
--(axis cs:0.06,0.186417889821219)
--(axis cs:0.05,0.0467769288751581)
--(axis cs:0.03,0.051192633119898)
--(axis cs:0.01,0.0169072694199763)
--(axis cs:0.008,0.0135885190758564)
--(axis cs:0.006,0.0128751867619852)
--(axis cs:0.004,0.00777170080570544)
--(axis cs:0.002,0.00379291720549701)
--(axis cs:0,7.83723763836465e-07)
--cycle;

\addplot [ultra thick, color0]
table {%
0 7.83723763836465e-07
0.002 0.00268711925436236
0.004 0.00528118716287026
0.006 0.00841019889486169
0.008 0.00696890205183888
0.01 0.0119518284980962
0.03 0.0279038631429625
0.05 0.0289241068082384
0.06 0.12674965790868
0.07 0.0679216255092434
0.075 0.129258171306894
0.08 0.103516142052838
0.09 0.203505945980391
0.1 0.118576511526168
};
\end{axis}

\end{tikzpicture}}
  \caption{Error versus noise levels for poroelasticity models. The shaded area shows the confidence interval for one standard deviation.}
  \label{equ:error}
\end{figure}
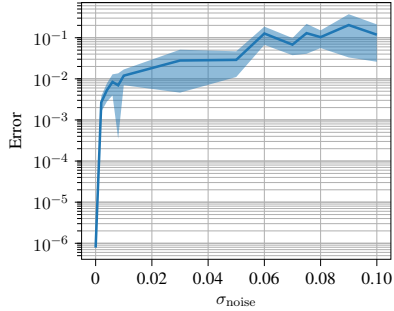

\Cref{fig:poro} shows a plot of the displacement, pressure and von Mises stress distribution at the terminal time. In order to test the robustness of our algorithm, we add Gaussian noise to our data
\begin{equation}\label{equ:noise}
	 u_i = u_i^{\mathrm{true}}(1+\sigma_{\mathrm{noise}}\epsilon_i)
\end{equation}
where $u_i^{\mathrm{true}}$ is the true horizontal displacement at location $i$, and $\epsilon_i$ are i.i.d. standard Gaussian noise, which has a unit variance and zero mean. The error is measured by 
\begin{equation}\label{equ:errorH}
	\mathrm{error} = \frac{\|H_{\mathrm{est}}-H^*\|_2}{\|H^*\|_2}
\end{equation}
\Cref{equ:error} shows a plot of the error versus noise levels. In the context of small noise, the estimation is quite accurate. For example, if the noise equals 0, the estimated and true linear elasticity matrices are
\begin{equation}
	H_{\mathrm{est}} = \begin{bmatrix}
		1.604938 & 0.864198 & 0.0 \\
0.864198 & 1.604938 & 0.0 \\
0.0 & 0.0 & 0.37037 
	\end{bmatrix} \qquad H^* = \begin{bmatrix}
		1.604938 & 0.864197 & 0.0 \\
0.864197 & 1.604938 & 0.0 \\
0.0 & 0.0 & 0.370371 
	\end{bmatrix}
\end{equation}

\subsection{Learning Nonlinear Viscoelasticity Models}
\label{sect:3}

We use the same geometry and governing equations as \Cref{sect:1} for the problem in this section, except that the viscosity parameter $\eta$ is now stress dependent \cite{solomatov1995scaling}, i.e., $\eta$ is a function of $\bsigma$. The two hypothetical viscosity parameters, which we use in benchmarks, are  
\begin{equation}\label{equ:eta12}
  \begin{aligned}
	\eta_1(\bsigma) &= 10 + \frac{5}{1+1000(\sigma_{11}^2 + \sigma_{22}^2)}\\
	\eta_2(\bsigma) &= 10 + \max \left(\frac{50}{1+1000(\sigma_{11}^2 + \sigma_{22}^2)}-10, 0\right)
\end{aligned}
\end{equation}
The coefficients in the equations are deliberately adjusted to make $\eta_1$ and $\eta_2$ effectively vary in the numerical simulations. The second function is continuous but nonsmooth. The nonlinearity in \Cref{equ:eta12} also makes the geomechanics equation highly nonlinear.

The mapping \Cref{equ:eta12} has 3D inputs and 1D output. However, because we assume we do not know the form of the representation, and the relation may be nonsmooth (e.g., $\eta_2$ is continuous but nonsmooth in $\bsigma$), we use a neural network to approximate the mapping between the stress tensor and viscosity parameter, i.e., 
\begin{equation}\label{equ:nlvnn}
	\eta(\bsigma) = \mathcal{NN}_{\bt}(\bsigma)
\end{equation}
where $\bt$ are the weights and biases of the neural network. The neural network consists of 3 hidden layers, 20 neurons per layer, and has tanh activation functions. We run 200 iterations of the L-BFGS-B optimization.

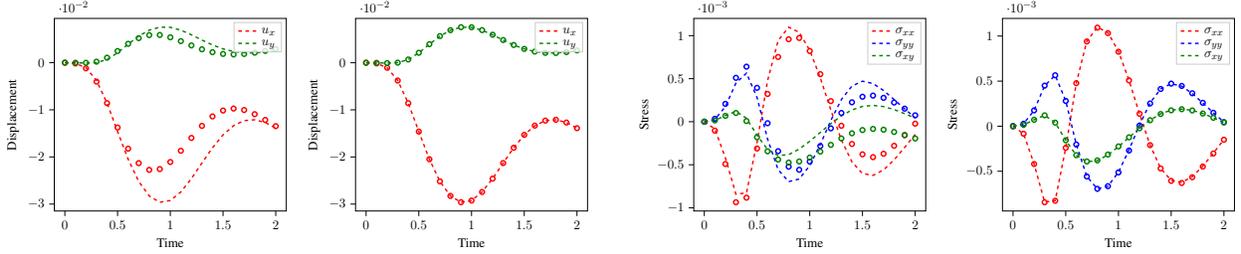
\begin{figure}[htpb]
\begin{subfigure}[t]{0.5\textwidth}
	\scalebox{0.45}{
\begin{tikzpicture}

\begin{axis}[
legend cell align={left},
legend style={fill opacity=0.8, draw opacity=1, text opacity=1, draw=white!80!black},
tick align=outside,
tick pos=left,
x grid style={white!69.0196078431373!black},
xlabel={Time},
xmin=-0.1, xmax=2.1,
xtick style={color=black},
y grid style={white!69.0196078431373!black},
ylabel={Displacement},
ymin=-0.0314658568539842, ymax=0.00939354511585883,
ytick style={color=black}
]
\addplot [very thick, red, dashed]
table {%
0 0
0.1 -0.000141730520011754
0.2 -0.00106641913088168
0.3 -0.0037716143855692
0.4 -0.00858196367520453
0.5 -0.0146070668844201
0.6 -0.0204750054891194
0.7 -0.0251838044409459
0.8 -0.0282714901154662
0.9 -0.0296086113099004
1 -0.0292467192823592
1.1 -0.0274276819475596
1.2 -0.024589480853923
1.3 -0.0212832831905945
1.4 -0.0180496279517417
1.5 -0.0153180198367408
1.6 -0.013361307009435
1.7 -0.0122965755074477
1.8 -0.0121117914566519
1.9 -0.0126981368913402
2 -0.013879690482403
};
\addlegendentry{$u_x$}
\addplot [very thick, red, mark=o, mark size=2, mark options={solid}, only marks, forget plot]
table {%
0 0
0.1 -0.000167363411391404
0.2 -0.00119956550033799
0.3 -0.00401140515749117
0.4 -0.00858949998802152
0.5 -0.0137620344222849
0.6 -0.0182421436179358
0.7 -0.0213127556135915
0.8 -0.0227402383281679
0.9 -0.022572506350587
1 -0.0210724154490806
1.1 -0.0186894597260115
1.2 -0.0159647203263289
1.3 -0.0134070667075268
1.4 -0.0113990340596459
1.5 -0.0101558679816431
1.6 -0.00973110841432128
1.7 -0.0100456242797907
1.8 -0.0109263994192349
1.9 -0.0121468230952621
2 -0.0134666884635271
};
\addplot [very thick, green!50!black, dashed]
table {%
0 0
0.1 -2.03490437097345e-05
0.2 -7.7099446892608e-05
0.3 1.73742276389046e-05
0.4 0.00068440467475761
0.5 0.0021027078004714
0.6 0.00392583288586885
0.7 0.00564068195519395
0.8 0.00690154155818607
0.9 0.00753629957177506
1 0.00750541295402453
1.1 0.00690015362926513
1.2 0.00591911716845484
1.3 0.00480270712425332
1.4 0.00376099740767201
1.5 0.00293137268910026
1.6 0.00237470390941913
1.7 0.00209627074903518
1.8 0.00207044545672343
1.9 0.00225754947078686
2 0.00261004157422391
};
\addlegendentry{$u_y$}
\addplot [very thick, green!50!black, mark=o, mark size=2, mark options={solid}, only marks, forget plot]
table {%
0 0
0.1 -1.36944858993409e-05
0.2 -1.90545175320099e-05
0.3 0.000219969628909723
0.4 0.00105739361275581
0.5 0.00247449274166008
0.6 0.00401001484284375
0.7 0.00521222535038073
0.8 0.00585310237755491
0.9 0.00588009511199154
1 0.00538785485959328
1.1 0.0045740866446039
1.2 0.00366365852460334
1.3 0.00284213125851402
1.4 0.00222485368526125
1.5 0.00185988400693438
1.6 0.00174632372840213
1.7 0.00185188269168545
1.8 0.00212403947777316
1.9 0.00249767464396619
2 0.00290273902098275
};
\end{axis}

\end{tikzpicture}}~
	\scalebox{0.45}{
\begin{tikzpicture}

\begin{axis}[
legend cell align={left},
legend style={fill opacity=0.8, draw opacity=1, text opacity=1, draw=white!80!black},
tick align=outside,
tick pos=left,
x grid style={white!69.0196078431373!black},
xlabel={Time},
xmin=-0.1, xmax=2.1,
xtick style={color=black},
y grid style={white!69.0196078431373!black},
ylabel={Displacement},
ymin=-0.03146656159152, ymax=0.00940834460411116,
ytick style={color=black}
]
\addplot [very thick, red, dashed]
table {%
0 0
0.1 -0.000141730520011754
0.2 -0.00106641913088168
0.3 -0.0037716143855692
0.4 -0.00858196367520453
0.5 -0.0146070668844201
0.6 -0.0204750054891194
0.7 -0.0251838044409459
0.8 -0.0282714901154662
0.9 -0.0296086113099004
1 -0.0292467192823592
1.1 -0.0274276819475596
1.2 -0.024589480853923
1.3 -0.0212832831905945
1.4 -0.0180496279517417
1.5 -0.0153180198367408
1.6 -0.013361307009435
1.7 -0.0122965755074477
1.8 -0.0121117914566519
1.9 -0.0126981368913402
2 -0.013879690482403
};
\addlegendentry{$u_x$}
\addplot [very thick, red, mark=o, mark size=2, mark options={solid}, only marks, forget plot]
table {%
0 0
0.1 -0.000141810432866871
0.2 -0.00106738848710186
0.3 -0.00377538630392223
0.4 -0.00858848947582795
0.5 -0.014611936760145
0.6 -0.0204745522827705
0.7 -0.0251793054729721
0.8 -0.0282657143010674
0.9 -0.0296038349086075
1 -0.0292463307935556
1.1 -0.0274337946197137
1.2 -0.0245991359749842
1.3 -0.021290366640806
1.4 -0.0180498408088453
1.5 -0.0153110885322544
1.6 -0.013350024570196
1.7 -0.0122850115989088
1.8 -0.0121039496976188
1.9 -0.0126971275534313
2 -0.0138870914363948
};
\addplot [very thick, green!50!black, dashed]
table {%
0 0
0.1 -2.03490437097345e-05
0.2 -7.7099446892608e-05
0.3 1.73742276389046e-05
0.4 0.00068440467475761
0.5 0.0021027078004714
0.6 0.00392583288586885
0.7 0.00564068195519395
0.8 0.00690154155818607
0.9 0.00753629957177506
1 0.00750541295402453
1.1 0.00690015362926513
1.2 0.00591911716845484
1.3 0.00480270712425332
1.4 0.00376099740767201
1.5 0.00293137268910026
1.6 0.00237470390941913
1.7 0.00209627074903518
1.8 0.00207044545672343
1.9 0.00225754947078686
2 0.00261004157422391
};
\addlegendentry{$u_y$}
\addplot [very thick, green!50!black, mark=o, mark size=2, mark options={solid}, only marks, forget plot]
table {%
0 0
0.1 -2.03321842723603e-05
0.2 -7.6884947567256e-05
0.3 1.88324837972465e-05
0.4 0.000689529253850607
0.5 0.00211298056965709
0.6 0.00393862280679008
0.7 0.0056519575164283
0.8 0.00691186318727435
0.9 0.00755039432249157
1 0.0075264953139251
1.1 0.00692668884791386
1.2 0.00594565623750318
1.3 0.00482312840415688
1.4 0.00377193780902561
1.5 0.00293344341015117
1.6 0.00237127336689559
1.7 0.00209148383877992
1.8 0.00206779360011903
1.9 0.00225902550579793
2 0.00261565208462884
};
\end{axis}

\end{tikzpicture}}
	\caption{Displacement. Left panel: initial guess; right panel: computed using the calibrated model.}
\end{subfigure}~
	\begin{subfigure}[t]{0.5\textwidth}
	\scalebox{0.45}{
\begin{tikzpicture}

\begin{axis}[
legend cell align={left},
legend style={fill opacity=0.8, draw opacity=1, text opacity=1, draw=white!80!black},
tick align=outside,
tick pos=left,
x grid style={white!69.0196078431373!black},
xlabel={Time},
xmin=-0.1, xmax=2.1,
xtick style={color=black},
y grid style={white!69.0196078431373!black},
ylabel={Stress},
ymin=-0.00103726975574867, ymax=0.00120027380312869,
ytick style={color=black}
]
\addplot [very thick, red, dashed]
table {%
0 0
0.1 -8.40761848897833e-05
0.2 -0.000419058553278349
0.3 -0.000842270468849947
0.4 -0.000826472091143151
0.5 -0.000241217550663016
0.6 0.000476766715294245
0.7 0.00094289724672934
0.8 0.00109856727772517
0.9 0.00103842425733627
1 0.000830043433986845
1.1 0.00050940700092794
1.2 0.000136141819806684
1.3 -0.00021137828616803
1.4 -0.0004670950015521
1.5 -0.00060213566235565
1.6 -0.000623637442963712
1.7 -0.000559760516947456
1.8 -0.000442687531381268
1.9 -0.000298933860781353
2 -0.000147583283219958
};
\addlegendentry{$\sigma_{xx}$}
\addplot [very thick, blue, dashed]
table {%
0 0
0.1 2.72497271017199e-05
0.2 0.000175339605586945
0.3 0.000449032005882322
0.4 0.000564554824513986
0.5 0.00028150937740707
0.6 -0.000201300271825239
0.7 -0.000557234665713862
0.8 -0.000696009345329593
0.9 -0.000665624468356435
1 -0.000509557949466447
1.1 -0.000263603727174275
1.2 1.49851900222686e-05
1.3 0.000257681273158202
1.4 0.000414665378039723
1.5 0.000471067678619722
1.6 0.000443496288503347
1.7 0.000362530908266973
1.8 0.000257556418942402
1.9 0.000148532432209162
2 4.62868040208091e-05
};
\addlegendentry{$\sigma_{yy}$}
\addplot [very thick, green!50!black, dashed]
table {%
0 0
0.1 1.46931085948668e-05
0.2 7.17651876730061e-05
0.3 0.000120723626259471
0.4 3.9738904497661e-05
0.5 -0.000156210017506617
0.6 -0.000322994624962358
0.7 -0.000393736606591584
0.8 -0.0003825495521468
0.9 -0.00031776889630459
1 -0.000227019316843051
1.1 -0.000124600816260352
1.2 -2.23505336478372e-05
1.3 6.80987065788389e-05
1.4 0.000136969044083033
1.5 0.000177965233035422
1.6 0.000189928875550257
1.7 0.000175613657600724
1.8 0.000141082118485671
1.9 9.37319349514393e-05
2 4.15593580389339e-05
};
\addlegendentry{$\sigma_{xy}$}
\addplot [very thick, red, mark=o, mark size=2, mark options={solid}, only marks, forget plot]
table {%
0 0
0.1 -0.000103434992456736
0.2 -0.000491352056497726
0.3 -0.000935563230345154
0.4 -0.000882400103956538
0.5 -0.000311234439276193
0.6 0.000322388714817689
0.7 0.00075230290968379
0.8 0.000958378035374612
0.9 0.000973123639195305
1 0.000821933445631216
1.1 0.000552504449826508
1.2 0.000239452395526682
1.3 -4.67432944031111e-05
1.4 -0.00025844611802368
1.5 -0.000378870168727656
1.6 -0.000411745494806743
1.7 -0.000372280619361181
1.8 -0.000279888931580777
1.9 -0.000155980816660203
2 -2.25052182106596e-05
};
\addplot [very thick, blue, mark=o, mark size=2, mark options={solid}, only marks, forget plot]
table {%
0 0
0.1 3.37617510153115e-05
0.2 0.00020710988289711
0.3 0.000510252682485677
0.4 0.000639699249147533
0.5 0.000393587736823366
0.6 -1.73594164667594e-05
0.7 -0.000342318203768923
0.8 -0.000522717859641346
0.9 -0.000558451756535726
1 -0.000463612728654029
1.1 -0.000282779511965736
1.2 -7.70192477102099e-05
1.3 0.000101664079574036
1.4 0.000225586606152861
1.5 0.000290118325250275
1.6 0.000303541536077736
1.7 0.00027741442191916
1.8 0.000223158018902559
1.9 0.000151615338129585
2 7.43077665463338e-05
};
\addplot [very thick, green!50!black, mark=o, mark size=2, mark options={solid}, only marks, forget plot]
table {%
0 0
0.1 1.49681223931187e-05
0.2 6.79048748520533e-05
0.3 0.000101213793041183
0.4 7.83325316536934e-06
0.5 -0.000180432893932385
0.6 -0.000343506657271889
0.7 -0.000439450955106888
0.8 -0.000474669088392892
0.9 -0.00046185110132337
1 -0.000416058847470978
1.1 -0.000347362552735828
1.2 -0.000269134924784611
1.3 -0.000194366167155962
1.4 -0.000134524508227465
1.5 -9.68716287649335e-05
1.6 -8.36295762368053e-05
1.7 -9.2821209051166e-05
1.8 -0.00011906675230299
1.9 -0.000155510745416679
2 -0.000194771717668091
};
\end{axis}

\end{tikzpicture}}~
	\scalebox{0.45}{
\begin{tikzpicture}

\begin{axis}[
legend cell align={left},
legend style={fill opacity=0.8, draw opacity=1, text opacity=1, draw=white!80!black},
tick align=outside,
tick pos=left,
x grid style={white!69.0196078431373!black},
xlabel={Time},
xmin=-0.1, xmax=2.1,
xtick style={color=black},
y grid style={white!69.0196078431373!black},
ylabel={Stress},
ymin=-0.00094084883263741, ymax=0.00119568233059958,
ytick style={color=black}
]
\addplot [very thick, red, dashed]
table {%
0 0
0.1 -8.40761848897833e-05
0.2 -0.000419058553278349
0.3 -0.000842270468849947
0.4 -0.000826472091143151
0.5 -0.000241217550663016
0.6 0.000476766715294245
0.7 0.00094289724672934
0.8 0.00109856727772517
0.9 0.00103842425733627
1 0.000830043433986845
1.1 0.00050940700092794
1.2 0.000136141819806684
1.3 -0.00021137828616803
1.4 -0.0004670950015521
1.5 -0.00060213566235565
1.6 -0.000623637442963712
1.7 -0.000559760516947456
1.8 -0.000442687531381268
1.9 -0.000298933860781353
2 -0.000147583283219958
};
\addlegendentry{$\sigma_{xx}$}
\addplot [very thick, blue, dashed]
table {%
0 0
0.1 2.72497271017199e-05
0.2 0.000175339605586945
0.3 0.000449032005882322
0.4 0.000564554824513986
0.5 0.00028150937740707
0.6 -0.000201300271825239
0.7 -0.000557234665713862
0.8 -0.000696009345329593
0.9 -0.000665624468356435
1 -0.000509557949466447
1.1 -0.000263603727174275
1.2 1.49851900222686e-05
1.3 0.000257681273158202
1.4 0.000414665378039723
1.5 0.000471067678619722
1.6 0.000443496288503347
1.7 0.000362530908266973
1.8 0.000257556418942402
1.9 0.000148532432209162
2 4.62868040208091e-05
};
\addlegendentry{$\sigma_{yy}$}
\addplot [very thick, green!50!black, dashed]
table {%
0 0
0.1 1.46931085948668e-05
0.2 7.17651876730061e-05
0.3 0.000120723626259471
0.4 3.9738904497661e-05
0.5 -0.000156210017506617
0.6 -0.000322994624962358
0.7 -0.000393736606591584
0.8 -0.0003825495521468
0.9 -0.00031776889630459
1 -0.000227019316843051
1.1 -0.000124600816260352
1.2 -2.23505336478372e-05
1.3 6.80987065788389e-05
1.4 0.000136969044083033
1.5 0.000177965233035422
1.6 0.000189928875550257
1.7 0.000175613657600724
1.8 0.000141082118485671
1.9 9.37319349514393e-05
2 4.15593580389339e-05
};
\addlegendentry{$\sigma_{xy}$}
\addplot [very thick, red, mark=o, mark size=2, mark options={solid}, only marks, forget plot]
table {%
0 0
0.1 -8.41354397685579e-05
0.2 -0.000419645610831092
0.3 -0.000843733779763001
0.4 -0.000827053284776356
0.5 -0.000239878642274871
0.6 0.000477167870865148
0.7 0.000940049855462689
0.8 0.00109342968219809
0.9 0.00103286320755868
1 0.000825886071147398
1.1 0.000508402393298595
1.2 0.00013770982702905
1.3 -0.000210278166009294
1.4 -0.000468857409834917
1.5 -0.000606844401553758
1.6 -0.000629894327466655
1.7 -0.000566063813653795
1.8 -0.000448060819208302
1.9 -0.000302865017165221
2 -0.000149806598173929
};
\addplot [very thick, blue, mark=o, mark size=2, mark options={solid}, only marks, forget plot]
table {%
0 0
0.1 2.72705525864166e-05
0.2 0.000175597539867412
0.3 0.000449753734678851
0.4 0.000564893038868961
0.5 0.000280291304195932
0.6 -0.000203329563882881
0.7 -0.000557916661813833
0.8 -0.000695096945735194
0.9 -0.000665613862472461
1 -0.000512778556695097
1.1 -0.000270180717865383
1.2 7.26736107209964e-06
1.3 0.000252005919805096
1.4 0.000412790537913518
1.5 0.000472487792806979
1.6 0.00044644079394643
1.7 0.000365344506573513
1.8 0.000259279517395167
1.9 0.000148753758780717
2 4.49432499817335e-05
};
\addplot [very thick, green!50!black, mark=o, mark size=2, mark options={solid}, only marks, forget plot]
table {%
0 0
0.1 1.46955806572018e-05
0.2 7.1745619539452e-05
0.3 0.00012049037880991
0.4 3.92047890927062e-05
0.5 -0.000156358442458028
0.6 -0.000321891949284309
0.7 -0.000391268640974457
0.8 -0.000379227234099848
0.9 -0.000314353530689909
1 -0.000224414304194834
1.1 -0.000123604761146277
1.2 -2.30519625807074e-05
1.3 6.63057924313557e-05
1.4 0.000134811382853835
1.5 0.000175966203634256
1.6 0.000188357801991538
1.7 0.000174547031954593
1.8 0.000140467264918192
1.9 9.34373552321627e-05
2 4.13615636222548e-05
};
\end{axis}

\end{tikzpicture}}
	\caption{Stress. Left panel: initial guess; right panel: computed using the calibrated model.}
\end{subfigure}
\caption{The initial and terminal quantities of the upper left point for the viscosity model $\eta_1$.  The dashed lines are true values, and the dots are reproduced values using the calibrated models.}
\label{fig:res1}
\end{figure}

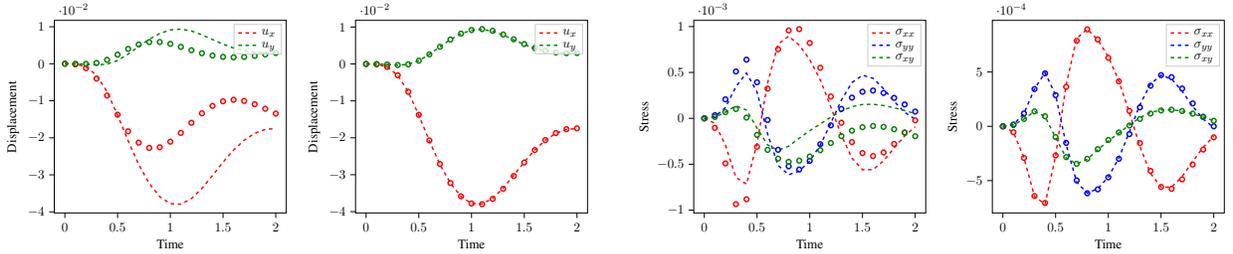
\begin{figure}[htpb]
\begin{subfigure}[t]{0.5\textwidth}
	\scalebox{0.45}{
\begin{tikzpicture}

\begin{axis}[
legend cell align={left},
legend style={fill opacity=0.8, draw opacity=1, text opacity=1, draw=white!80!black},
tick align=outside,
tick pos=left,
x grid style={white!69.0196078431373!black},
xlabel={Time},
xmin=-0.1, xmax=2.1,
xtick style={color=black},
y grid style={white!69.0196078431373!black},
ylabel={Displacement},
ymin=-0.04028679554023, ymax=0.0117202572478728,
ytick style={color=black}
]
\addplot [very thick, red, dashed]
table {%
0 0
0.1 -9.95948177808011e-05
0.2 -0.000803396476905668
0.3 -0.00306748651266507
0.4 -0.0075564188569748
0.5 -0.0138778859365614
0.6 -0.0208231157112783
0.7 -0.0272005793083875
0.8 -0.032321408498514
0.9 -0.0359043525358027
1 -0.0377963308769123
1.1 -0.0379228385953162
1.2 -0.0364147269561444
1.3 -0.0336549484994603
1.4 -0.030188839016599
1.5 -0.0265823093670443
1.6 -0.0233027631516495
1.7 -0.020657060059876
1.8 -0.01879266248127
1.9 -0.0177406878815139
2 -0.0174661154135103
};
\addlegendentry{$u_x$}
\addplot [very thick, red, mark=o, mark size=2, mark options={solid}, only marks, forget plot]
table {%
0 0
0.1 -0.000167363411391404
0.2 -0.0011995564929555
0.3 -0.00401122028639252
0.4 -0.00858832047600062
0.5 -0.0137580919360129
0.6 -0.0182335546392594
0.7 -0.0212987785815573
0.8 -0.0227215051003115
0.9 -0.0225504535589285
1 -0.0210488517799503
1.1 -0.0186663620856215
1.2 -0.0159439587744773
1.3 -0.0133900212630883
1.4 -0.0113862906826005
1.5 -0.010147146041404
1.6 -0.00972541962178312
1.7 -0.0100415588543043
1.8 -0.0109224281582721
1.9 -0.0121415521086423
2 -0.0134590485299376
};
\addplot [very thick, green!50!black, dashed]
table {%
0 0
0.1 -2.59283668751985e-05
0.2 -0.000147604310734974
0.3 -0.000303564735515872
0.4 -0.000100817161348252
0.5 0.000880782691703432
0.6 0.00260126341946369
0.7 0.00464757870947332
0.8 0.00660709149449667
0.9 0.00817843107794753
1 0.00913754579678381
1.1 0.00935630030295904
1.2 0.00885627732930645
1.3 0.00782167844664107
1.4 0.00653619790797634
1.5 0.00528215818855298
1.6 0.00425993090963975
1.7 0.00355491345413332
1.8 0.0031538212595286
1.9 0.0029948230474668
2 0.00301508396435256
};
\addlegendentry{$u_y$}
\addplot [very thick, green!50!black, mark=o, mark size=2, mark options={solid}, only marks, forget plot]
table {%
0 0
0.1 -1.36944858993409e-05
0.2 -1.90585482406414e-05
0.3 0.000219920572782648
0.4 0.00105709918235956
0.5 0.00247344585005371
0.6 0.00400754339539194
0.7 0.00520792148408095
0.8 0.00584703274587503
0.9 0.00587268246489769
1 0.00537974783994764
1.1 0.00456606461068941
1.2 0.00365647250583536
1.3 0.00283631686691294
1.4 0.00222061034035187
1.5 0.00185707162174221
1.6 0.00174455546832831
1.7 0.00185065203611874
1.8 0.00212282980934953
1.9 0.00249603300992162
2 0.00290032188053359
};
\end{axis}

\end{tikzpicture}}~
	\scalebox{0.45}{
\begin{tikzpicture}

\begin{axis}[
legend cell align={left},
legend style={fill opacity=0.8, draw opacity=1, text opacity=1, draw=white!80!black},
tick align=outside,
tick pos=left,
x grid style={white!69.0196078431373!black},
xlabel={Time},
xmin=-0.1, xmax=2.1,
xtick style={color=black},
y grid style={white!69.0196078431373!black},
ylabel={Displacement},
ymin=-0.0403567638952452, ymax=0.0118138419090113,
ytick style={color=black}
]
\addplot [very thick, red, dashed]
table {%
0 0
0.1 -9.95948177808011e-05
0.2 -0.000803396476905668
0.3 -0.00306748651266507
0.4 -0.0075564188569748
0.5 -0.0138778859365614
0.6 -0.0208231157112783
0.7 -0.0272005793083875
0.8 -0.032321408498514
0.9 -0.0359043525358027
1 -0.0377963308769123
1.1 -0.0379228385953162
1.2 -0.0364147269561444
1.3 -0.0336549484994603
1.4 -0.030188839016599
1.5 -0.0265823093670443
1.6 -0.0233027631516495
1.7 -0.020657060059876
1.8 -0.01879266248127
1.9 -0.0177406878815139
2 -0.0174661154135103
};
\addlegendentry{$u_x$}
\addplot [very thick, red, mark=o, mark size=2, mark options={solid}, only marks, forget plot]
table {%
0 0
0.1 -9.97788626689834e-05
0.2 -0.000803696309918647
0.3 -0.00306174017488702
0.4 -0.00752276081149372
0.5 -0.013794436678556
0.6 -0.0207069214041007
0.7 -0.0270975391017667
0.8 -0.0322509828479612
0.9 -0.0358582598917125
1 -0.0377868888232678
1.1 -0.0379853727223245
1.2 -0.0365587623736849
1.3 -0.0338411271076007
1.4 -0.0303570410698479
1.5 -0.026688196101168
1.6 -0.0233313440365979
1.7 -0.020619516853906
1.8 -0.0187196212678701
1.9 -0.0176736306815501
2 -0.0174454467149507
};
\addplot [very thick, green!50!black, dashed]
table {%
0 0
0.1 -2.59283668751985e-05
0.2 -0.000147604310734974
0.3 -0.000303564735515872
0.4 -0.000100817161348252
0.5 0.000880782691703432
0.6 0.00260126341946369
0.7 0.00464757870947332
0.8 0.00660709149449667
0.9 0.00817843107794753
1 0.00913754579678381
1.1 0.00935630030295904
1.2 0.00885627732930645
1.3 0.00782167844664107
1.4 0.00653619790797634
1.5 0.00528215818855298
1.6 0.00425993090963975
1.7 0.00355491345413332
1.8 0.0031538212595286
1.9 0.0029948230474668
2 0.00301508396435256
};
\addlegendentry{$u_y$}
\addplot [very thick, green!50!black, mark=o, mark size=2, mark options={solid}, only marks, forget plot]
table {%
0 0
0.1 -2.59190156482594e-05
0.2 -0.000147288431823308
0.3 -0.0003026545727759
0.4 -0.000101619291297156
0.5 0.000875812010284843
0.6 0.00259811857843158
0.7 0.00465657244265665
0.8 0.00663097432670901
0.9 0.00821543542680101
1 0.00919262306963021
1.1 0.00944245073609057
1.2 0.00898151385400761
1.3 0.00797415311066069
1.4 0.00668545222590707
1.5 0.0053955037260577
1.6 0.00431700353269677
1.7 0.00355257468270622
1.8 0.00310623812280881
1.9 0.00292941046967413
2 0.00296201982194988
};
\end{axis}

\end{tikzpicture}}
	\caption{Displacement. Left panel: initial guess; right panel: computed using the calibrated model.}
\end{subfigure}~
	\begin{subfigure}[t]{0.5\textwidth}
	\scalebox{0.45}{
\begin{tikzpicture}

\begin{axis}[
legend cell align={left},
legend style={fill opacity=0.8, draw opacity=1, text opacity=1, draw=white!80!black},
tick align=outside,
tick pos=left,
x grid style={white!69.0196078431373!black},
xlabel={Time},
xmin=-0.1, xmax=2.1,
xtick style={color=black},
y grid style={white!69.0196078431373!black},
ylabel={Stress},
ymin=-0.00103084803638525, ymax=0.00106765804141765,
ytick style={color=black}
]
\addplot [very thick, red, dashed]
table {%
0 0
0.1 -5.40215732442364e-05
0.2 -0.000291169431858713
0.3 -0.000645538807551012
0.4 -0.000716906855492529
0.5 -0.00027244352816315
0.6 0.000383454724948506
0.7 0.000809266444479423
0.8 0.0008948148397385
0.9 0.000797872802454762
1 0.000639708896169788
1.1 0.000417717004387516
1.2 0.000124510474593818
1.3 -0.000185862386399486
1.4 -0.000432582077359381
1.5 -0.0005601062373667
1.6 -0.00055963135624632
1.7 -0.000463134200525664
1.8 -0.000324195592410646
1.9 -0.000192096622286243
2 -9.16120971929781e-05
};
\addlegendentry{$\sigma_{xx}$}
\addplot [very thick, blue, dashed]
table {%
0 0
0.1 1.61055404938681e-05
0.2 0.00011854641880097
0.3 0.000345415356521786
0.4 0.000494136006646537
0.5 0.000287294089399512
0.6 -0.00016928815066185
0.7 -0.000513807260640971
0.8 -0.000615101967777028
0.9 -0.000573753938693589
1 -0.00046674045790449
1.1 -0.000291522405107775
1.2 -5.1529507539005e-05
1.3 0.000199774860631218
1.4 0.00038968492063139
1.5 0.000469679861340607
1.6 0.000436288187116814
1.7 0.000323345294154287
1.8 0.000184252035986012
1.9 6.61943003940618e-05
2 -1.01381478151038e-05
};
\addlegendentry{$\sigma_{yy}$}
\addplot [very thick, green!50!black, dashed]
table {%
0 0
0.1 1.21364586666931e-05
0.2 6.76134815754413e-05
0.3 0.000136905839570173
0.4 8.9892112870764e-05
0.5 -0.000104991746110473
0.6 -0.000289263045875102
0.7 -0.000347611677417378
0.8 -0.000303159161210242
0.9 -0.000216603822905519
1 -0.000134310895912814
1.1 -6.14214726889601e-05
1.2 1.20937208215865e-05
1.3 7.89106662774188e-05
1.4 0.000127136854378501
1.5 0.00015116094737427
1.6 0.000153374805659098
1.7 0.000139562626682816
1.8 0.000115351312917618
1.9 8.48267908352476e-05
2 5.16918187454772e-05
};
\addlegendentry{$\sigma_{xy}$}
\addplot [very thick, red, mark=o, mark size=2, mark options={solid}, only marks, forget plot]
table {%
0 0
0.1 -0.000103434992456736
0.2 -0.000491346127084675
0.3 -0.000935461396485116
0.4 -0.000881929805763214
0.5 -0.000310298961394885
0.6 0.000323303418787673
0.7 0.000752634713960163
0.8 0.000958032713466345
0.9 0.000972271401517523
1 0.000820742659831226
1.1 0.000551154029176015
1.2 0.000238179051346441
1.3 -4.76971136876944e-05
1.4 -0.000258926776509281
1.5 -0.000378862181624334
1.6 -0.000411350107040929
1.7 -0.000371657010335648
1.8 -0.000279200397007652
1.9 -0.000155366445887432
2 -2.20714317341938e-05
};
\addplot [very thick, blue, mark=o, mark size=2, mark options={solid}, only marks, forget plot]
table {%
0 0
0.1 3.37617510153115e-05
0.2 0.00020710451831381
0.3 0.00051018672542164
0.4 0.000639414420415043
0.5 0.000392992161646363
0.6 -1.8031751637527e-05
0.7 -0.000342721718826358
0.8 -0.000522741197606052
0.9 -0.000558144190592338
1 -0.000463045791791583
1.1 -0.000282068717540034
1.2 -7.63375189701814e-05
1.3 0.000102141476089825
1.4 0.000225753881165022
1.5 0.000289970514960369
1.6 0.000303152835777923
1.7 0.000276891653813751
1.8 0.000222602247022678
1.9 0.000151105594643558
2 7.39012077658746e-05
};
\addplot [very thick, green!50!black, mark=o, mark size=2, mark options={solid}, only marks, forget plot]
table {%
0 0
0.1 1.49681223931187e-05
0.2 6.7901915740875e-05
0.3 0.000101184671342686
0.4 7.73686064655319e-06
0.5 -0.000180571179822811
0.6 -0.000343548387633671
0.7 -0.000439277106492022
0.8 -0.000474292495197527
0.9 -0.000461358620034905
1 -0.000415528618649085
1.1 -0.00034684987150019
1.2 -0.000268686521231682
1.3 -0.000194022602926831
1.4 -0.00013430977260867
1.5 -9.67842399175748e-05
1.6 -8.3643440276392e-05
1.7 -9.28939784195245e-05
1.8 -0.000119150420985901
1.9 -0.000155560994915385
2 -0.000194754438529363
};
\end{axis}

\end{tikzpicture}}~
	\scalebox{0.45}{
\begin{tikzpicture}

\begin{axis}[
legend cell align={left},
legend style={fill opacity=0.8, draw opacity=1, text opacity=1, draw=white!80!black},
tick align=outside,
tick pos=left,
x grid style={white!69.0196078431373!black},
xlabel={Time},
xmin=-0.1, xmax=2.1,
xtick style={color=black},
y grid style={white!69.0196078431373!black},
ylabel={Stress},
ymin=-0.00079749294025408, ymax=0.000975400924500052,
ytick style={color=black}
]
\addplot [very thick, red, dashed]
table {%
0 0
0.1 -5.40215732442364e-05
0.2 -0.000291169431858713
0.3 -0.000645538807551012
0.4 -0.000716906855492529
0.5 -0.00027244352816315
0.6 0.000383454724948506
0.7 0.000809266444479423
0.8 0.0008948148397385
0.9 0.000797872802454762
1 0.000639708896169788
1.1 0.000417717004387516
1.2 0.000124510474593818
1.3 -0.000185862386399486
1.4 -0.000432582077359381
1.5 -0.0005601062373667
1.6 -0.00055963135624632
1.7 -0.000463134200525664
1.8 -0.000324195592410646
1.9 -0.000192096622286243
2 -9.16120971929781e-05
};
\addlegendentry{$\sigma_{xx}$}
\addplot [very thick, blue, dashed]
table {%
0 0
0.1 1.61055404938681e-05
0.2 0.00011854641880097
0.3 0.000345415356521786
0.4 0.000494136006646537
0.5 0.000287294089399512
0.6 -0.00016928815066185
0.7 -0.000513807260640971
0.8 -0.000615101967777028
0.9 -0.000573753938693589
1 -0.00046674045790449
1.1 -0.000291522405107775
1.2 -5.1529507539005e-05
1.3 0.000199774860631218
1.4 0.00038968492063139
1.5 0.000469679861340607
1.6 0.000436288187116814
1.7 0.000323345294154287
1.8 0.000184252035986012
1.9 6.61943003940618e-05
2 -1.01381478151038e-05
};
\addlegendentry{$\sigma_{yy}$}
\addplot [very thick, green!50!black, dashed]
table {%
0 0
0.1 1.21364586666931e-05
0.2 6.76134815754413e-05
0.3 0.000136905839570173
0.4 8.9892112870764e-05
0.5 -0.000104991746110473
0.6 -0.000289263045875102
0.7 -0.000347611677417378
0.8 -0.000303159161210242
0.9 -0.000216603822905519
1 -0.000134310895912814
1.1 -6.14214726889601e-05
1.2 1.20937208215865e-05
1.3 7.89106662774188e-05
1.4 0.000127136854378501
1.5 0.00015116094737427
1.6 0.000153374805659098
1.7 0.000139562626682816
1.8 0.000115351312917618
1.9 8.48267908352476e-05
2 5.16918187454772e-05
};
\addlegendentry{$\sigma_{xy}$}
\addplot [very thick, red, mark=o, mark size=2, mark options={solid}, only marks, forget plot]
table {%
0 0
0.1 -5.41468681842903e-05
0.2 -0.000291156146311633
0.3 -0.000641751692840703
0.4 -0.000705435333897939
0.5 -0.000268893029000147
0.6 0.000363175826812333
0.7 0.000787105877554489
0.8 0.000892396913832716
0.9 0.000799084236583553
1 0.000630468877743596
1.1 0.00041376534572363
1.2 0.000142437692781981
1.3 -0.000154779355665758
1.4 -0.000410173942721773
1.5 -0.000558223918398447
1.6 -0.000575526067090724
1.7 -0.000488461737693198
1.8 -0.000351507484425843
1.9 -0.000213975304194602
2 -0.000101675179758687
};
\addplot [very thick, blue, mark=o, mark size=2, mark options={solid}, only marks, forget plot]
table {%
0 0
0.1 1.61539932632015e-05
0.2 0.000118528278562188
0.3 0.000343103672347466
0.4 0.000487400427347728
0.5 0.000287352532819237
0.6 -0.000152530125753977
0.7 -0.000498044669275114
0.8 -0.000616078000547926
0.9 -0.000580542303836162
1 -0.000469852140084626
1.1 -0.000299301390525997
1.2 -7.16358844452918e-05
1.3 0.000174526768548195
1.4 0.000373729236286176
1.5 0.000470838544653861
1.6 0.000452272525647038
1.7 0.000347612312949853
1.8 0.000210599459194406
1.9 8.76401685739032e-05
2 -1.28411676988624e-07
};
\addplot [very thick, green!50!black, mark=o, mark size=2, mark options={solid}, only marks, forget plot]
table {%
0 0
0.1 1.21528651968171e-05
0.2 6.76185534587711e-05
0.3 0.000137227935058711
0.4 9.23362433274154e-05
0.5 -9.89437341003968e-05
0.6 -0.00028364323209461
0.7 -0.000344863642689675
0.8 -0.00029957211815137
0.9 -0.000209839632471931
1 -0.000125633180359805
1.1 -5.67848920501053e-05
1.2 7.39273270228359e-06
1.3 6.78962182756187e-05
1.4 0.000116943310435899
1.5 0.000145802699041033
1.6 0.000152379908433379
1.7 0.000140658801578191
1.8 0.000116611648070491
1.9 8.51953728836199e-05
2 5.09099988297858e-05
};
\end{axis}

\end{tikzpicture}}
	\caption{Stress. Left panel: initial guess; right panel: computed using the calibrated model.}
\end{subfigure}
\caption{The initial and terminal quantities of the upper left point for the viscosity model $\eta_2$.  The dashed lines are true values, and the dots are reproduced values using the calibrated models.}
\label{fig:res2}
\end{figure}

\begin{figure}[htpb]
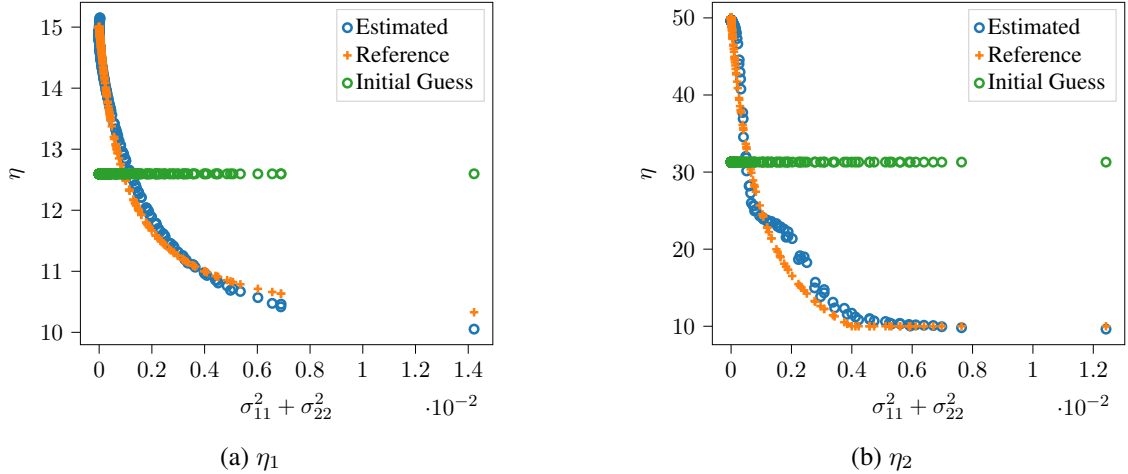

\begin{subfigure}[t]{0.5\textwidth}
\centering
	\scalebox{0.8}{\input{figures/nlv1.tex}}
	\caption{$\eta_1$}
\end{subfigure}~
	\begin{subfigure}[t]{0.5\textwidth}
	\centering
	\scalebox{0.8}{\input{figures/nlv2.tex}}
	\caption{$\eta_2$}
\end{subfigure}
\caption{The viscosity parameter $\eta$ versus ${\sigma_{11}^2+\sigma_{22}^2}$. We collect the stress tensors in the last time step, where the data are simulated using the calibrated neural network, and compute the corresponding $\eta$ using the initial guess of the neural network, the calibrated neural network, and the true relations \Cref{equ:eta12}.}
\label{fig:nlv12}
\end{figure}

Our results are shown in \Cref{fig:res1,fig:res2}, respectively. In the left panel, we show both the initial and terminal displacement of the upper left point. In the right panel, we show both the initial and terminal stress of the upper left point. Since both the $y$ direction displacement and the stress tensor are not present in the training set, they serve as a verification of the estimated constitutive relation \Cref{equ:nlvnn}. The results indicate that the neural-network-based constitutive relations produce similar displacement and stress estimation to the true data. We also show the mapping between estimated $\eta$ and ${\sigma_{11}^2+\sigma_{22}^2}$ in \Cref{fig:nlv12}, where we collect the stress tensors in the last time step, and the data are simulated using the calibrated neural network. We compute the corresponding $\eta$ using the initial guess of the neural network, the calibrated neural network, and the true relations \Cref{equ:eta12}. The reference points correspond to the values calculated according to \Cref{equ:eta12}.

\begin{figure}[htpb]
\begin{subfigure}[t]{0.48\textwidth}
\centering
	\includegraphics[width=1.0\textwidth]{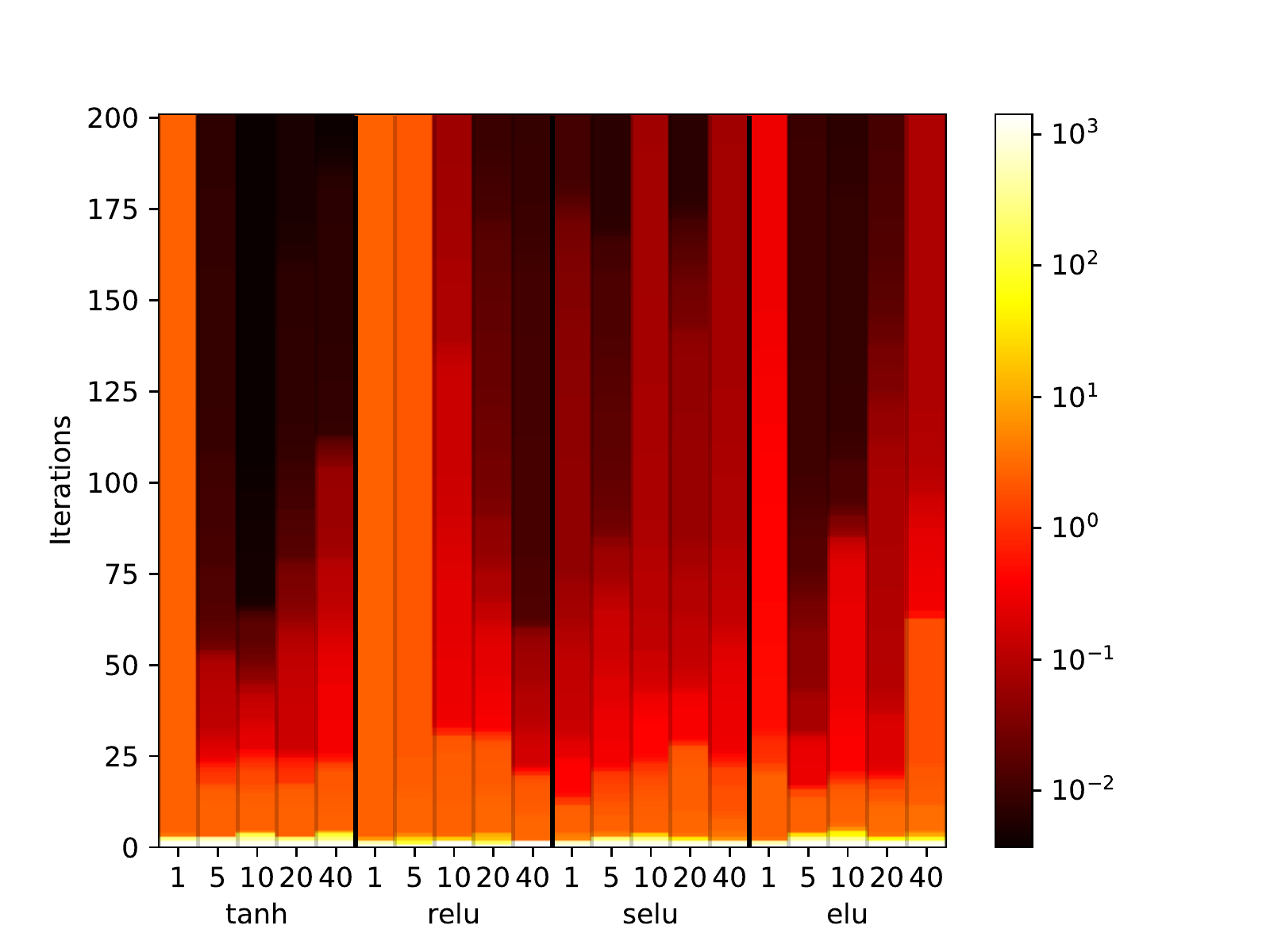}
	\caption{$\eta_1$, depth $=3$}
\end{subfigure}~
	\begin{subfigure}[t]{0.48\textwidth}
\centering
	\includegraphics[width=1.0\textwidth]{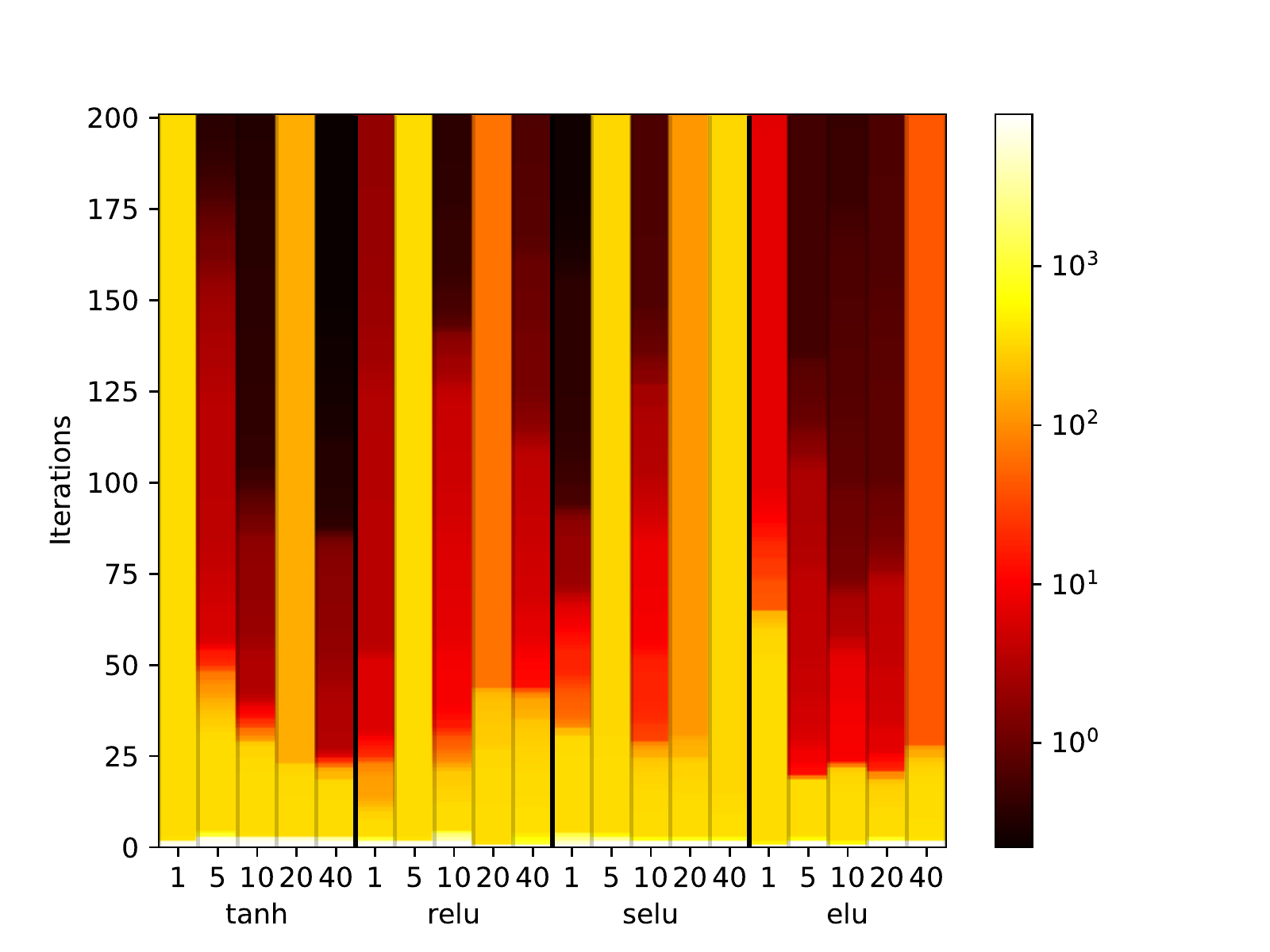}
	\caption{$\eta_2$, depth $=3$}
\end{subfigure}

\begin{subfigure}[t]{0.48\textwidth}
\centering
	\includegraphics[width=1.0\textwidth]{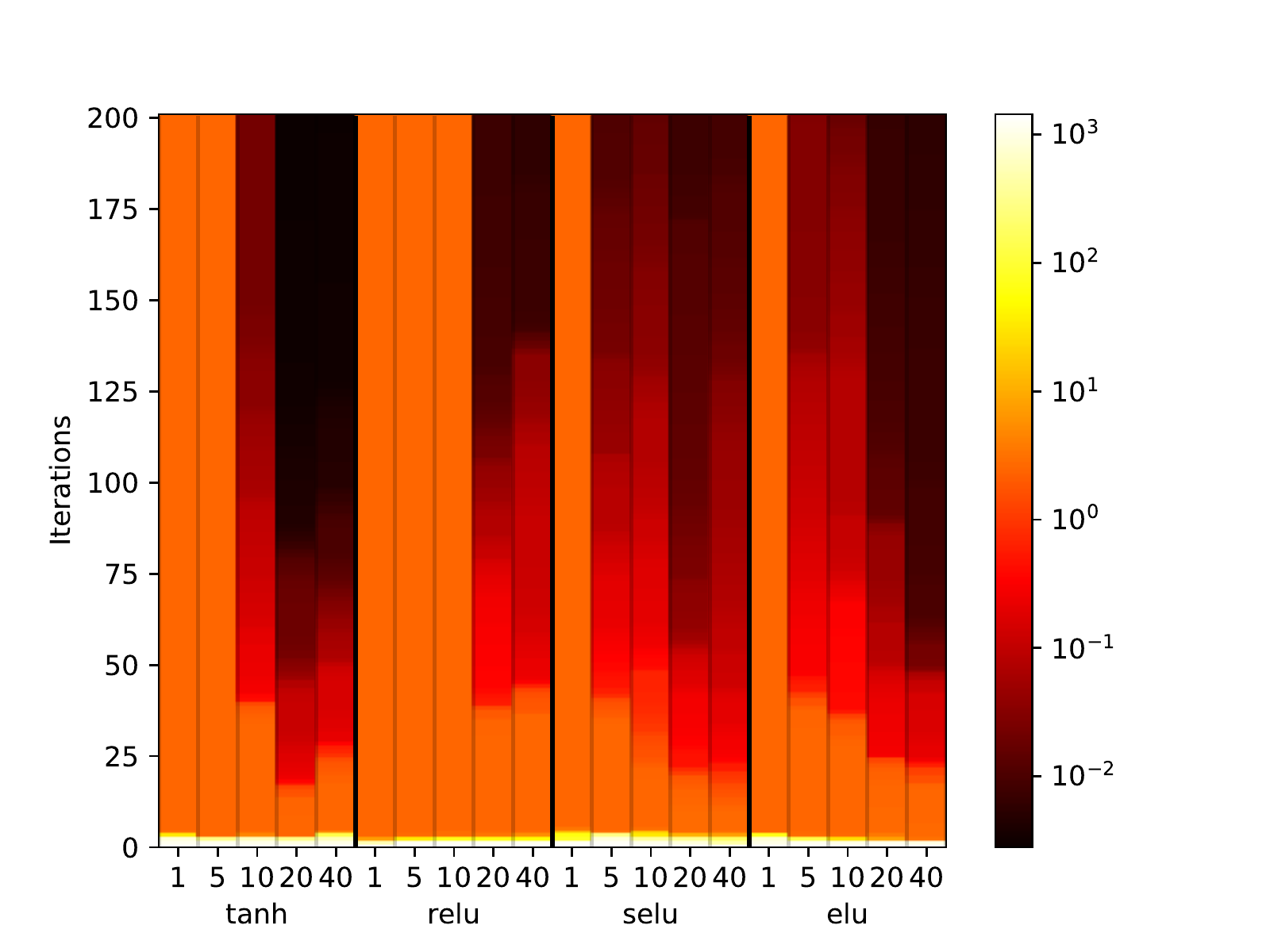}
	\caption{$\eta_1$, depth $=10$}
\end{subfigure}~
	\begin{subfigure}[t]{0.48\textwidth}
\centering
	\includegraphics[width=1.0\textwidth]{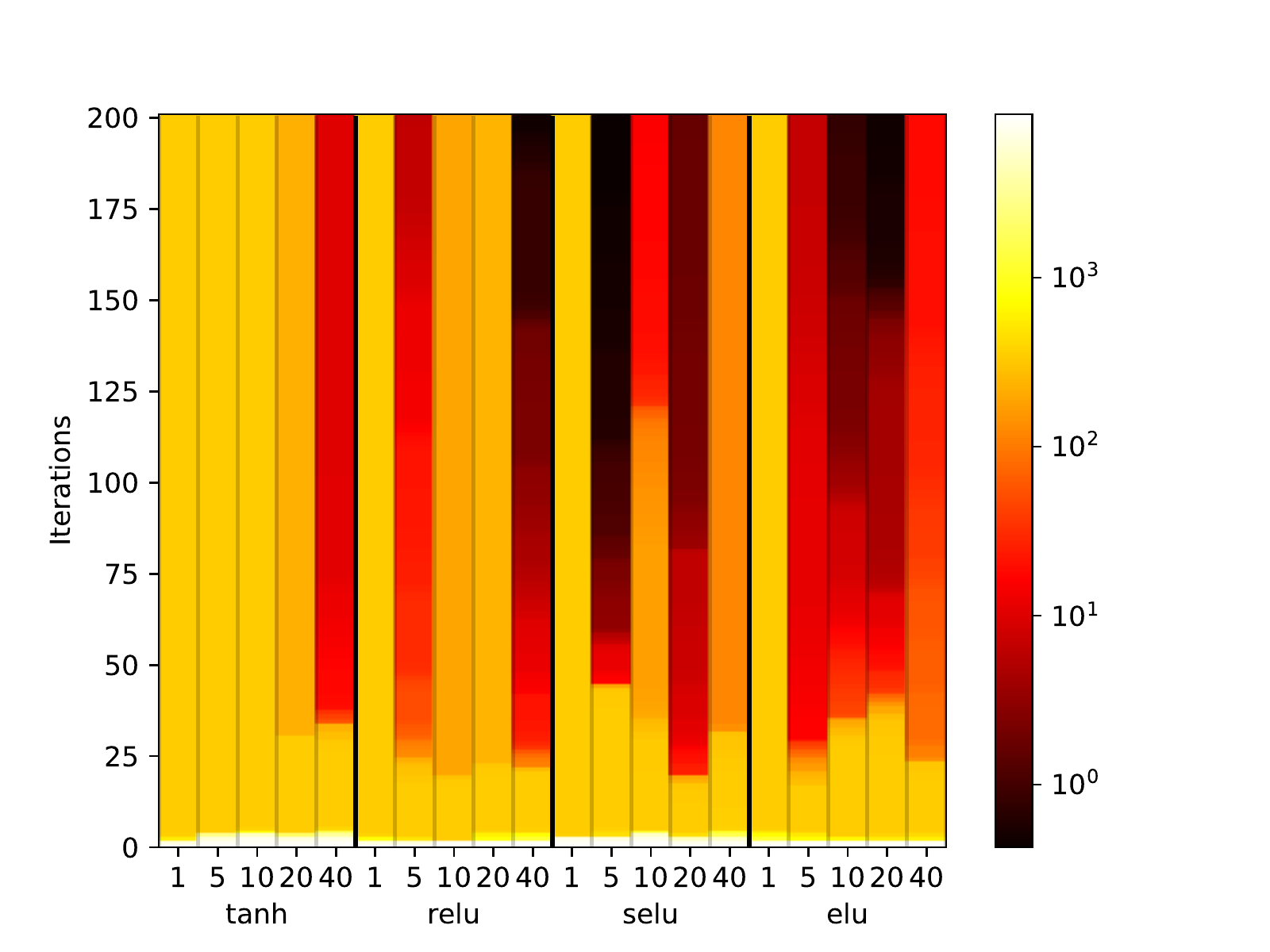}
	\caption{$\eta_2$, depth $=10$}
\end{subfigure}
\caption{Convergence of loss functions using different neural network architectures. The minor ticks 1, 5, 10, 20, and 40 represent the width of the neural network, i.e., number of neurons per layer.}
\label{fig:nn12}
\end{figure}

\paragraph{Convergence Behavior of Different Neural Network Architecture} For fully connected neural network, the neural network architecture---namely the width (number of neurons per layer) and depth (number of hidden layers) of the neural network---plays a crucial role in the convergence of the training. We benchmarked the effect of the choice of the neural network architecture in this problem. The problem setting is as follows: for each test case, we  consider a fully connected neural network with different widths $1$, $5$, $10$, $20$, and $40$;  different depths $1$, $3$, $5$, $10$; and different activation functions $\tanh$, ReLU (rectified linear units \cite{glorot2011deep}), ELU (exponential linear units \cite{clevert2015fast}), and SeLU (scaled exponential linear units \cite{klambauer2017self}). We report the convergence in terms of losses. The results are shown in \Cref{fig:nn12}; it is interesting to see that in general, accuracy increases as width increases, although the convergence is slower. In the context of fully connected neuron network, increasing the depth along makes the optimization more difficult, and we see a better performance for depth 3 compared to depth 10 in general. One potential approach to train deeper neural networks is using a bottleneck architecture similar to residual networks \cite{he2016identity}. This alternative architecture is left to future work. The $\tanh$ activation function leads to more accurate results than the others. Therefore, we stick to this activation function throughout the paper. We also note that the approximating $\eta_2$ is much more difficult than approximating $\eta_1$, which may be due to the non-smoothness of the viscosity parameter. For comprehensiveness, we also report additional reports in \ref{sect:additional}.

\subsection{Learning Constitutive Relations from Experimental Data}\label{sect:4a}

In this example, we consider learning a neural-network-based constitutive relation from 1D experimental data and compare it with parametric Kevin-Voigt models. The data are time series data  of strain-stress pairs. The data is extracted from the mechanical behavior of a viscoelastic polymer composite under different excitation frequencies and strains. For a detailed description of the experiment setup and data collection methods, see \cite{javidan2020experimental}.

In \cite{javidan2020experimental}, the authors fitted the strain-stress data using the Kevin-Voigt model, where the strain stress relation is expressed by a parallel combination of linear elastic and linear viscous terms 
\begin{equation}\label{equ:kevinvoigt}
  \sigma = G'\epsilon + \frac{G''}{\omega}\dot\epsilon
\end{equation}
where $\sigma$ is the shear stress and $\epsilon$ is the shear strain, $G'$ and $G''$ are respectively storage and loss moduli with the unit of shear stress (MPa), and $\omega$ is the angular load frequency. 

\begin{figure}
\centering
\begin{subfigure}[t]{0.33\textwidth}
\includegraphics[width=1.0\textwidth]{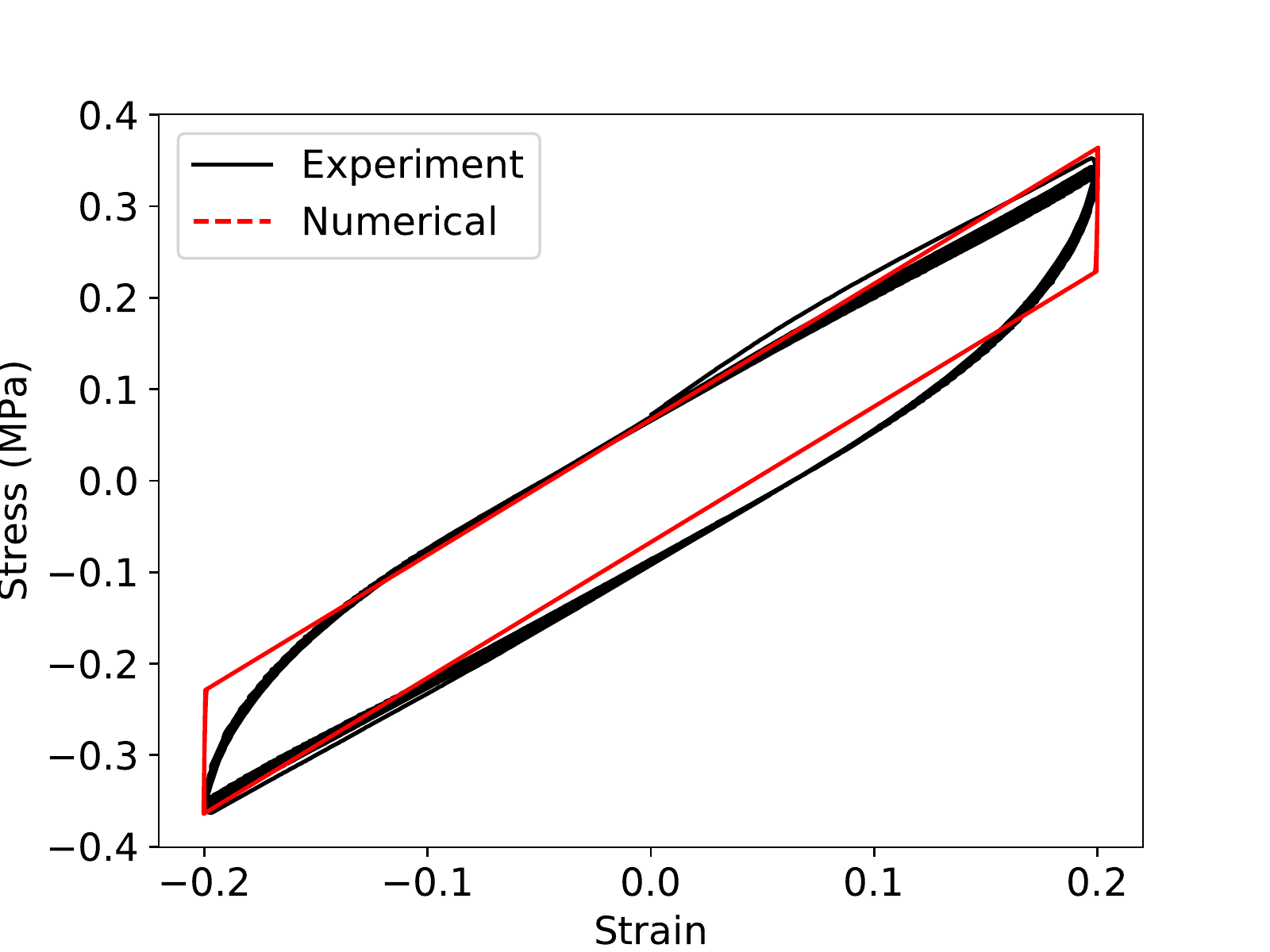}
\end{subfigure}~
	\begin{subfigure}[t]{0.33\textwidth}
\includegraphics[width=1.0\textwidth]{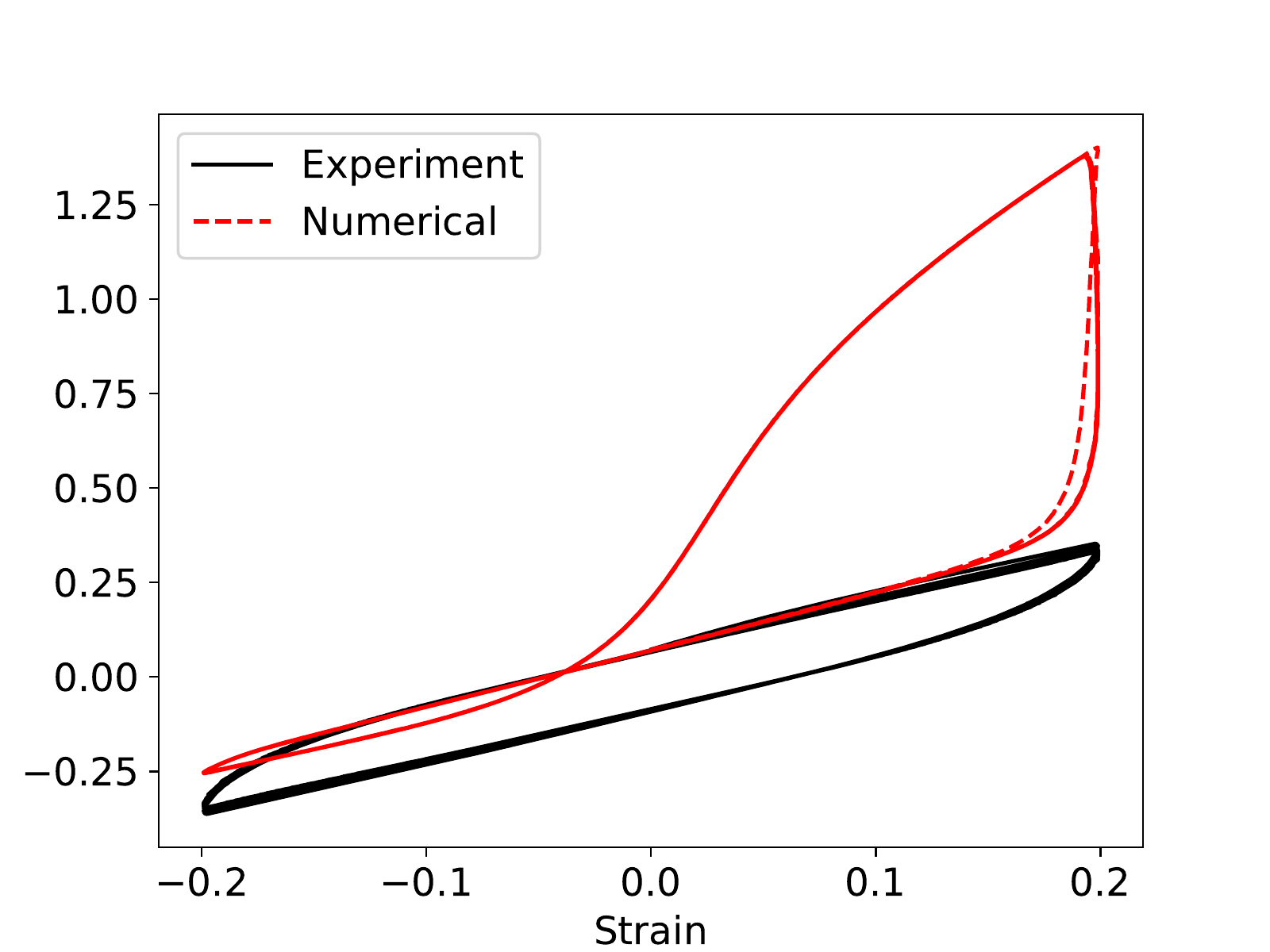}
\end{subfigure}~
\begin{subfigure}[t]{0.33\textwidth}
\includegraphics[width=1.0\textwidth]{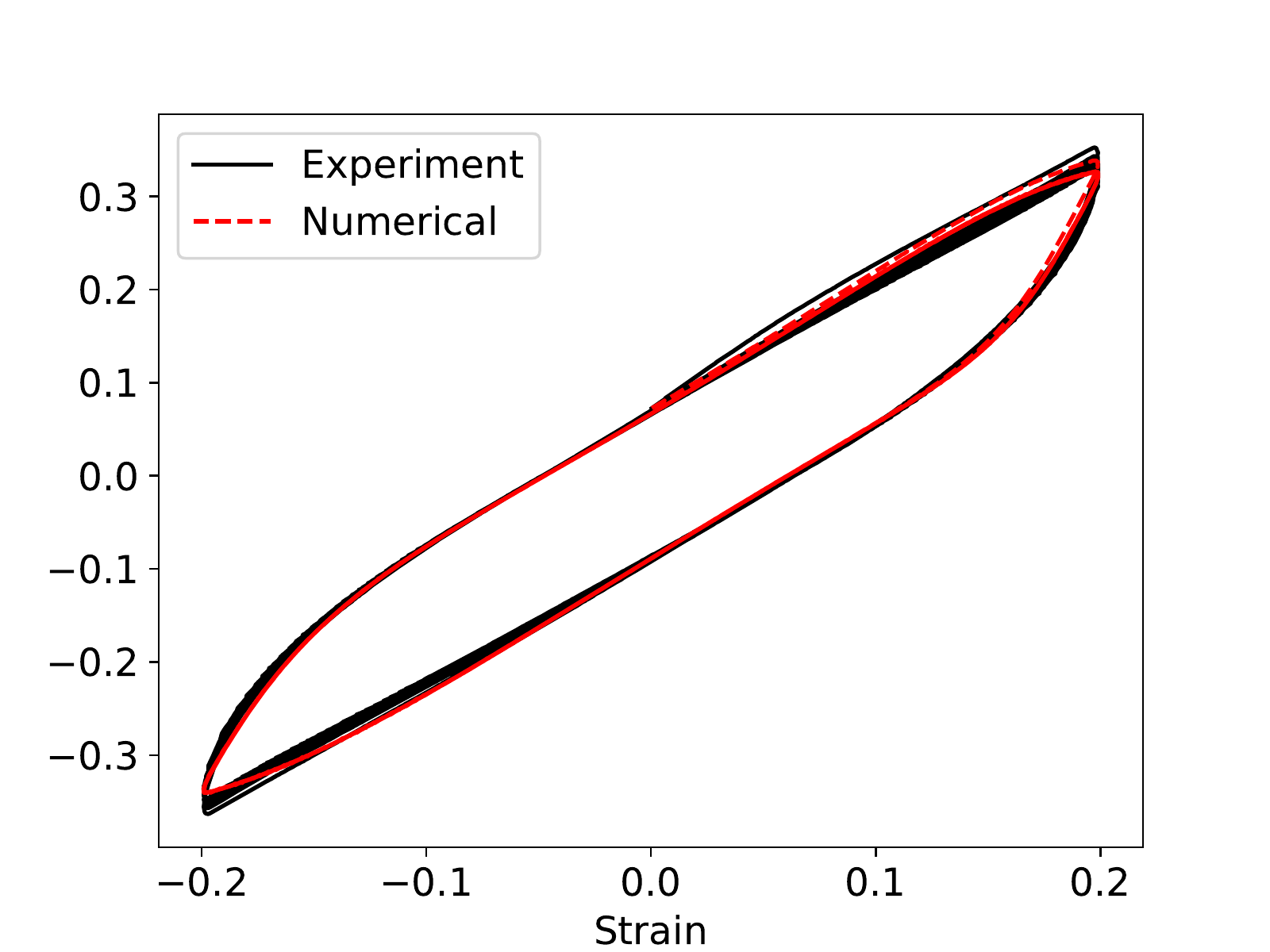}
\end{subfigure}

\begin{subfigure}[t]{0.33\textwidth}
\includegraphics[width=1.0\textwidth]{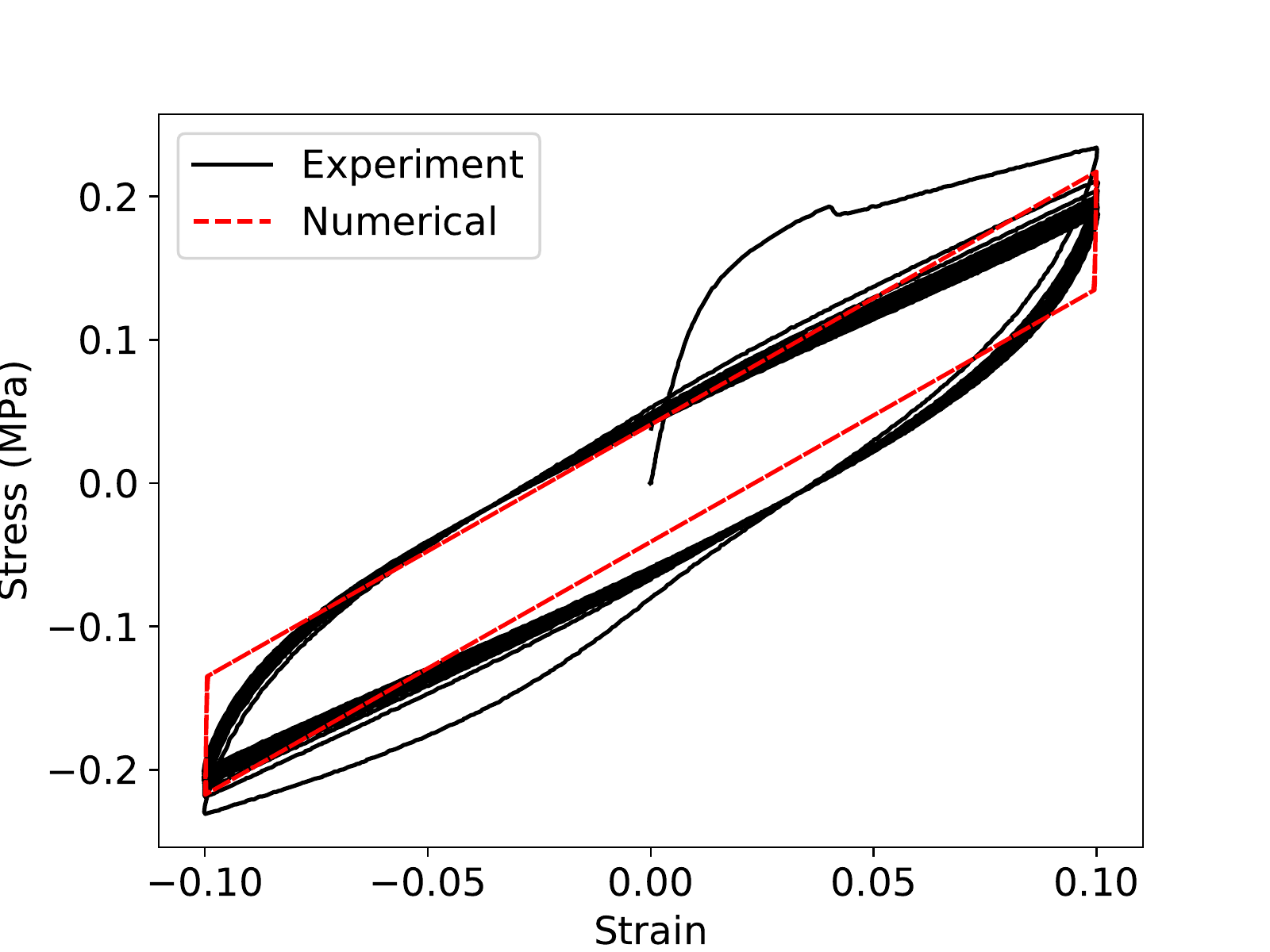}
	\caption{Kevin-Voigt model}
\end{subfigure}~
	\begin{subfigure}[t]{0.33\textwidth}
\includegraphics[width=1.0\textwidth]{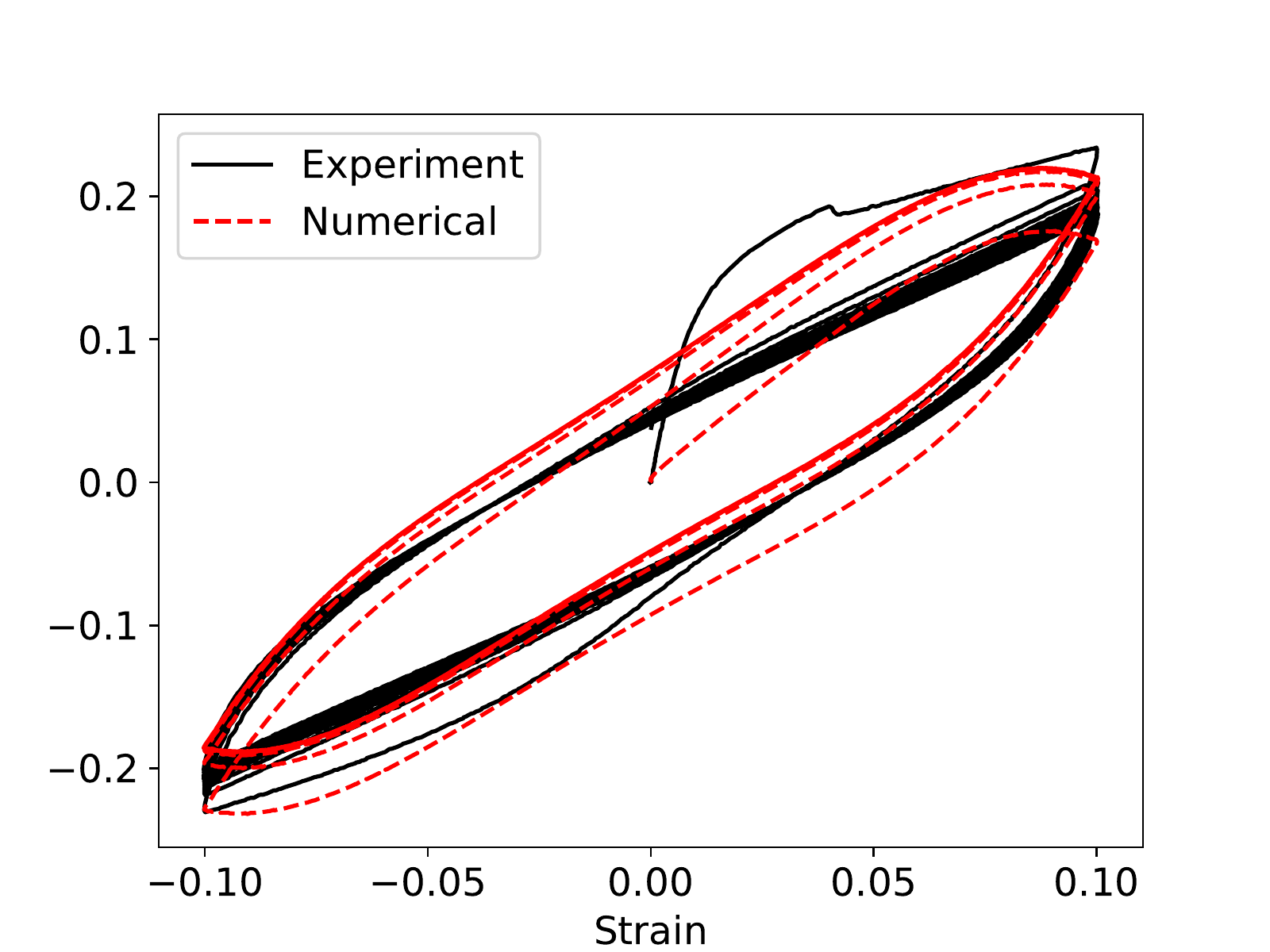}
	\caption{Neural-network-based constitutive model trained with input-output pairs (\Cref{equ:supervised})}
\end{subfigure}~
\begin{subfigure}[t]{0.33\textwidth}
\includegraphics[width=1.0\textwidth]{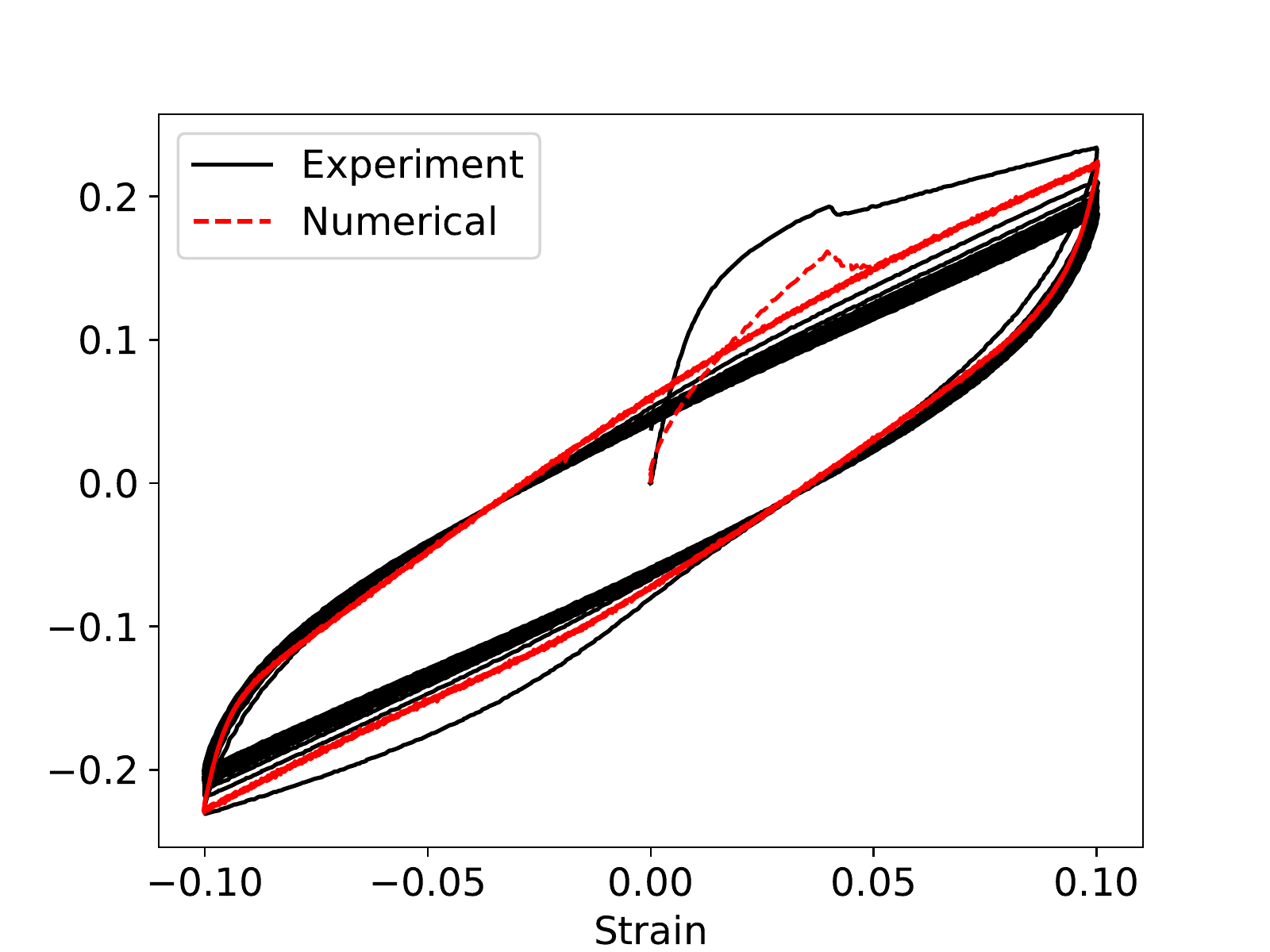}
	\caption{Neural-network-based constitutive model trained with \Cref{equ:rnn} (ours).}
\end{subfigure}
	\caption{Comparison of parametric and neural-network-based constitutive models. The top and bottom rows are results for two different datasets.}
\label{fig:nn_kevin}
\end{figure}

Rather than using a parametric model \Cref{equ:kevinvoigt}, we propose using the neural-network-based constitutive relation \Cref{equ:constitutive_relation}, which is more expressive in representing function forms. The neural network is trained with two approaches, \Cref{equ:supervised,equ:rnn}. \Cref{fig:nn_kevin} shows the results for two typical dataset. In the numerical examples, we use the first 20\% strain-stress pairs as training data. Once the neural network is trained, it is plugged into \Cref{equ:constitutive_relation} to predict the full trajectory. Therefore, the results in \Cref{fig:nn_kevin} demonstrate the models' generalization capability. We can see that the parametric Kevin-Voigt model, despite robustness, fails to capture the creep and recovery processes accurately. The NN model trained with input-output pairs exhibits instability in the prediction. Our NN model outperforms the alternatives, exhibiting accuracy and generalization for predictive modeling. The conclusion from this section justifies the form of neural network we use in \Cref{sect:4}.

\subsection{Learning Constitutive Relations in Coupled Geomechanics and Single-Phase Flow Equations}
\label{sect:4}

\begin{figure}[htpb]
\centering
  \includegraphics[width=0.4\textwidth]{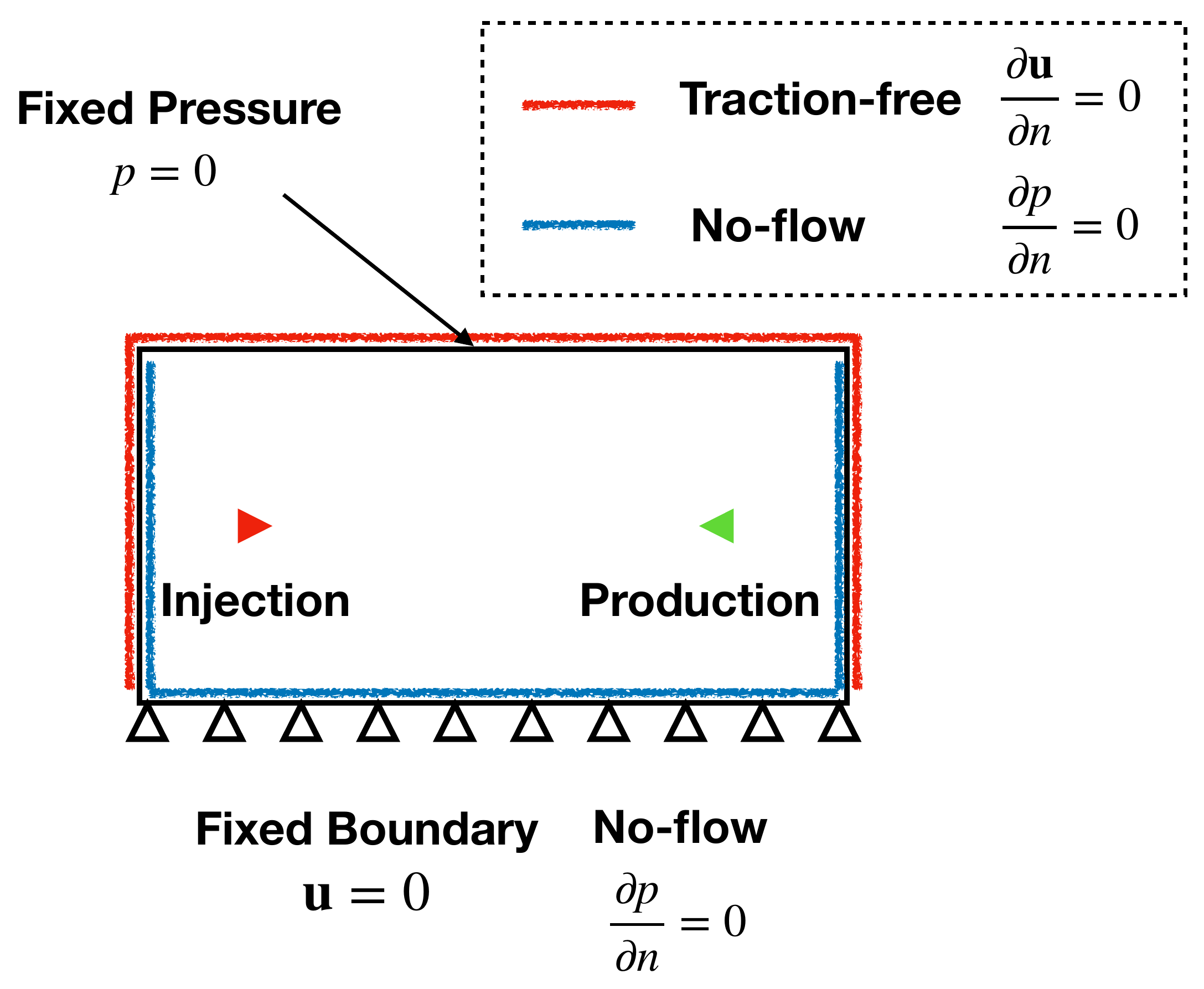}
  \includegraphics[width=0.25\textwidth]{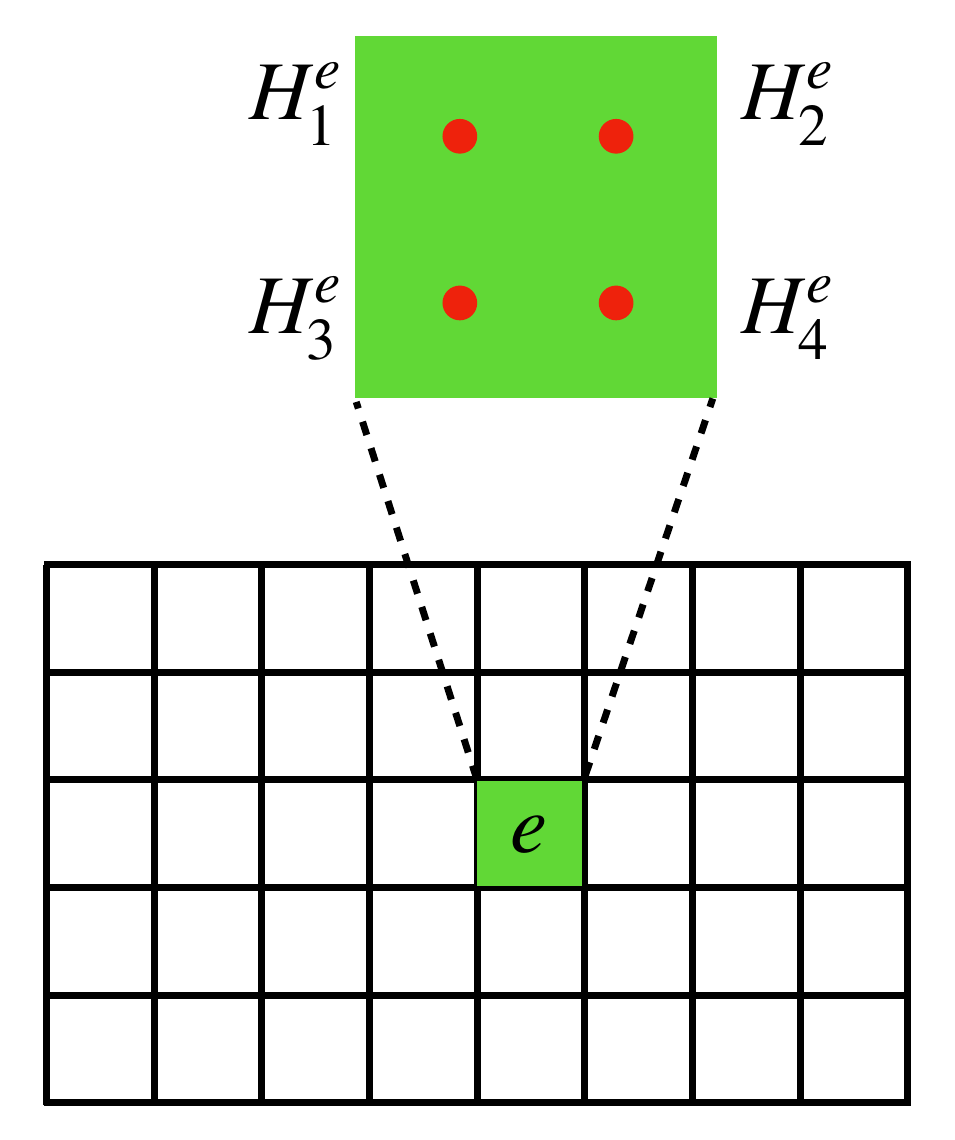}
  \caption{Left panel: geometric setting for \Cref{sect:4}; right panel: each gauss quadrature point  possesses a different linear elastic matrix $H^e_i$, where $e$ is the element index and $i$ is the Gauss quadrature point index in the element $e$.}
  \label{fig:spaceH}
\end{figure}

In this case, we consider learning a neural-network-based constitutive relation in the system of coupled geomechanics and single-phase flow equations. The geometry setting is shown in the left panel of \Cref{fig:spaceH}. The geomechanics equation is the Maxwell model  \Cref{equ:maxwell}. 

We compare two methods:
\begin{enumerate}
	\item Space varying linear elasticity. We assume that the materials are linear elastic, but each gauss quadrature point  possesses a different linear elastic matrix $H^e_i$, where $e$ is the element index and $i$ is the Gauss quadrature point index in the element $e$. See \Cref{fig:spaceH} for a schematic illustration.
	\item The constitutive relation is approximated by 
	\begin{equation}
		\bsigma^{n+1} = \mathcal{NN}_{\bt} (\bsigma^n, \bepsilon^n) + H\bepsilon^{n+1}
	\end{equation}
	where $\bt$ are the weights and biases of the neural network, and $H$ is an unknown tangent elastic matrix, which is symmetric positive definite. Note $H\bepsilon^{n+1}$ is indispensible here for incorporating the strain rate information (together with $\bepsilon^n$)\footnote{Additionally, assuming that $\bsigma^{n+1}$ is linear in $H\bepsilon^{n+1}$ allows for simpler implementation. }. The neural network consists of 3 hidden layers, 20 neurons per layer, and has tanh activation functions. 
\end{enumerate}

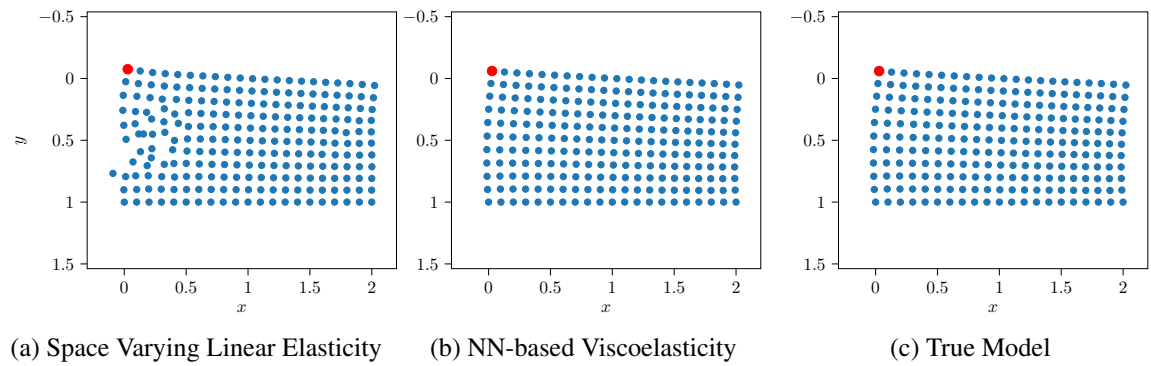
\begin{figure}[htpb]
\begin{subfigure}[t]{0.3\textwidth}
	\scalebox{0.6}{
\begin{tikzpicture}

\definecolor{color0}{rgb}{0.12156862745098,0.466666666666667,0.705882352941177}

\begin{axis}[
tick align=outside,
tick pos=left,
x grid style={white!69.01960784313725!black},
xlabel={$x$},
xmin=-0.3, xmax=2.2,
xtick style={color=black},
y dir=reverse,
y grid style={white!69.01960784313725!black},
ylabel={$y$},
axis equal,
ymin=-0.2, ymax=1.2,
ytick style={color=black}
]
\addplot [only marks,  draw=color0, fill=color0, colormap/viridis]
table{%
x                      y
0.0279517446686601 -0.0749553592887729
0.128475453158363 -0.0605699439274414
0.229410261888852 -0.0477888804136968
0.330825526239181 -0.039010266634784
0.432408596982598 -0.0317466947204165
0.533696337274082 -0.0256116491528994
0.635072835191939 -0.0206799673074868
0.736033602815457 -0.0148776973143798
0.836729785163401 -0.00946608368629199
0.936884084161175 -0.0039106788956536
1.03685836279921 0.00156279800169264
1.13648189850279 0.00681410354730295
1.23549911237969 0.0118514976818002
1.3341600759894 0.0169139380764111
1.43230580167478 0.0216878883930129
1.53001407132732 0.0264256583126613
1.6276216435773 0.0311294929935035
1.72549324338071 0.0363049304622933
1.82396438153172 0.0416870278644418
1.92323787563199 0.0485980694514573
2.02276038356363 0.0569106286864285
0.0119127184462704 0.0260644716878856
0.113961234733571 0.0417903508013043
0.217709153294317 0.0515830095705116
0.323469916539823 0.0615485495501692
0.425085780562907 0.069568674553373
0.528175815735225 0.0742069484773472
0.629374196089692 0.0801847906443547
0.730811397302793 0.085152523262145
0.831390512205546 0.0907626887482205
0.931839208458175 0.0962093922695793
1.03179370741443 0.101534981931864
1.13124778774305 0.10672080465513
1.23034513897592 0.111748259018272
1.32903304676735 0.116566822408248
1.42733146506892 0.121246354436701
1.52523028750317 0.125763650928739
1.6223633500877 0.130505663751007
1.71971639301272 0.135180899875062
1.81724212689811 0.140823405700841
1.91507159822807 0.14721537262476
2.01375623007642 0.155823430424247
-0.00848340701125286 0.136392149961883
0.101997003121939 0.144632130102685
0.208307493994002 0.155353564199738
0.306744195417377 0.167178003652215
0.421008591410724 0.169738170194933
0.521492182608883 0.176866462912749
0.624697427975359 0.180651529407428
0.725651675662943 0.186062241631077
0.826606470320862 0.1912354248783
0.927241004412956 0.196356225109095
1.02732301699193 0.201495870585763
1.12676730263816 0.206481820604994
1.22575614784555 0.211245335823447
1.32432650277487 0.215651931620822
1.42236971907196 0.219748882992393
1.5199654793115 0.224104278632205
1.61713139683866 0.228145241007041
1.71288731888356 0.233049830961542
1.80959908048744 0.237428722259951
1.90796637467172 0.242698136972509
2.00350303805034 0.251041995807131
-0.0116222876296989 0.25753346986463
0.0906804510388916 0.267367198554629
0.180483583152164 0.273761191264714
0.32337587244522 0.244918073458225
0.407664385962368 0.287546088356242
0.519241907852899 0.277121385249799
0.619703435465837 0.283263183966003
0.7211973534111 0.287540663527999
0.822409023701281 0.291882469131136
0.923109942691378 0.2965819286873
1.02323615216596 0.301397494362925
1.12276883697625 0.306064142669786
1.22169693141526 0.310274043790461
1.32025154712073 0.314306633571949
1.41800917160615 0.317617257758435
1.51498309898286 0.320588139089955
1.61088998724332 0.325351415802658
1.7115466042246 0.328183105386489
1.80652749787251 0.332956729921399
1.89891089509015 0.341848971582213
1.99630435047189 0.339161797741526
-0.00373554930338797 0.37851998244415
0.0871661224423727 0.366333191046678
0.220665257144592 0.328110273366211
0.320250423090433 0.350449310882571
0.436631196794803 0.363063915608698
0.518000248265196 0.388453971148032
0.616953654638563 0.387430549122178
0.717937218855606 0.389373071511134
0.818901382539075 0.392866683196182
0.919498889742686 0.396961195170296
1.01962162882354 0.401264398228635
1.11921664553603 0.405412324545347
1.21821355232077 0.409140677007144
1.31634919083799 0.411892996480101
1.41537531681579 0.414051521592979
1.5140812439862 0.418202395075797
1.60875853897351 0.419916692781846
1.70365907138624 0.423606378313848
1.7935044817535 0.438467591279151
1.89667730511293 0.429811285606732
1.99620844893575 0.432290141976472
0.0144084046079175 0.492172096547045
0.154684893366524 0.448969372994598
0.228278734358922 0.451225171949071
0.32801411874207 0.43601914267239
0.405377288313997 0.501175478433518
0.510838676275263 0.489003098798786
0.614731591407661 0.489698858949437
0.715410528076669 0.491592975918236
0.815986892550242 0.49420187261891
0.916396544027935 0.497523753694366
1.01646746207948 0.501116479025444
1.11610280857082 0.504567070368177
1.21526692036589 0.50748724122326
1.31444440404054 0.50913167535684
1.41254525929838 0.511529668165968
1.50915794035367 0.514067747887042
1.61050781981285 0.517921232582604
1.70277553725258 0.523268164635113
1.80570337574006 0.52596341090323
1.89933770431858 0.525232570365318
1.99785678517751 0.523037497240636
0.117202686395937 0.448796335501956
0.13033955021096 0.592273884504749
0.223225424371599 0.568426251316856
0.217759167906765 0.642515023101182
0.391184895176246 0.576822043925548
0.507071393418139 0.587066075063536
0.610266651331556 0.590784908229466
0.712336320882015 0.593262038034419
0.813197681719352 0.595579158474898
0.91356946410338 0.598165792451034
1.01360449456647 0.600926735208022
1.11329894992964 0.603532594898678
1.21271076589326 0.605579490370318
1.31137088138696 0.607536258835397
1.40916001025914 0.609434424093585
1.50684328627428 0.611229510638444
1.59846571713911 0.613606860786149
1.70455053125831 0.618360610397671
1.80195974010193 0.619615394585412
1.90178325423974 0.620324239992131
2.00047464454811 0.618379654795447
-0.0916014136781384 0.767884350740203
0.0702945113791765 0.675404079459449
0.184820909822521 0.704881780384746
0.323161423410233 0.694117578880652
0.401431932236999 0.687804901387419
0.500535818978099 0.686866088833136
0.605922714431408 0.691390400622368
0.708618810532779 0.694359625466759
0.810025908202455 0.696673674663528
0.91061012834371 0.698715642662993
1.01069813439767 0.700659232778871
1.1104093625444 0.702432844662736
1.20961321038069 0.704097885991935
1.30824274985826 0.705800805964132
1.40654880235064 0.707648785993645
1.50465178725671 0.709755371875817
1.60405385173813 0.710197688145206
1.70013224249431 0.712127001128937
1.80060589293261 0.713780739346272
1.90044693058422 0.713804596974714
2.00114089118949 0.714212352158893
0.010807661726554 0.79375725324851
0.0922647069166583 0.787550800671774
0.193254593192582 0.786698310756445
0.302552187458387 0.793901416830373
0.40649819018075 0.797470030489926
0.50283101323278 0.793078695788322
0.603078388996237 0.79319915323274
0.705174892660524 0.795591296108955
0.806550351195505 0.797594516784057
0.907248705381383 0.799114023759542
1.00740447728404 0.800338722000719
1.10710172884884 0.801487041527011
1.20639644385783 0.802718670885173
1.3054136020938 0.804046297635939
1.40431683561647 0.805290521076957
1.50333205032183 0.80600024958523
1.6017911861319 0.807226266121819
1.7001744513857 0.807864455112132
1.79915787583567 0.807741206628781
1.89920148879095 0.808030707100621
2.00023037936973 0.808680234198349
-0.00349054811200798 0.900439656340927
0.100261865498069 0.898129416642735
0.197021110377874 0.894119043085851
0.298290127205526 0.895139094326446
0.401872039128463 0.898088261718448
0.503159140507508 0.897757552256643
0.60234125835797 0.896889831526315
0.702456962787398 0.897580002591464
0.803125350109105 0.898667372819924
0.903591435371038 0.899479815586849
1.00372033064508 0.900084822938921
1.1035497097355 0.900670676761754
1.20318301659883 0.901338741550278
1.30272894497688 0.902021049455781
1.40221948672323 0.902564540494713
1.50143202121439 0.903118025206271
1.60053599932016 0.903612028311564
1.69987261003843 0.90374315797483
1.79947299457185 0.903525229009554
1.89906976318009 0.903293491658867
1.99917420877099 0.903390500411393
0 1
0.1 1
0.2 1
0.3 1
0.4 1
0.5 1
0.6 1
0.7 1
0.8 1
0.9 1
1 1
1.1 1
1.2 1
1.3 1
1.4 1
1.5 1
1.6 1
1.7 1
1.8 1
1.9 1
2 1
};
\addplot [semithick, red, mark=*, mark size=3, mark options={solid}, only marks]
table {%
0.0279517446686601 -0.0749553592887729
};
\end{axis}

\end{tikzpicture}}
	\caption{Space Varying Linear Elasticity}
\end{subfigure}~
	\begin{subfigure}[t]{0.3\textwidth}
	\scalebox{0.6}{
\begin{tikzpicture}

\definecolor{color0}{rgb}{0.12156862745098,0.466666666666667,0.705882352941177}

\begin{axis}[
tick align=outside,
tick pos=left,
x grid style={white!69.01960784313725!black},
xlabel={$x$},
axis equal,
xmin=-0.3, xmax=2.2,
xtick style={color=black},
y dir=reverse,
y grid style={white!69.01960784313725!black},
ymin=-0.2, ymax=1.2,
ytick style={color=black}
]
\addplot [only marks,  draw=color0, fill=color0, colormap/viridis]
table{%
x                      y
0.0279800530890972 -0.0604013718537024
0.128437868339413 -0.0520642105462108
0.229188707120128 -0.0448281596446044
0.33046822563519 -0.0384575405638968
0.432057591257339 -0.0325100984646737
0.533626925526361 -0.0267085504636846
0.63493980744896 -0.02095065013986
0.735896155477977 -0.0152344773236224
0.836484151771263 -0.00957984055903371
0.936718760376426 -0.00398994372576622
1.03660400193535 0.00154715249292252
1.13612340080896 0.00704098554180959
1.23523935338103 0.0124815387002322
1.33390253867038 0.0178187023869515
1.43209360607002 0.022957414571394
1.52989409433095 0.027800016790146
1.62754553840703 0.0323358477869371
1.72542621244094 0.036726712550676
1.82391070922606 0.0413335858170921
1.92311531372165 0.0466341969845876
2.02270000637567 0.0528498366421506
0.0191250119849748 0.0403479063810958
0.120379194069597 0.0487077147377297
0.222106266543189 0.0557244870475743
0.324174097837844 0.0619096672415176
0.426233761854622 0.0677345431849475
0.528024115677526 0.0734664891828415
0.629436136701697 0.0791808520071709
0.730458836988923 0.0848596751016724
0.831113904221275 0.0904783730857143
0.931414356438667 0.0960363807112366
1.0313641693079 0.101544735599044
1.1309580345927 0.107011764334697
1.23018140379863 0.112428150708147
1.32901052022139 0.117750630151662
1.42741204740123 0.122886157628927
1.52538651316783 0.127719339128768
1.62304170119871 0.132207145058552
1.72065084052653 0.136488812167496
1.81857950582373 0.140923769246586
1.9170547777841 0.14601701574477
2.01602920030789 0.152240766758917
0.00974865505633387 0.143159365620618
0.111972676330654 0.151253372670778
0.21477997813354 0.157729667954376
0.31793808875147 0.16349059624729
0.420733064301938 0.16902549890402
0.522884153571208 0.174528855486695
0.624437153806904 0.180011322285672
0.725537389705304 0.185439521770855
0.826252648744536 0.190816726935929
0.926613328397269 0.196147418385082
1.02663235821854 0.201435597707303
1.12631372533745 0.206679827475508
1.22565264100963 0.21186751448
1.3246268095879 0.216969368734041
1.42317408504106 0.22193642147373
1.52114827430572 0.226657970280996
1.61849066099935 0.231027181527653
1.71544880676962 0.235109128280948
1.81268859769953 0.239226361032111
1.91052544811506 0.243921716435672
2.00887322157299 0.249952270614594
0.000878684326152306 0.248580763136122
0.103804931498127 0.255632057756392
0.207737537573349 0.26112327883799
0.312254638915447 0.266196138774059
0.415825731753395 0.271323947016541
0.518214940245894 0.27647356342857
0.619910077299831 0.281497365160537
0.721095824823211 0.286469619368071
0.821869989412862 0.291417074524949
0.922283949979346 0.296343323998085
1.02236416733175 0.301238286317801
1.12212101523916 0.306084648853164
1.22154904938136 0.310860634618245
1.32061819525516 0.31554655718977
1.41925188692788 0.320134091429329
1.51728101354252 0.324649304932304
1.61422563885465 0.328871475172854
1.7102427681686 0.332630766745709
1.80675907645414 0.336185501258067
1.90412947866244 0.340230446064614
2.0020232125262 0.345609196687301
-0.00645067330559578 0.356498015469165
0.0970136688095184 0.362046603613372
0.201277297098835 0.3660716271101
0.307815741625685 0.370164132277106
0.411475190610647 0.374903834115345
0.514046060538411 0.379245278268095
0.615860565131988 0.383549028857167
0.7171165196791 0.387877191133705
0.817933446209861 0.392241348127488
0.918384420932091 0.396621188165205
1.01850916651024 0.400984695323959
1.11832120091364 0.405294433555199
1.21780818747062 0.40951195205589
1.31692525332073 0.413604759543656
1.41558090757375 0.417566908033343
1.51360958011182 0.421438287962413
1.61071096476932 0.42554413519155
1.70512400257884 0.428856154122018
1.80156161891414 0.431548186060428
1.89861591808141 0.435005441483089
1.99616082765156 0.439371157203941
-0.011452847554101 0.466012231525949
0.0923926078132215 0.469848098077313
0.197010148346289 0.473877442091147
0.30386385471469 0.476961440145564
0.407760400741126 0.47948667768177
0.51044008005086 0.482604702071959
0.612306214940782 0.485989875601849
0.713578117574279 0.489558313579549
0.814399070985041 0.4932397759543
0.914859236435776 0.496976472083396
1.01500763331595 0.500712521035843
1.11485805937819 0.50439199116934
1.21438957016479 0.507957667656973
1.31354107351842 0.511350825618713
1.412198908076 0.514504405982963
1.51016492694262 0.517343958196466
1.60704113497084 0.519596238374931
1.70129604734597 0.522358667726873
1.7976523583679 0.52571730126906
1.89448322839337 0.528711306060624
1.99178658257637 0.531849370287959
-0.0139377915658964 0.575518537052751
0.0900834306978003 0.577346191554612
0.195069105373143 0.57956328342627
0.300476815736026 0.581735349457216
0.404675739021282 0.583793161975222
0.507393884260088 0.586017527567621
0.609216218410633 0.588588999059132
0.710428755509044 0.591400113437734
0.811201218506549 0.59436418842837
0.91163592029772 0.597405691548928
1.01178659229806 0.600457799332665
1.11166532389155 0.603456799713143
1.2112432479103 0.606338665090314
1.31044531495742 0.609036355830551
1.40913338400819 0.611483617575456
1.50706173592918 0.613620538897297
1.6037094562803 0.615672542901387
1.6993342685808 0.617866976582672
1.79534504484697 0.620143486578092
1.89201872878679 0.622209262107297
1.98916550899381 0.623871800713423
-0.0141631646633523 0.684199211113744
0.0899261779814595 0.683962564583451
0.194448375494962 0.684657410418688
0.298848029769847 0.685952905639966
0.402375298019862 0.687547251722361
0.504874728033399 0.689293710677391
0.606519331098768 0.691220318237025
0.707586995626495 0.693336523653762
0.808254351582153 0.695586963921362
0.908627750109852 0.697910511926115
1.00876119529672 0.70024809990017
1.1086638425627 0.702543768524825
1.20830071347286 0.704743665007159
1.30758720047697 0.706799442832915
1.40637589073441 0.708681750619259
1.50443711696623 0.710440111428554
1.60166160966831 0.712151599171226
1.69825193623515 0.713808383455644
1.79473199833231 0.715206533010531
1.89138018794123 0.716025394345656
1.9884002657613 0.716069198181734
-0.0123513098560072 0.791827750822548
0.0916109687174408 0.789768396762345
0.195421881203817 0.789394192881529
0.298514565919847 0.790146772821091
0.401008619377805 0.791271686615725
0.502867932490618 0.792543312099056
0.604142499651681 0.79390926679327
0.704962459782713 0.795369670483651
0.805465079779878 0.79691108135341
0.905742662655704 0.798500085253066
1.00584262079898 0.800099370780088
1.10577070005586 0.801673017847737
1.20549151509492 0.803189238967481
1.30492692613412 0.804626279525416
1.40395954605153 0.805987861855408
1.50248818589776 0.807296671734201
1.60050860740887 0.808557508539005
1.69815950121893 0.809708140070854
1.79559381696405 0.810551103729798
1.89262539139168 0.810435423154261
1.98959014601938 0.808837648172147
-0.00841063401505938 0.898095527081093
0.0955864711732018 0.894560216913914
0.197529260707134 0.894341916074973
0.299039623124377 0.89472045392681
0.400308175840254 0.895332762553826
0.501301240747663 0.8960352630613
0.602007006476436 0.896781996861445
0.702465878091961 0.897561055377535
0.802742549104138 0.898367215323219
0.90289288320119 0.899190374318373
1.00294705895946 0.900016874032367
1.10290760228993 0.900833761747976
1.20274978277525 0.901631179300579
1.302424328971 0.902405779056833
1.40187393999861 0.903160075853559
1.50106499030772 0.903893245612574
1.60002117734205 0.904587702608341
1.6988250927297 0.905199290654587
1.79755215967561 0.905616924638782
1.89602123011546 0.905556668928296
1.99285326461602 0.902754876138457
0 1
0.1 1
0.2 1
0.3 1
0.4 1
0.5 1
0.6 1
0.7 1
0.8 1
0.9 1
1 1
1.1 1
1.2 1
1.3 1
1.4 1
1.5 1
1.6 1
1.7 1
1.8 1
1.9 1
2 1
};
\addplot [semithick, red, mark=*, mark size=3, mark options={solid}, only marks]
table {%
0.0279800530890972 -0.0604013718537024
};
\end{axis}

\end{tikzpicture}}
	\caption{NN-based Viscoelasticity}
\end{subfigure}~
\begin{subfigure}[t]{0.3\textwidth}
	\scalebox{0.6}{
\begin{tikzpicture}

\definecolor{color0}{rgb}{0.12156862745098,0.466666666666667,0.705882352941177}

\begin{axis}[
tick align=outside,
tick pos=left,
axis equal,
x grid style={white!69.01960784313725!black},
xlabel={$x$},
xmin=-0.3, xmax=2.2,
xtick style={color=black},
y dir=reverse,
y grid style={white!69.01960784313725!black},
ymin=-0.2, ymax=1.2,
ytick style={color=black}
]
\addplot [only marks,  draw=color0, fill=color0, colormap/viridis]
table{%
x                      y
0.0279924303683926 -0.0598660472222394
0.128484600632421 -0.051212716370722
0.229223480420536 -0.0439356448106439
0.330473864388915 -0.0376694781585028
0.432043884819334 -0.0318296204435259
0.533611274157473 -0.0261236363856527
0.634944000511321 -0.0204780334936946
0.735936862321805 -0.0148992179736619
0.836559372499505 -0.0093951355894959
0.936808161175263 -0.00395846965260399
1.03667948800776 0.00142706692266165
1.13616018969773 0.00677555147817822
1.2352267708402 0.0120883985448515
1.3338530741529 0.0173400605788432
1.43203808584051 0.022465145966067
1.52985829226316 0.0273648973623594
1.62753550591254 0.031973277683202
1.72544126237196 0.0363955378583076
1.823962165116 0.0410480139822284
1.92320373330495 0.0466114057532113
2.02278450056662 0.0534550188394798
0.0186914109731642 0.0408089731095454
0.120174295048369 0.0495158228160227
0.222055253183021 0.0565359184355898
0.324204456623993 0.0626065304991711
0.426320225647283 0.068336772343218
0.528174610997712 0.0739881530954834
0.629665340243207 0.0796021949655966
0.730761635643098 0.0851603849651656
0.831464876524661 0.0906514948160654
0.931779385020764 0.0960822209638106
1.03170621745782 0.101466241404499
1.13124008404495 0.106815032794165
1.23036951961412 0.112127872582824
1.32908062141428 0.117379481540663
1.42736490738961 0.122502345348835
1.52525181763939 0.127386744866583
1.62285255285595 0.131935317314878
1.72041556639552 0.136218665656808
1.81825692181142 0.140657024617941
1.91657383312701 0.145982415371796
2.01534580580979 0.152891516704076
0.00853503883105441 0.143436980795109
0.111307449644394 0.151951507383142
0.214460868392427 0.158394407833341
0.31785517300596 0.164064605809878
0.420858848590867 0.169564764768278
0.523231215825177 0.175018427452027
0.624986925603231 0.180414102443214
0.726234403192812 0.185741002388114
0.82703687902006 0.191012770505841
0.927428701879876 0.196238235836567
1.02742539437621 0.201424245622235
1.12702929638141 0.206572782506288
1.22623000833895 0.211677356224006
1.32499954423676 0.216717434257898
1.42328264348109 0.221655716946747
1.52097934907124 0.22640059343922
1.61807900896509 0.230820617409719
1.71481025691716 0.23487392637514
1.8117838201211 0.238940326201619
1.90924497900205 0.243832310265607
2.00709353181143 0.25065779261364
-0.00134318639692438 0.248604658185873
0.102515582403938 0.256106754073979
0.207079063666944 0.261577002984763
0.311989280508038 0.266624539149237
0.415989062037173 0.271804061716564
0.518792046110604 0.276925209330522
0.620800909080539 0.281885965812839
0.722205632514615 0.286783846595978
0.823114931381898 0.291649712519383
0.923587921653943 0.296491680253277
1.02365433542096 0.301303696663317
1.1233225441008 0.306072853395254
1.22257987773306 0.310781947622217
1.32138488220771 0.315415190541662
1.41964618801643 0.319963497747748
1.51718793712262 0.324463541106766
1.61361328203111 0.328724941626523
1.70920288474412 0.332428464673235
1.80520001684652 0.335889346517166
1.90191381577585 0.340124531833132
1.99900585740195 0.346318891382832
-0.00967599654921518 0.356327618169606
0.0949867940563402 0.362190215141901
0.200194736812409 0.36622895927495
0.307439824485393 0.370480911154343
0.411757834272961 0.375319260997375
0.514819268317741 0.379646270188436
0.617016512811419 0.383920384907594
0.718556967079289 0.388200974633428
0.819562518181421 0.392507400603764
0.92010721594466 0.396823714219347
1.02023117419596 0.401122415701358
1.11994668859423 0.405370025448962
1.21923841851437 0.409530974260782
1.31805721223862 0.413572912740316
1.41630487896935 0.417489997648277
1.51379311549113 0.421306120794415
1.61014027348119 0.425429545017289
1.70367685648477 0.428682112087627
1.79935139065134 0.431292373940575
1.8955557902432 0.434951483784299
1.99207616763381 0.439961329955184
-0.0151192777293308 0.465773781378334
0.0899752418308809 0.469641791876843
0.195632146387453 0.473851805240421
0.303279664561089 0.477186163080466
0.407914827959975 0.479765752339392
0.511204461417589 0.482932922043729
0.61355591558759 0.48631983269722
0.715197264734245 0.489870010576469
0.816269503704838 0.493519639315414
0.916861606470031 0.497216431752884
1.01702251760051 0.500908101380562
1.11676625775239 0.504541081481232
1.21607208447521 0.508058265105743
1.31487938623984 0.511398790720947
1.41307453675774 0.514488819179626
1.51045236003586 0.51726415808303
1.60661484834754 0.519351643300368
1.69988144873385 0.522085219711582
1.79520513448546 0.525573969813799
1.89101466939723 0.528734236317417
1.9871612827009 0.532225424398743
-0.0174087873154682 0.575182584186662
0.087757098219831 0.576907723898132
0.193653947874452 0.579365631357524
0.299713432520674 0.581751018340957
0.404586748355282 0.583910569267154
0.507998452881193 0.586215723743474
0.610408665139704 0.588833190016376
0.71206993620923 0.591661621621685
0.813147860515592 0.594626164498055
0.913745438336203 0.597656828645215
1.01391827322634 0.600690110231492
1.1136800216276 0.603663359683777
1.21300258506167 0.606510941886082
1.31181090948254 0.609160584992134
1.40996743535167 0.61153990466598
1.50724501766027 0.613563434888644
1.60319681222877 0.615444391254225
1.69805921607989 0.617633447200316
1.79306837321847 0.620095214259892
1.88861795278061 0.622323887892781
1.9845801671916 0.624052965680382
-0.0172835833653965 0.683721778499927
0.0877155706352356 0.68336905705555
0.193001503228747 0.684306215797693
0.298076608916229 0.685842425115396
0.40224638560482 0.687539811331881
0.505380023519401 0.689357186439028
0.607602551376838 0.691350500541622
0.709127898145796 0.693516321800183
0.810110825144185 0.69579954138325
0.910654130678378 0.698142671655583
1.01081322656198 0.700489231457466
1.11059962484332 0.702783676402168
1.20998107806496 0.704969098106291
1.30887846565423 0.706989449120204
1.40716182973111 0.708798660648187
1.50464384803553 0.710423312664755
1.60122524726926 0.711993701716878
1.69704716461085 0.71366327700974
1.79262472204583 0.715248492004147
1.88825484155262 0.716217139826542
1.98418321503826 0.716128407510811
-0.0150970279564239 0.791283216045141
0.0895373191952183 0.789090435606265
0.194067750732488 0.788992238430789
0.297868635678453 0.789963687696289
0.400974077910865 0.79118681032263
0.5033585073988 0.792517042157703
0.605087882939154 0.793943231097076
0.706276372898818 0.795463534443008
0.807038868149111 0.797055404554259
0.90745825396627 0.79868384785752
1.00758052088426 0.800312263631921
1.10741379613186 0.801903991990186
1.20692830503801 0.803422765451979
1.30605745553429 0.804836081258547
1.4047036343196 0.806131005825309
1.50277421182493 0.807333096997679
1.60023849481836 0.808507372044641
1.69719511708526 0.809675692398294
1.79378253088423 0.810663707268657
1.88988382476701 0.810679206830652
1.98594845143005 0.808764766158114
-0.0106244602047332 0.897601024032288
0.0939882252511863 0.893990712781948
0.196644872043739 0.894032212249191
0.298707684654996 0.894541247137237
0.400392874295634 0.895232173069972
0.501706482073602 0.89598153216893
0.60266816843453 0.896769970649067
0.703329825130473 0.897592117657695
0.803751608506135 0.898439065450059
0.903982283376451 0.899297676089926
1.00404947793148 0.900153273252445
1.1039573446271 0.900991186769761
1.20368637705907 0.901797149597778
1.30319547248883 0.902559822991193
1.40243192239982 0.903275944422949
1.50135306126081 0.903955980353826
1.59995985724049 0.904617399932477
1.69830697835415 0.905255444303603
1.79644242612079 0.905752323646052
1.89414870602822 0.905837694499008
1.99012766248703 0.902511520869614
0 1
0.1 1
0.2 1
0.3 1
0.4 1
0.5 1
0.6 1
0.7 1
0.8 1
0.9 1
1 1
1.1 1
1.2 1
1.3 1
1.4 1
1.5 1
1.6 1
1.7 1
1.8 1
1.9 1
2 1
};
\addplot [semithick, red, mark=*, mark size=3, mark options={solid}, only marks]
table {%
0.0279924303683926 -0.0598660472222394
};
\end{axis}

\end{tikzpicture}}
	\caption{True Model}
\end{subfigure}
  \caption{Displacement at the terminal time.}
  \label{fig:spacedisp}
\end{figure}

\begin{figure}
\centering
\begin{subfigure}[t]{0.33\textwidth}
\includegraphics[width=1.0\textwidth]{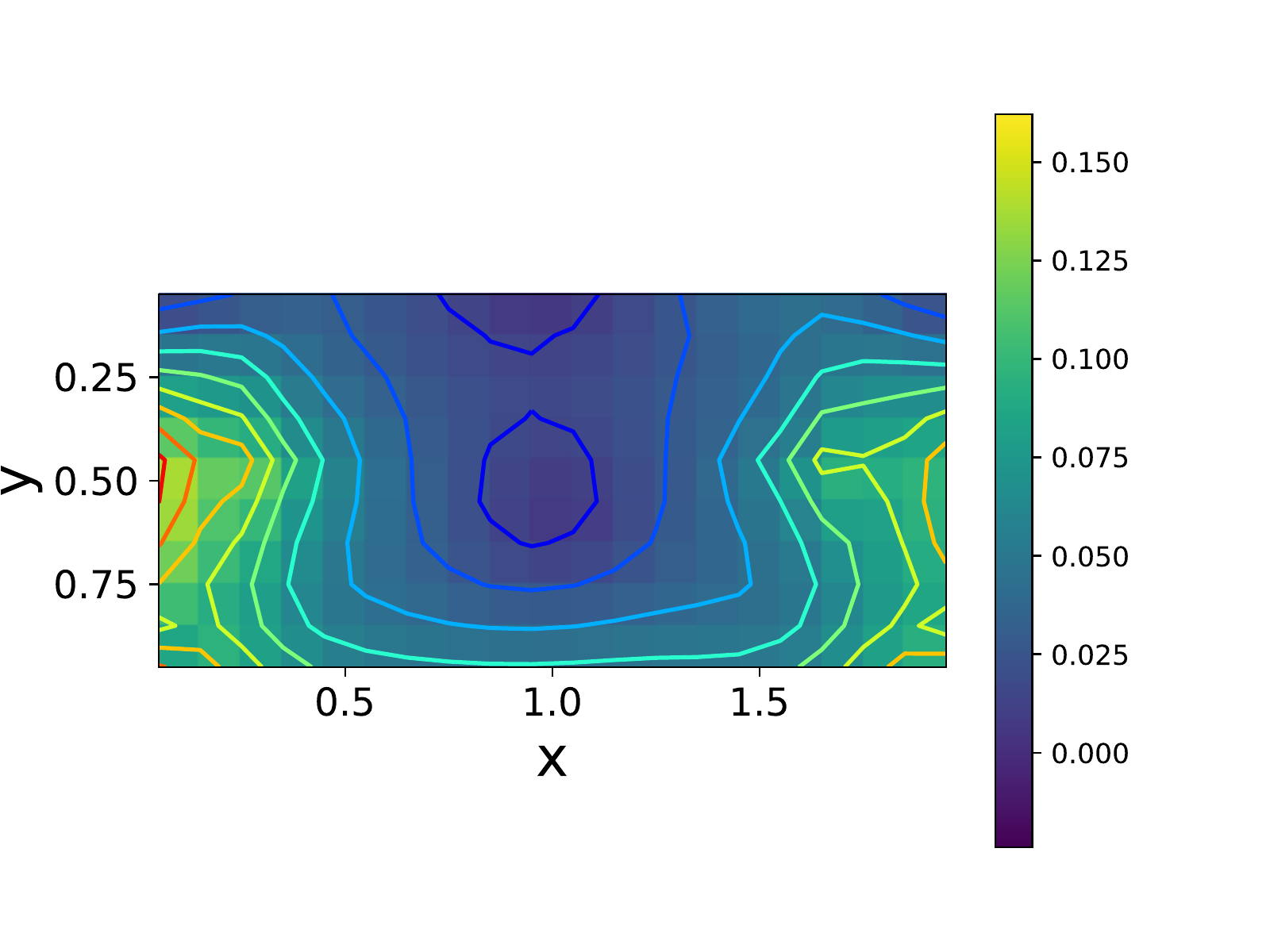}
	\includegraphics[width=1.0\textwidth]{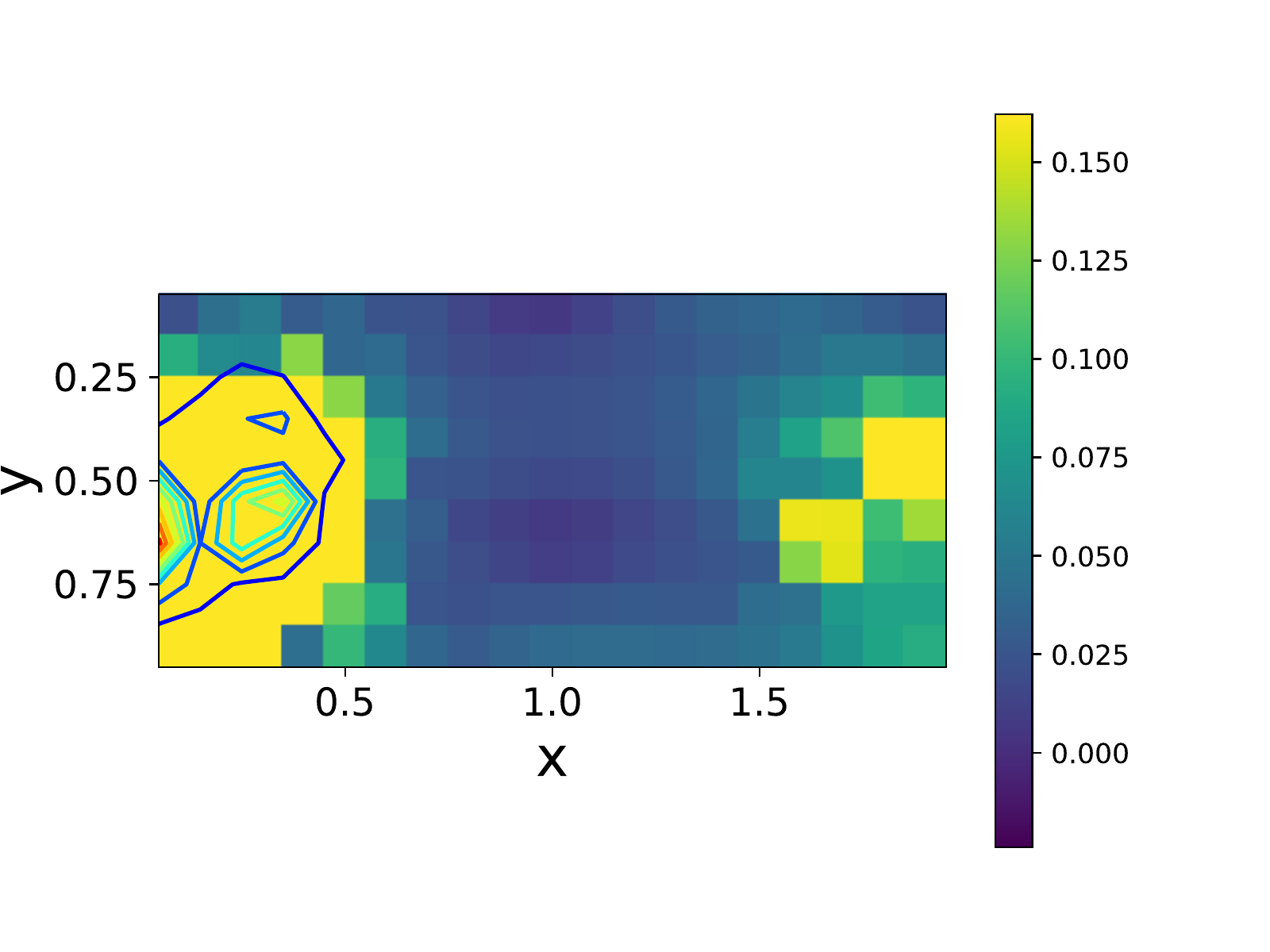}
	\caption{Space Varying Linear Elasticity}
\end{subfigure}~
	\begin{subfigure}[t]{0.33\textwidth}
	\includegraphics[width=1.0\textwidth]{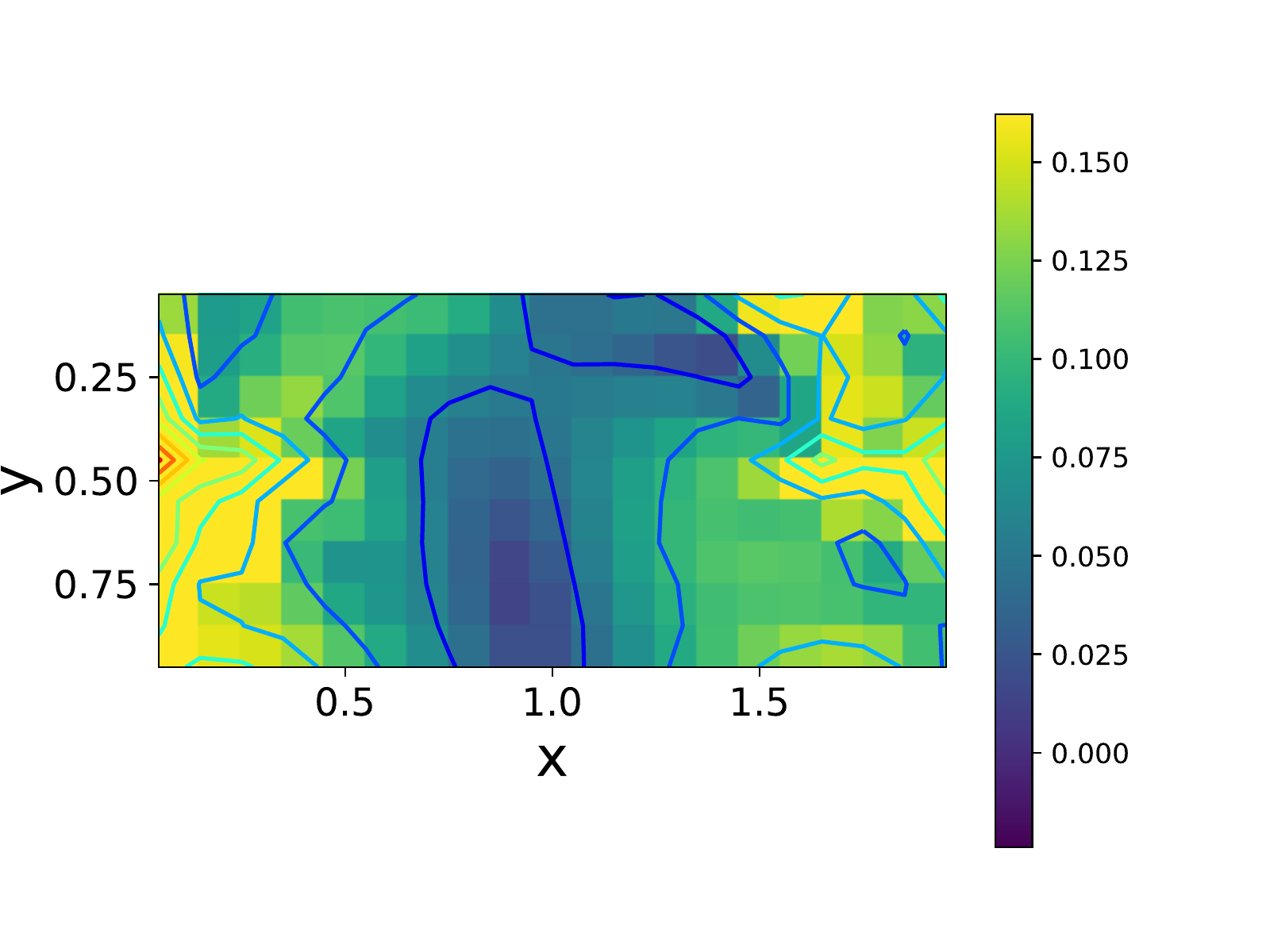}
	\includegraphics[width=1.0\textwidth]{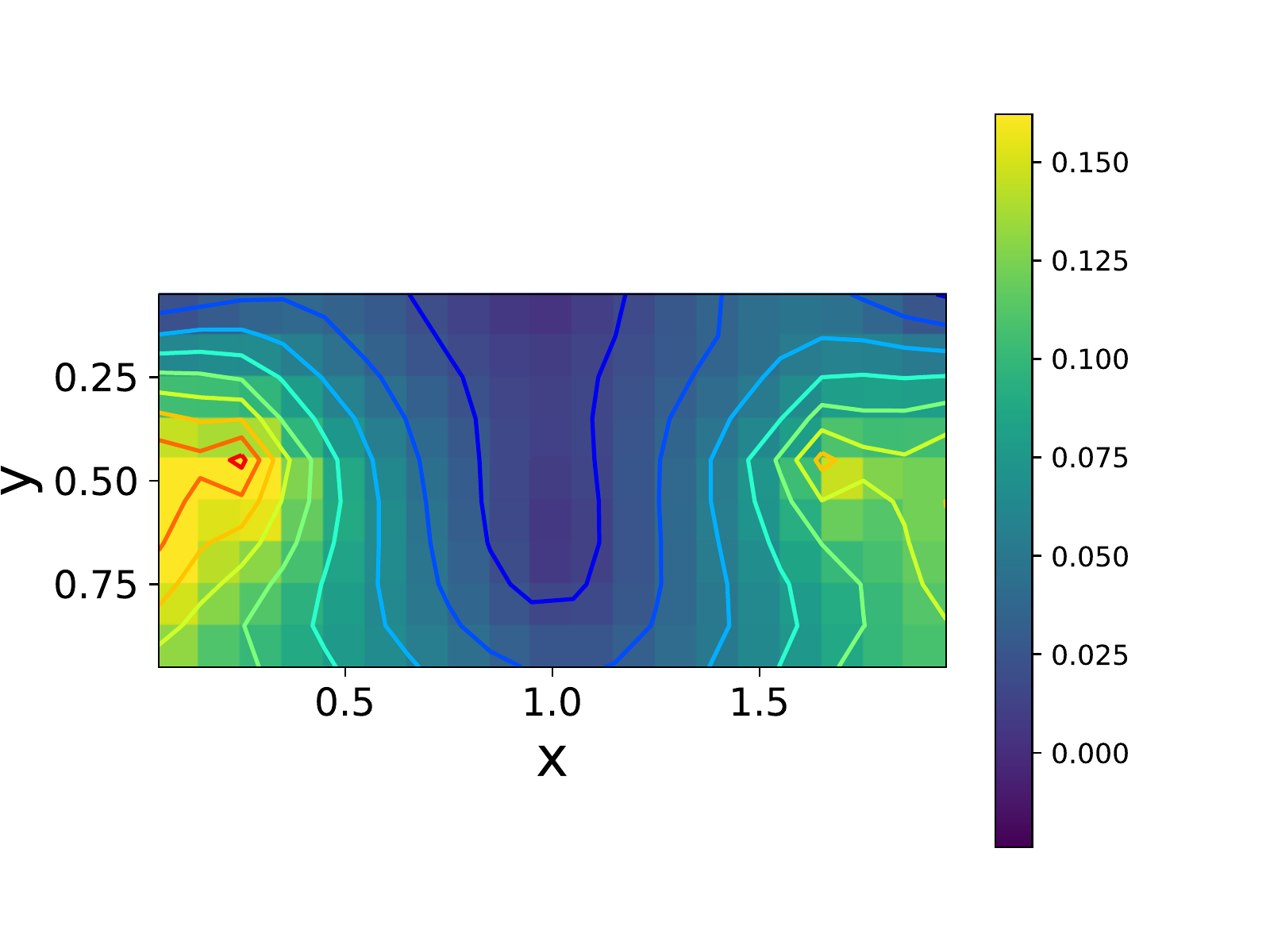}
	\caption{NN-based Viscoelasticity}
\end{subfigure}~
\begin{subfigure}[c]{0.33\textwidth}
	\includegraphics[width=1.0\textwidth]{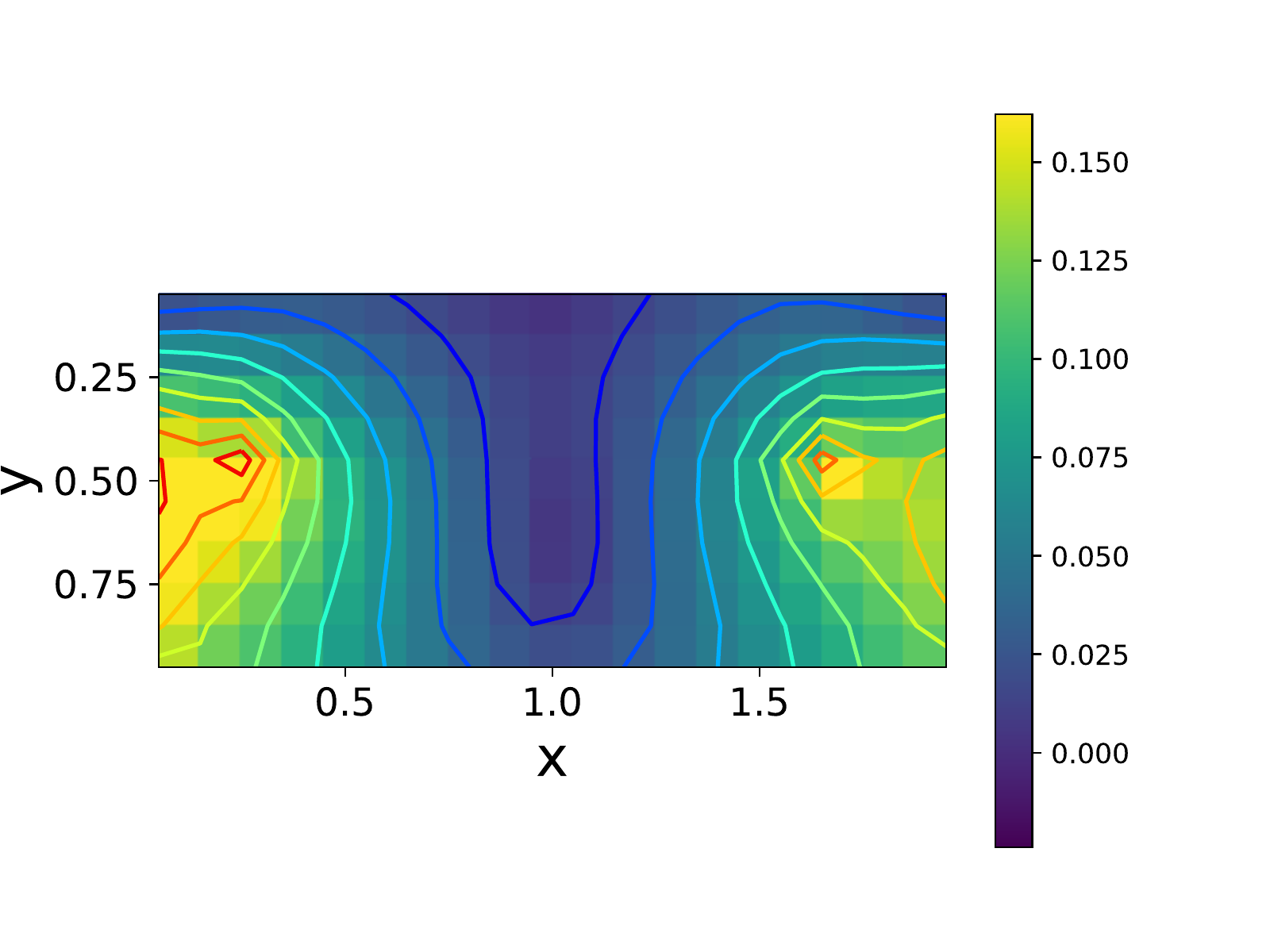}
	\caption{True Model}
\end{subfigure}
\caption{Von Mises strain distributions at terminal time. On the left and middle panels, the top plots show the distribution generated from the initial guess in both methods. The bottom plots show the distribution generated from the calibrated models.}
\label{fig:spaces}
\end{figure}

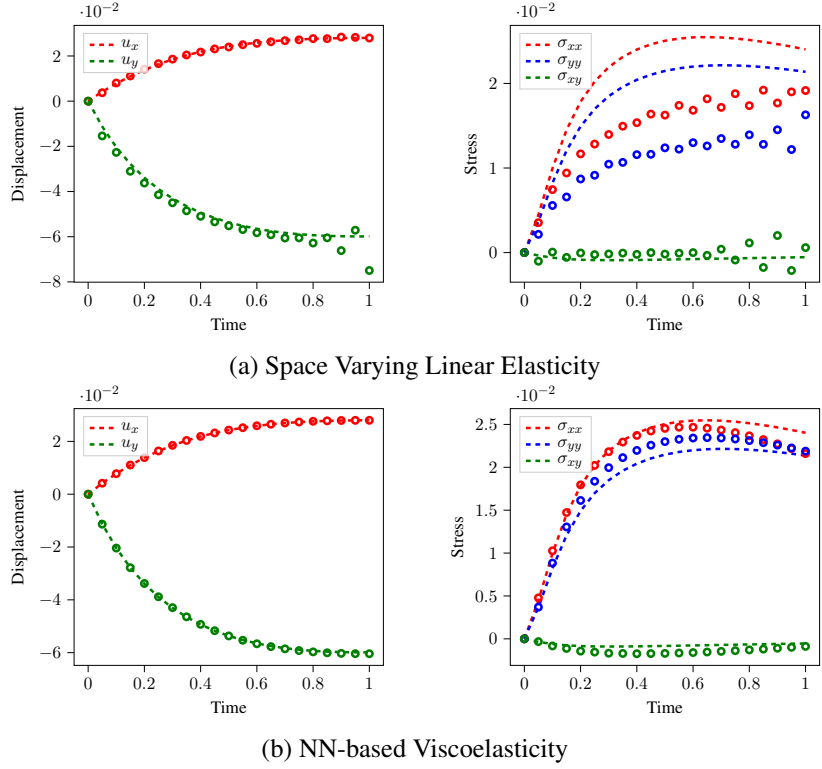
\begin{figure}[htpb]
\begin{subfigure}[t]{1.0\textwidth}
\centering
	\scalebox{0.6}{
\begin{tikzpicture}

\begin{groupplot}[group style={group size=2 by 1, horizontal sep=80pt}]
\nextgroupplot[
legend cell align={left},
legend style={fill opacity=0.8, draw opacity=1, text opacity=1, at={(0.03,0.97)}, anchor=north west, draw=white!80.0!black},
tick align=outside,
tick pos=left,
x grid style={white!69.01960784313725!black},
xlabel={Time},
xmin=-0.05, xmax=1.05,
xtick style={color=black},
y grid style={white!69.01960784313725!black},
ylabel={Displacement},
ymin=-0.0801213883868493, ymax=0.03353125177083,
ytick style={color=black}
]
\addplot [ultra thick, red, dashed]
table {%
0 0
0.05 0.00417823202893617
0.1 0.00778262106008245
0.15 0.0109902002809395
0.2 0.0138348578022242
0.25 0.016325252799982
0.3 0.0184773122694861
0.35 0.0203163667384792
0.4 0.0218729748660487
0.45 0.0231790795743023
0.5 0.0242655768798464
0.55 0.0251610446857104
0.6 0.0258912189143624
0.65 0.0264788995264547
0.7 0.0269440818730851
0.75 0.027304192081293
0.8 0.0275743585904904
0.85 0.027767684011373
0.9 0.0278954998906127
0.95 0.0279675971497867
1 0.0279924303683926
};
\addlegendentry{$u_x$}
\addplot [ultra thick, red, mark=o, mark size=2, mark options={solid}, only marks, forget plot]
table {%
0 0
0.05 0.00373076449841727
0.1 0.0080297901022909
0.15 0.0111442432828665
0.2 0.0142001093558407
0.25 0.0166036431165658
0.3 0.0185974401554295
0.35 0.0204430708747142
0.4 0.0217328145867045
0.45 0.0231245870986168
0.5 0.0239743091723988
0.55 0.0249934370471451
0.6 0.0255866419706541
0.65 0.0262836922835764
0.7 0.0267707931053311
0.75 0.0271394102383636
0.8 0.0276865875162755
0.85 0.0276643023327844
0.9 0.0283652226727536
0.95 0.0282361361138348
1 0.0279517446686601
};
\addplot [ultra thick, green!50.0!black, dashed]
table {%
0 0
0.05 -0.0112353840426635
0.1 -0.0203866614938739
0.15 -0.0278190583325006
0.2 -0.0338985949473629
0.25 -0.0389086376422945
0.3 -0.0430571854905126
0.35 -0.0464986008877637
0.4 -0.0493509736268258
0.45 -0.0517075403755356
0.5 -0.0536438695792025
0.55 -0.0552224266177515
0.6 -0.0564955753518074
0.65 -0.057507641448364
0.7 -0.0582963977994067
0.75 -0.0588941801703612
0.8 -0.059328755565301
0.85 -0.0596240174955178
0.9 -0.0598005547270633
0.95 -0.0598761239604021
1 -0.0598660472222394
};
\addlegendentry{$u_y$}
\addplot [ultra thick, green!50.0!black, mark=o, mark size=2, mark options={solid}, only marks, forget plot]
table {%
0 0
0.05 -0.015404565364234
0.1 -0.0226842710563333
0.15 -0.0310286969263911
0.2 -0.0362497597747967
0.25 -0.0414563591188589
0.3 -0.044984277282659
0.35 -0.0485309048908543
0.4 -0.0509669547866934
0.45 -0.053428436503345
0.5 -0.0551889448801136
0.55 -0.0568359640684951
0.6 -0.0582411022266569
0.65 -0.0591495219984261
0.7 -0.0605833070401631
0.75 -0.0604949990243936
0.8 -0.0627866314665423
0.85 -0.060484332726404
0.9 -0.0661458902114231
0.95 -0.0571223196105983
1 -0.0749553592887729
};

\nextgroupplot[
legend cell align={left},
legend style={fill opacity=0.8, draw opacity=1, text opacity=1, at={(0.03,0.97)}, anchor=north west, draw=white!80.0!black},
tick align=outside,
tick pos=left,
x grid style={white!69.01960784313725!black},
xlabel={Time},
xmin=-0.05, xmax=1.05,
xtick style={color=black},
y grid style={white!69.01960784313725!black},
ylabel={Stress},
ymin=-0.00349710545866601, ymax=0.0268612847654836,
ytick style={color=black}
]
\addplot [ultra thick, red, dashed]
table {%
0 0
0.05 0.00447661446243479
0.1 0.00995445394259398
0.15 0.0144985721497447
0.2 0.0178400231864664
0.25 0.020222044687035
0.3 0.021919900353217
0.35 0.0231368568056531
0.4 0.0240084575123435
0.45 0.0246237466629114
0.5 0.0250424421115917
0.55 0.0253060845464823
0.6 0.0254448046948328
0.65 0.0254813579371131
0.7 0.0254335414106765
0.75 0.025315669278529
0.8 0.0251394995146542
0.85 0.0249148380044124
0.9 0.0246499493322016
0.95 0.0243518487694495
1 0.0240265188398888
};
\addlegendentry{$\sigma_{xx}$}
\addplot [ultra thick, blue, dashed]
table {%
0 0
0.05 0.00366525528775264
0.1 0.00823030315266473
0.15 0.0120619732394573
0.2 0.014917007775115
0.25 0.0169864713315238
0.3 0.0184931403130063
0.35 0.0196025478115995
0.4 0.0204252351010343
0.45 0.0210337905489412
0.5 0.0214768419642486
0.55 0.0217881819358846
0.6 0.0219923140599285
0.65 0.0221077583012859
0.7 0.0221490277764376
0.75 0.0221278317130024
0.8 0.0220538280985743
0.85 0.0219351118584006
0.9 0.0217785450340394
0.95 0.0215899902399097
1 0.0213744830148816
};
\addlegendentry{$\sigma_{yy}$}
\addplot [ultra thick, green!50.0!black, dashed]
table {%
0 0
0.05 -0.00034700589725159
0.1 -0.000601779768426984
0.15 -0.000758114720409178
0.2 -0.000844237101152967
0.25 -0.000884581361412765
0.3 -0.000896439306598954
0.35 -0.000890991440300645
0.4 -0.00087514237580473
0.45 -0.000853077329151867
0.5 -0.000827330542943508
0.55 -0.00079945398030136
0.6 -0.000770417829386695
0.65 -0.000740845498653251
0.7 -0.000711150453474966
0.75 -0.000681615964695061
0.8 -0.000652441948690276
0.85 -0.000623772881570668
0.9 -0.000595714801190315
0.95 -0.000568345972591857
1 -0.000541723829043222
};
\addlegendentry{$\sigma_{xy}$}
\addplot [ultra thick, red, mark=o, mark size=2, mark options={solid}, only marks, forget plot]
table {%
0 0
0.05 0.00354142356020415
0.1 0.00744664854185497
0.15 0.00941010572546121
0.2 0.0116698288256716
0.25 0.012832129377
0.3 0.0139633684134598
0.35 0.0149441338520072
0.4 0.0153577040738771
0.45 0.0163725862184822
0.5 0.0162534692302899
0.55 0.017401445143293
0.6 0.0168256433997746
0.65 0.0181799969946985
0.7 0.0171702022796343
0.75 0.0187884543292734
0.8 0.0173767865829721
0.85 0.0192033600879651
0.9 0.0176851706800677
0.95 0.0190020356144554
1 0.0191641192380631
};
\addplot [ultra thick, blue, mark=o, mark size=2, mark options={solid}, only marks, forget plot]
table {%
0 0
0.05 0.00215541538050174
0.1 0.00556883947095832
0.15 0.0065749375632695
0.2 0.00870411434043962
0.25 0.009143581025216
0.3 0.0104499661793784
0.35 0.0106655528309637
0.4 0.0115801654178885
0.45 0.0116208797035435
0.5 0.0123864915328336
0.55 0.0122305351107999
0.6 0.0129960875750752
0.65 0.0126076591891554
0.7 0.0134822836145408
0.75 0.0128082858171121
0.8 0.0139195138725963
0.85 0.0128001563414567
0.9 0.0145204055807971
0.95 0.0121857982216425
1 0.0162849850656455
};
\addplot [ultra thick, green!50.0!black, mark=o, mark size=2, mark options={solid}, only marks, forget plot]
table {%
0 0
0.05 -0.00102201430960144
0.1 6.08602451273593e-05
0.15 -0.000568289112270222
0.2 -4.53076077411481e-05
0.25 -0.000236246502091495
0.3 -0.000159631250861562
0.35 -6.03027897689572e-05
0.4 -0.000210642878066931
0.45 -3.12824619895405e-06
0.5 -0.000173316807409952
0.55 -7.18197437360273e-05
0.6 -3.39581126952357e-07
0.65 -0.000326548301143477
0.7 0.000401847154464331
0.75 -0.000876691757012589
0.8 0.00114322550018289
0.85 -0.00175902129813426
0.9 0.00201959063886623
0.95 -0.00211717863029558
1 0.0005902920121644
};
\end{groupplot}

\end{tikzpicture}}
	\caption{Space Varying Linear Elasticity}
\end{subfigure}
	\begin{subfigure}[t]{1.0\textwidth}
	\centering
	\scalebox{0.6}{
\begin{tikzpicture}

\begin{groupplot}[group style={group size=2 by 1, horizontal sep=80pt}]
\nextgroupplot[
legend cell align={left},
legend style={fill opacity=0.8, draw opacity=1, text opacity=1, at={(0.03,0.97)}, anchor=north west, draw=white!80.0!black},
tick align=outside,
tick pos=left,
x grid style={white!69.01960784313725!black},
xlabel={Time},
xmin=-0.05, xmax=1.05,
xtick style={color=black},
y grid style={white!69.01960784313725!black},
ylabel={Displacement},
ymin=-0.0648210619648071, ymax=0.0324121204794973,
ytick style={color=black}
]
\addplot [ultra thick, red, dashed]
table {%
0 0
0.05 0.00417823202893617
0.1 0.00778262106008245
0.15 0.0109902002809395
0.2 0.0138348578022242
0.25 0.016325252799982
0.3 0.0184773122694861
0.35 0.0203163667384792
0.4 0.0218729748660487
0.45 0.0231790795743023
0.5 0.0242655768798464
0.55 0.0251610446857104
0.6 0.0258912189143624
0.65 0.0264788995264547
0.7 0.0269440818730851
0.75 0.027304192081293
0.8 0.0275743585904904
0.85 0.027767684011373
0.9 0.0278954998906127
0.95 0.0279675971497867
1 0.0279924303683926
};
\addlegendentry{$u_x$}
\addplot [ultra thick, red, mark=o, mark size=2, mark options={solid}, only marks, forget plot]
table {%
0 0
0.05 0.00417042839967369
0.1 0.00777529473567472
0.15 0.011048339283332
0.2 0.0138665978666862
0.25 0.0163751916069114
0.3 0.0185075483125997
0.35 0.0203477364616554
0.4 0.0218940482399183
0.45 0.0231999583378458
0.5 0.0242837399964738
0.55 0.0251813339972308
0.6 0.0259131052430044
0.65 0.0265040013022705
0.7 0.0269713822298395
0.75 0.027332922226287
0.8 0.0276021992168473
0.85 0.0277918954691408
0.9 0.0279122566901214
0.95 0.027972445055836
1 0.0279800530890972
};
\addplot [ultra thick, green!50.0!black, dashed]
table {%
0 0
0.05 -0.0112353840426635
0.1 -0.0203866614938739
0.15 -0.0278190583325006
0.2 -0.0338985949473629
0.25 -0.0389086376422945
0.3 -0.0430571854905126
0.35 -0.0464986008877637
0.4 -0.0493509736268258
0.45 -0.0517075403755356
0.5 -0.0536438695792025
0.55 -0.0552224266177515
0.6 -0.0564955753518074
0.65 -0.057507641448364
0.7 -0.0582963977994067
0.75 -0.0588941801703612
0.8 -0.059328755565301
0.85 -0.0596240174955178
0.9 -0.0598005547270633
0.95 -0.0598761239604021
1 -0.0598660472222394
};
\addlegendentry{$u_y$}
\addplot [ultra thick, green!50.0!black, mark=o, mark size=2, mark options={solid}, only marks, forget plot]
table {%
0 0
0.05 -0.011353426658663
0.1 -0.0203913712688422
0.15 -0.0278300961449616
0.2 -0.0338140284216805
0.25 -0.0388361753636636
0.3 -0.0429531035611755
0.35 -0.0464237972709672
0.4 -0.0492899567903692
0.45 -0.051691996191205
0.5 -0.0536695051436681
0.55 -0.0553048870925602
0.6 -0.0566334826371014
0.65 -0.0577073490341997
0.7 -0.058555787291752
0.75 -0.0592135607636411
0.8 -0.0597038133508143
0.85 -0.0600506302436655
0.9 -0.0602717478252616
0.95 -0.0603842203815988
1 -0.0604013718537024
};

\nextgroupplot[
legend cell align={left},
legend style={fill opacity=0.8, draw opacity=1, text opacity=1, at={(0.03,0.97)}, anchor=north west, draw=white!80.0!black},
tick align=outside,
tick pos=left,
x grid style={white!69.01960784313725!black},
xlabel={Time},
xmin=-0.05, xmax=1.05,
xtick style={color=black},
y grid style={white!69.01960784313725!black},
ylabel={Stress},
ymin=-0.00309393197445362, ymax=0.0268420860281401,
ytick style={color=black}
]
\addplot [ultra thick, red, dashed]
table {%
0 0
0.05 0.00447661446243479
0.1 0.00995445394259398
0.15 0.0144985721497447
0.2 0.0178400231864664
0.25 0.020222044687035
0.3 0.021919900353217
0.35 0.0231368568056531
0.4 0.0240084575123435
0.45 0.0246237466629114
0.5 0.0250424421115917
0.55 0.0253060845464823
0.6 0.0254448046948328
0.65 0.0254813579371131
0.7 0.0254335414106765
0.75 0.025315669278529
0.8 0.0251394995146542
0.85 0.0249148380044124
0.9 0.0246499493322016
0.95 0.0243518487694495
1 0.0240265188398888
};
\addlegendentry{$\sigma_{xx}$}
\addplot [ultra thick, blue, dashed]
table {%
0 0
0.05 0.00366525528775264
0.1 0.00823030315266473
0.15 0.0120619732394573
0.2 0.014917007775115
0.25 0.0169864713315238
0.3 0.0184931403130063
0.35 0.0196025478115995
0.4 0.0204252351010343
0.45 0.0210337905489412
0.5 0.0214768419642486
0.55 0.0217881819358846
0.6 0.0219923140599285
0.65 0.0221077583012859
0.7 0.0221490277764376
0.75 0.0221278317130024
0.8 0.0220538280985743
0.85 0.0219351118584006
0.9 0.0217785450340394
0.95 0.0215899902399097
1 0.0213744830148816
};
\addlegendentry{$\sigma_{yy}$}
\addplot [ultra thick, green!50.0!black, dashed]
table {%
0 0
0.05 -0.00034700589725159
0.1 -0.000601779768426984
0.15 -0.000758114720409178
0.2 -0.000844237101152967
0.25 -0.000884581361412765
0.3 -0.000896439306598954
0.35 -0.000890991440300645
0.4 -0.00087514237580473
0.45 -0.000853077329151867
0.5 -0.000827330542943508
0.55 -0.00079945398030136
0.6 -0.000770417829386695
0.65 -0.000740845498653251
0.7 -0.000711150453474966
0.75 -0.000681615964695061
0.8 -0.000652441948690276
0.85 -0.000623772881570668
0.9 -0.000595714801190315
0.95 -0.000568345972591857
1 -0.000541723829043222
};
\addlegendentry{$\sigma_{xy}$}
\addplot [ultra thick, red, mark=o, mark size=2, mark options={solid}, only marks, forget plot]
table {%
0 0
0.05 0.00478077729407618
0.1 0.0102672187120263
0.15 0.0147270333216637
0.2 0.0179530719157509
0.25 0.0202215798136911
0.3 0.0218208886948118
0.35 0.0229414599162989
0.4 0.0237167094617784
0.45 0.0242289840701366
0.5 0.0245353619900243
0.55 0.0246749298515982
0.6 0.0246762087883216
0.65 0.0245610114050725
0.7 0.0243462357681521
0.75 0.0240457471874207
0.8 0.0236707773860069
0.85 0.0232308534614964
0.9 0.0227339023429152
0.95 0.0221867310454986
1 0.0215950853000314
};
\addplot [ultra thick, blue, mark=o, mark size=2, mark options={solid}, only marks, forget plot]
table {%
0 0
0.05 0.00369720569940358
0.1 0.00883727686295203
0.15 0.0130462861101968
0.2 0.0161349643867507
0.25 0.0183715851566468
0.3 0.0199597679478174
0.35 0.0211283053847298
0.4 0.0219667479346413
0.45 0.0225760068959853
0.5 0.022994020150851
0.55 0.0232665171458019
0.6 0.0234137148421192
0.65 0.0234591681175742
0.7 0.0234159911110299
0.75 0.0232980616747294
0.8 0.0231147349885039
0.85 0.0228750687688884
0.9 0.0225859581325922
0.95 0.0222537018531817
1 0.0218834232547522
};
\addplot [ultra thick, green!50.0!black, mark=o, mark size=2, mark options={solid}, only marks, forget plot]
table {%
0 0
0.05 -0.000321031352914378
0.1 -0.000827844988554562
0.15 -0.00112037678822032
0.2 -0.00143109124926108
0.25 -0.00154672154307872
0.3 -0.00167357180532328
0.35 -0.00170142972770634
0.4 -0.00173320388342663
0.45 -0.0017162163970575
0.5 -0.00169648606315658
0.55 -0.00165100849249778
0.6 -0.00159977017896252
0.65 -0.00153350034439356
0.7 -0.00146089927922692
0.75 -0.00137858805798475
0.8 -0.00129047251611521
0.85 -0.00119554565459243
0.9 -0.00109567604315507
0.95 -0.000990742687988959
1 -0.000881728286810345
};
\end{groupplot}

\end{tikzpicture}}
	\caption{NN-based Viscoelasticity}
\end{subfigure}
  \caption{Displacement and stress tensors of the left top point. The dashed lines are true values, and the dots are reproduced values using the calibrated models.}
  \label{fig:spacecurve}
\end{figure}

We generate five sets of training data with different source function scales. The training data are used to calibrate the linear elastic matrices or training the neural network. Then we test on the test data, whose source function scale is within the ranges of training set force scales.  
The results are shown in \Cref{fig:spacedisp,fig:spaces,fig:spacecurve}. We can see that the neural network outperforms the space varying linear elastic models. It is interesting to see from \Cref{fig:spacecurve} that the space varying linear elasticity model actually reproduces the  observed data (red dots and lines on the left) very well. However, its prediction of vertical displacement and stress tensors indicates that the space varying linear elasticity model overfits. The space varying linear elasticity model has too much freedom to select specific linear elastic matrices for each Gauss quadrature points to fit the observed data.

\subsection{Learning Constitutive Relations in Coupled Geomechanics and Two-Phase Flow Equations}
\label{sect:5}

\begin{figure}[htpb]
\centering
  \includegraphics[width=0.45\textwidth]{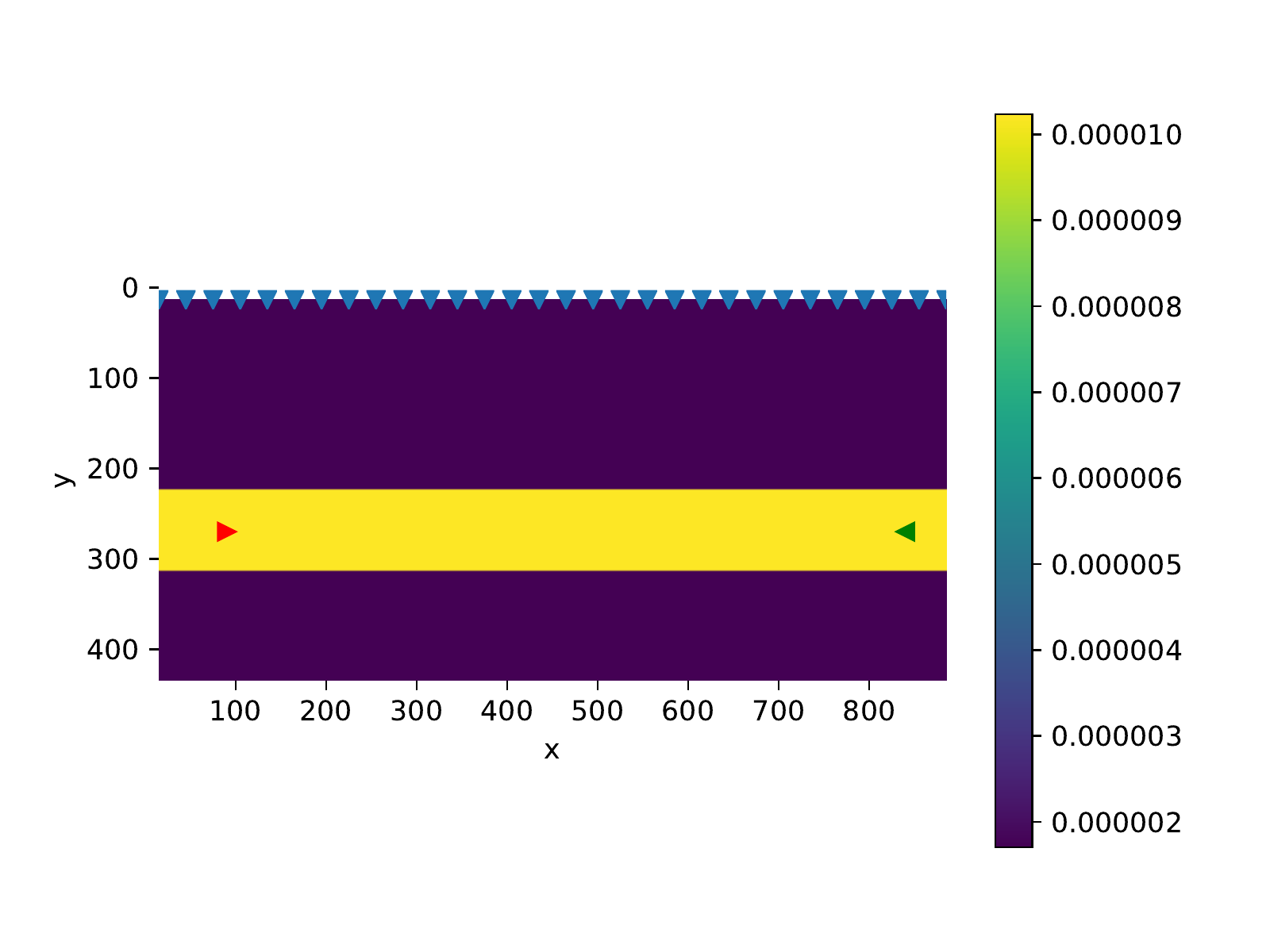}
    \includegraphics[width=0.45\textwidth]{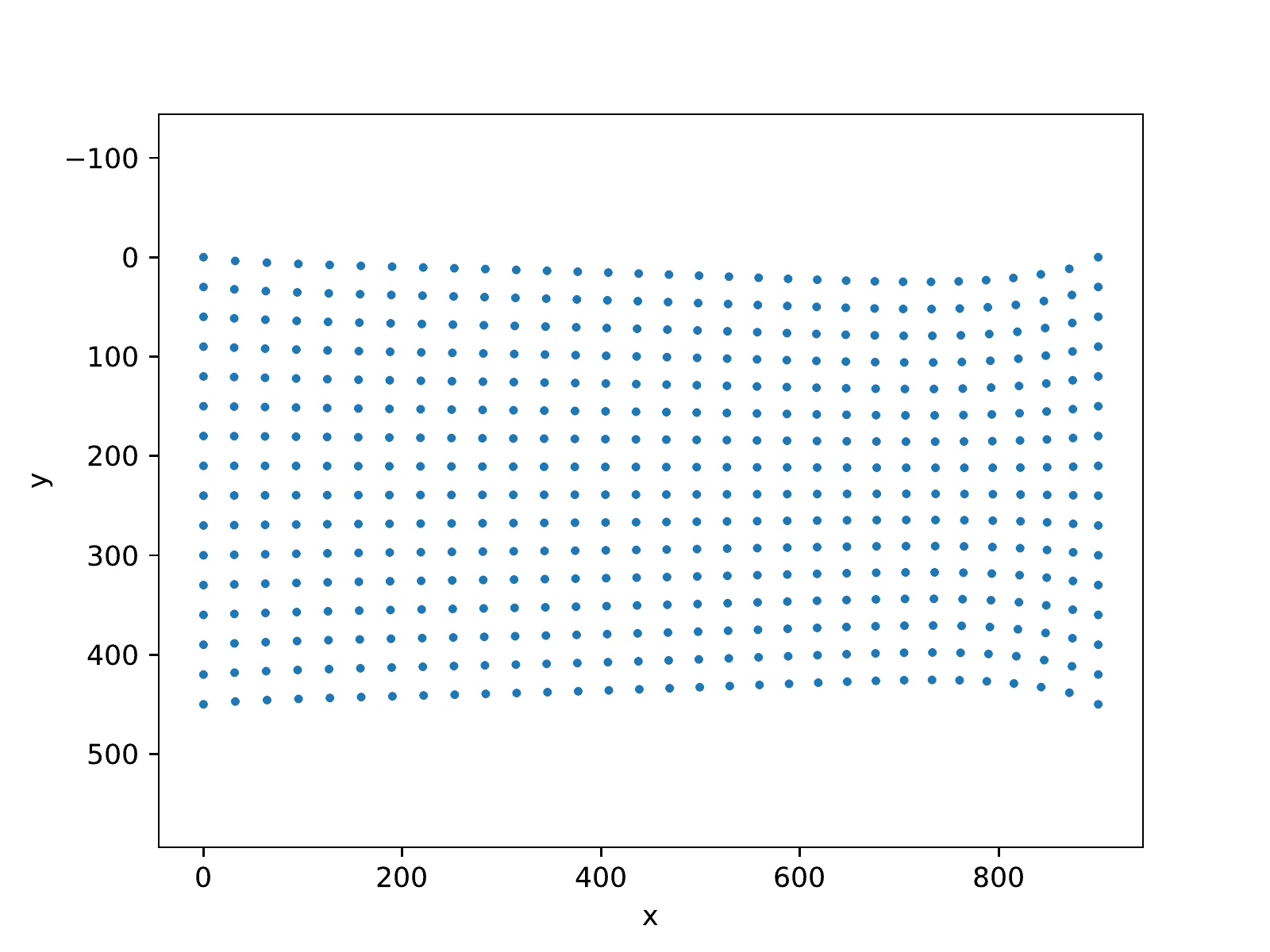}
  \caption{Left: the distribution of the relative permeability $K$ in the numerical simulation. Right: the displacement at the terminal time. We magnify the displacements by 10 times for clarity.}
  \label{fig:setting}
\end{figure}

\begin{figure}[htpb]
\centering
  \includegraphics[width=0.33\textwidth]{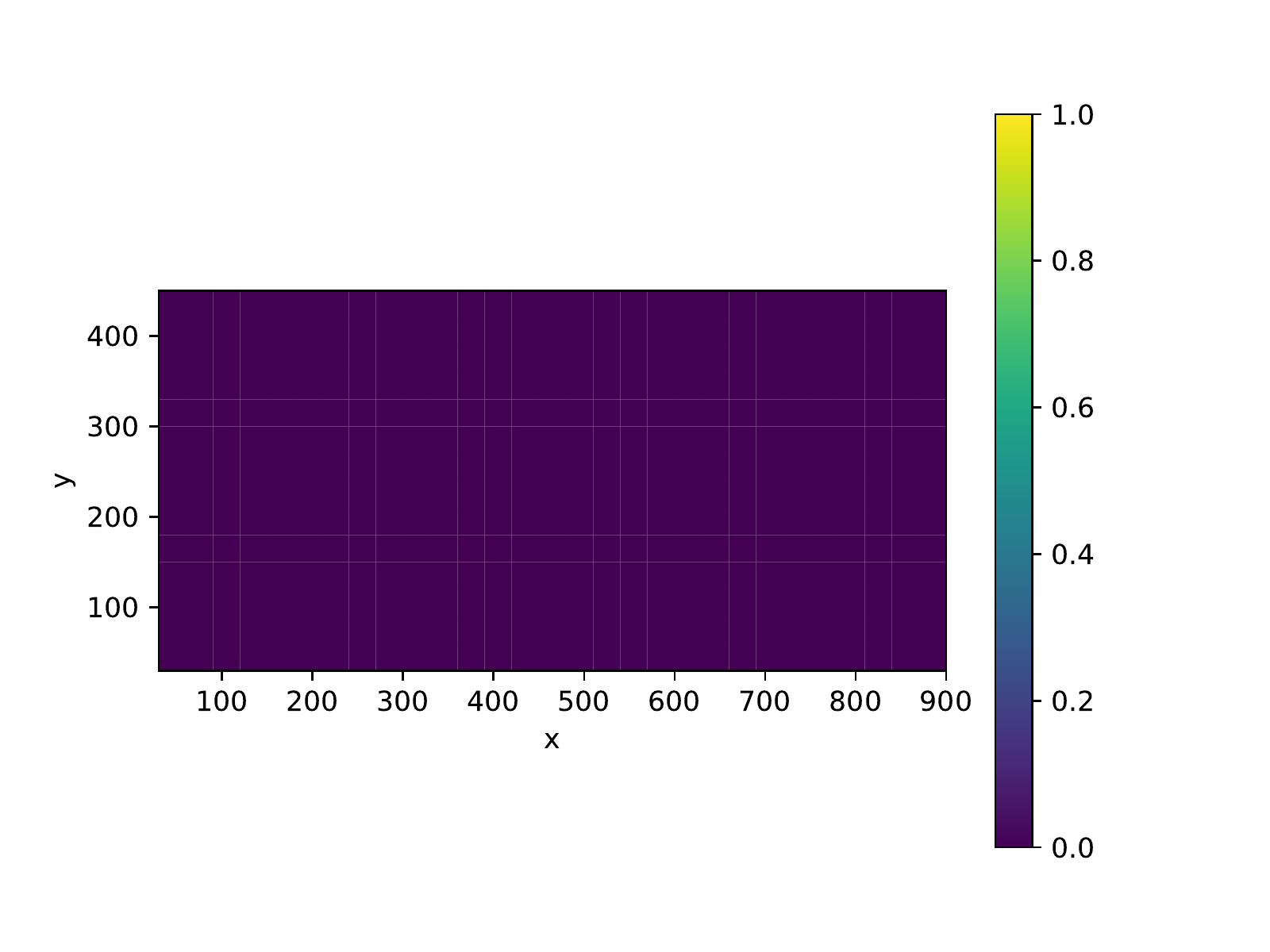}~
    \includegraphics[width=0.33\textwidth]{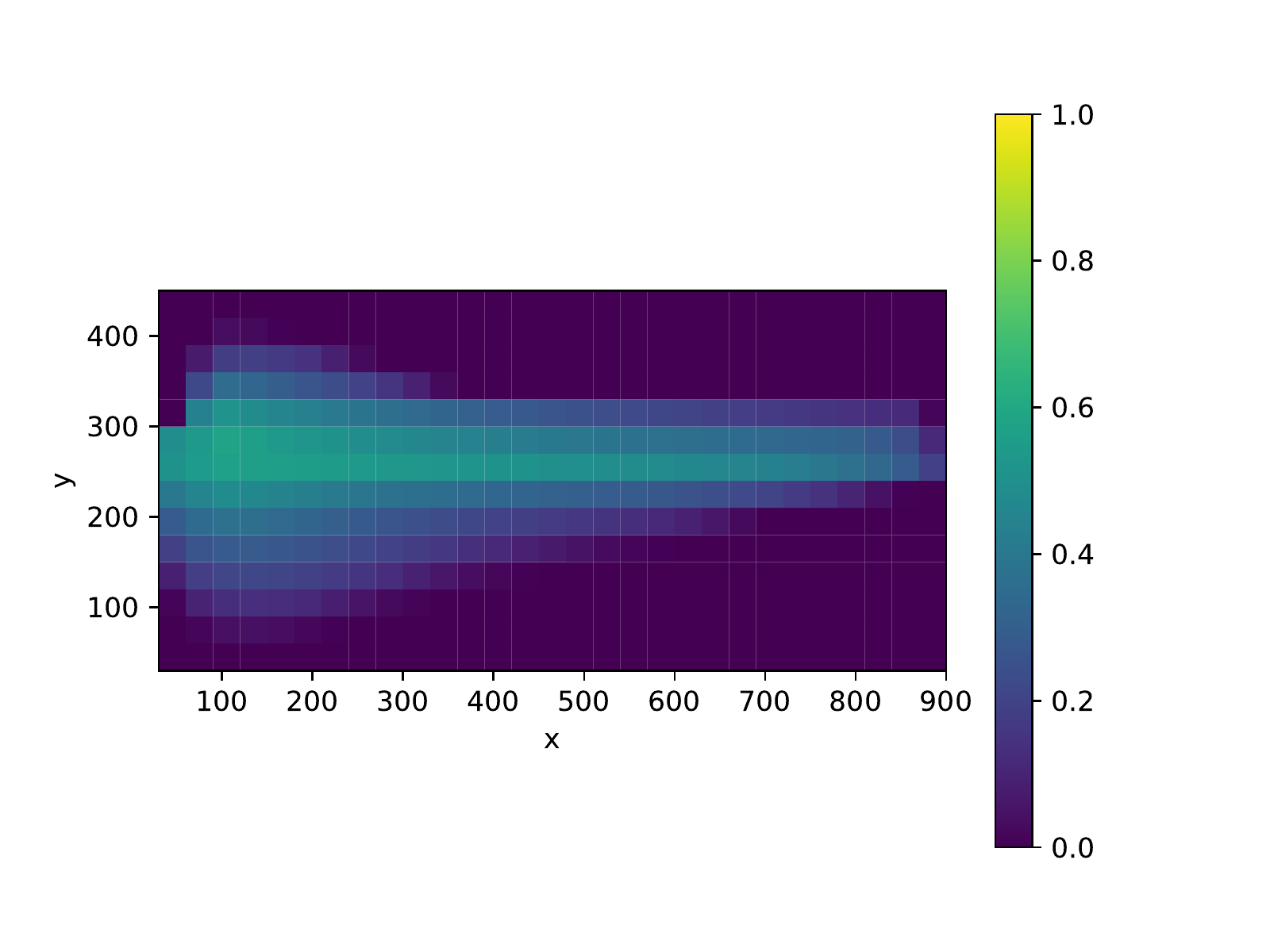}~
      \includegraphics[width=0.33\textwidth]{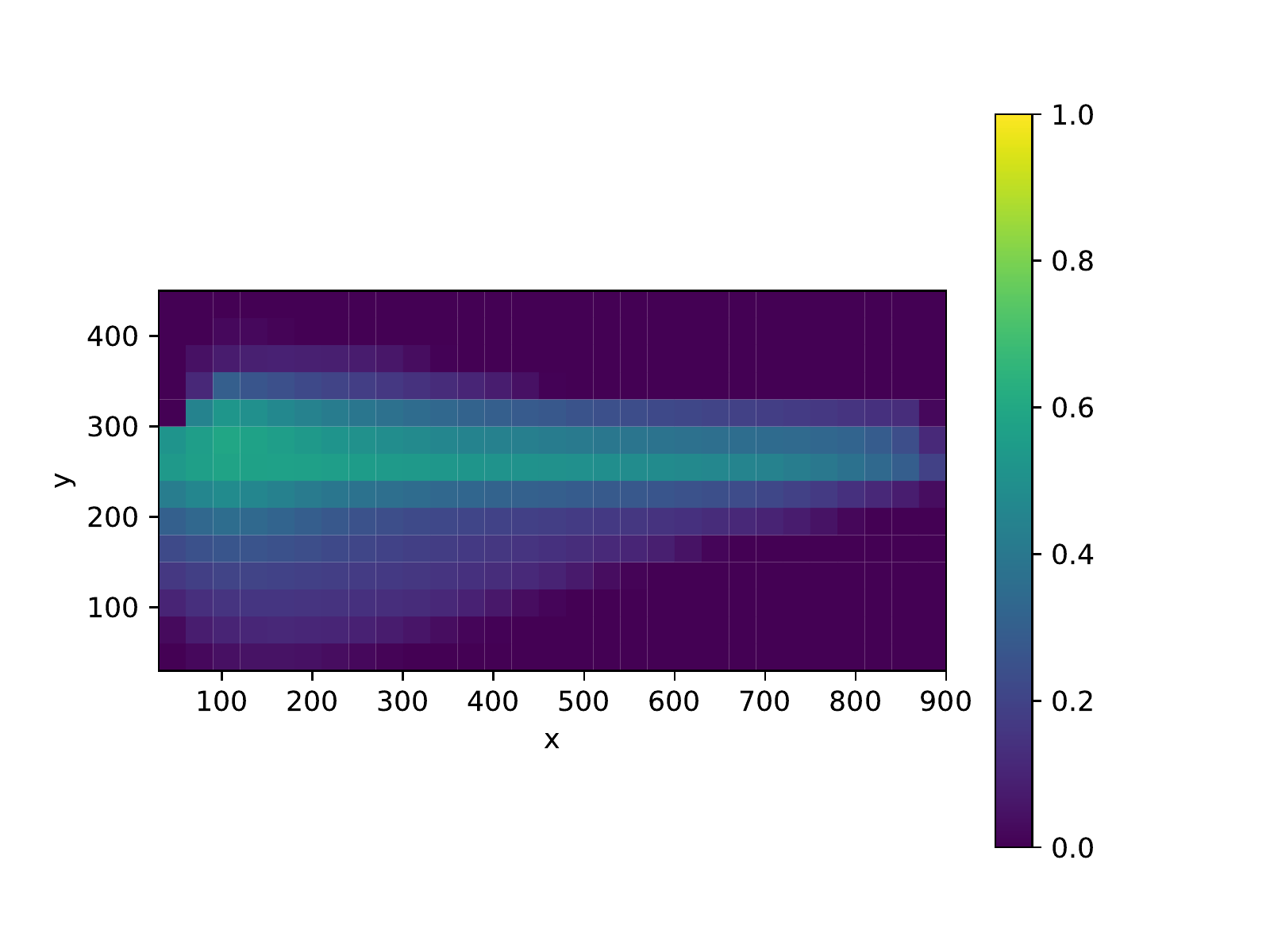}
  \caption{Evolution of the saturation $S_2$.}
  \label{fig:two_s}
\end{figure}

We consider the inverse computation in the coupled geomechanics and two-phase flow equations. The governing equations are given in \Cref{sect:twophaseflow}, and the geometric setting is given by \Cref{fig:ip}. We consider two models for the constitutive relations: the linear elasticity \Cref{equ:linear2} and the viscoelasticity \Cref{equ:maxwell2}. In the linear elasticity model, the unknown is the linear elasticity matrix, and the error is reported by \Cref{equ:errorH}. In the viscoelasticity model, the unknowns are $\mu$, $\lambda$, and $\eta$. The error is reported by 
\begin{equation}\label{equ:errorMu}
	E_\mu = \frac{|\mu-\mu^*|}{|\mu^*|}\qquad 	E_\lambda = \frac{|\lambda-\lambda^*|}{|\lambda^*|}\qquad 	E_\eta = \frac{|\eta-\eta^*|}{|\eta^*|}
\end{equation} 
where $\mu^*$, $\lambda^*$, and $\eta^*$ are the true values. \Cref{fig:setting} shows the distribution of the relative permeability $K$ in the numerical simulation and the displacement at the terminal time. \Cref{fig:two_s} shows the evolution of the saturation $S_2$. Since we have included the gravity constant, the evolution of the saturation is not symmetric vertically.

\begin{figure}[htpb]
\centering
	\scalebox{0.6}{
\begin{tikzpicture}

\definecolor{color0}{rgb}{0.12156862745098,0.466666666666667,0.705882352941177}

\begin{axis}[
log basis y={10},
tick align=outside,
tick pos=left,
grid=both,
x grid style={white!69.0196078431373!black},
xlabel={\(\displaystyle \sigma_{\mathrm{noise}}\)},
xtick={0,0.02,0.04,0.06,0.08,0.10},
 xticklabels={0,0.02,0.04,0.06,0.08,0.10},
xmin=-0.005, xmax=0.105,
xtick style={color=black},
y grid style={white!69.0196078431373!black},
ylabel={Error},
ymin=0.000112056504174044, ymax=0.372659694423655,
ymode=log,
ytick style={color=black}
]
\path [fill=color0, fill opacity=0.5]
(axis cs:0,0.000162002722491279)
--(axis cs:0,0.000162002722491279)
--(axis cs:0.001,0.000666525323816818)
--(axis cs:0.003,0.00327374700424619)
--(axis cs:0.009,0.00663080334776873)
--(axis cs:0.01,0.0058723571280109)
--(axis cs:0.02,0.00987687683308293)
--(axis cs:0.03,0.00848022463674844)
--(axis cs:0.04,0.0269491907634617)
--(axis cs:0.05,0.0189232418155159)
--(axis cs:0.06,0.0504287889638052)
--(axis cs:0.1,0.118689090109433)
--(axis cs:0.1,0.257766918737617)
--(axis cs:0.1,0.257766918737617)
--(axis cs:0.06,0.121020746960289)
--(axis cs:0.05,0.0619324540234882)
--(axis cs:0.04,0.0631443442851182)
--(axis cs:0.03,0.0452135493380398)
--(axis cs:0.02,0.025960122785019)
--(axis cs:0.01,0.0133011756871595)
--(axis cs:0.009,0.0115519042688197)
--(axis cs:0.003,0.00532186375320271)
--(axis cs:0.001,0.00140019729920141)
--(axis cs:0,0.000162002722491279)
--cycle;

\addplot [very thick, color0]
table {%
0 0.000162002722491279
0.001 0.00103336131150911
0.003 0.00429780537872445
0.009 0.00909135380829424
0.01 0.00958676640758522
0.02 0.017918499809051
0.03 0.0268468869873941
0.04 0.04504676752429
0.05 0.040427847919502
0.06 0.0857247679620471
0.1 0.188228004423525
};
\end{axis}

\end{tikzpicture}}~
	\scalebox{0.6}{
\begin{tikzpicture}

\definecolor{color0}{rgb}{0.12156862745098,0.466666666666667,0.705882352941177}
\definecolor{color1}{rgb}{1,0.498039215686275,0.0549019607843137}
\definecolor{color2}{rgb}{0.172549019607843,0.627450980392157,0.172549019607843}

\begin{axis}[
legend cell align={left},
legend style={fill opacity=0.8, draw opacity=1, text opacity=1, at={(0.97,0.03)}, anchor=south east, draw=white!80!black},
log basis y={10},
tick align=outside,
tick pos=left,
grid=both,
x grid style={white!69.0196078431373!black},
xlabel={$\sigma_{\mathrm{noise}}$},
xmin=-0.005, xmax=0.105,
xtick style={color=black},
xtick={0,0.02,0.04,0.06,0.08,0.10},
 xticklabels={0,0.02,0.04,0.06,0.08,0.10},
y grid style={white!69.0196078431373!black},
ylabel={Error},
ymin=7.0174147964051e-08, ymax=0.0739802612144677,
ymode=log,
ytick style={color=black},
ytick={1e-09,1e-08,1e-07,1e-06,1e-05,0.0001,0.001,0.01,0.1,1},
yticklabels={\(\displaystyle {10^{-9}}\),\(\displaystyle {10^{-8}}\),\(\displaystyle {10^{-7}}\),\(\displaystyle {10^{-6}}\),\(\displaystyle {10^{-5}}\),\(\displaystyle {10^{-4}}\),\(\displaystyle {10^{-3}}\),\(\displaystyle {10^{-2}}\),\(\displaystyle {10^{-1}}\),\(\displaystyle {10^{0}}\)}
]
\path [fill=color0, fill opacity=0.5]
(axis cs:0,2.00916133239269e-06)
--(axis cs:0,2.00916133239269e-06)
--(axis cs:0.001,1.01820635994524e-05)
--(axis cs:0.003,0.000116457536626234)
--(axis cs:0.009,0.00271361833070112)
--(axis cs:0.01,0.00313397721469435)
--(axis cs:0.02,0.00409135073041029)
--(axis cs:0.03,0.00318722914936966)
--(axis cs:0.04,0.0124227025499684)
--(axis cs:0.05,0.00960966164045476)
--(axis cs:0.06,0.00730763494051498)
--(axis cs:0.1,0.00869756542524509)
--(axis cs:0.1,0.0318269820130215)
--(axis cs:0.1,0.0318269820130215)
--(axis cs:0.06,0.0260394656940426)
--(axis cs:0.05,0.0185728722282232)
--(axis cs:0.04,0.0238701509569159)
--(axis cs:0.03,0.00850923447210321)
--(axis cs:0.02,0.00840991909076533)
--(axis cs:0.01,0.00644516996011387)
--(axis cs:0.009,0.00458130869003062)
--(axis cs:0.003,0.000673010604573402)
--(axis cs:0.001,0.000235464678143063)
--(axis cs:0,2.00916133239269e-06)
--cycle;

\path [fill=color1, fill opacity=0.5]
(axis cs:0,1.31809614453997e-07)
--(axis cs:0,1.31809614453997e-07)
--(axis cs:0.001,8.97360657892922e-05)
--(axis cs:0.003,0.000145328610317809)
--(axis cs:0.009,0.00142857991126305)
--(axis cs:0.01,0.000508892759026035)
--(axis cs:0.02,0.000643091772964121)
--(axis cs:0.03,0.00186193908824688)
--(axis cs:0.04,0.00386781216559051)
--(axis cs:0.05,0.00109046744442326)
--(axis cs:0.06,0.00100423750167516)
--(axis cs:0.1,0.00908035317026491)
--(axis cs:0.1,0.018849159527074)
--(axis cs:0.1,0.018849159527074)
--(axis cs:0.06,0.0144073951366435)
--(axis cs:0.05,0.00886444458426886)
--(axis cs:0.04,0.0105490787842811)
--(axis cs:0.03,0.00656962185221309)
--(axis cs:0.02,0.00495242442959129)
--(axis cs:0.01,0.00275632446892263)
--(axis cs:0.009,0.00236693950718273)
--(axis cs:0.003,0.000348580115521588)
--(axis cs:0.001,0.000178359635822595)
--(axis cs:0,1.31809614453997e-07)
--cycle;

\path [fill=color2, fill opacity=0.5]
(axis cs:0,8.01619580788571e-07)
--(axis cs:0,8.0161958078857e-07)
--(axis cs:0.001,0.00011069784029991)
--(axis cs:0.003,0.000646617782696149)
--(axis cs:0.009,0.00243084149096726)
--(axis cs:0.01,0.00153880851158646)
--(axis cs:0.02,0.00216268126622033)
--(axis cs:0.03,0.00388001481826823)
--(axis cs:0.04,0.0140736747110811)
--(axis cs:0.05,0.00219006063024832)
--(axis cs:0.06,0.00487879523390922)
--(axis cs:0.1,0.0230764870748074)
--(axis cs:0.1,0.0393863666045021)
--(axis cs:0.1,0.0393863666045021)
--(axis cs:0.06,0.0290696355423287)
--(axis cs:0.05,0.0102321350569891)
--(axis cs:0.04,0.0197362545776002)
--(axis cs:0.03,0.0105385496431669)
--(axis cs:0.02,0.00857929483605271)
--(axis cs:0.01,0.00553313449632353)
--(axis cs:0.009,0.00543725795986733)
--(axis cs:0.003,0.000851467282266381)
--(axis cs:0.001,0.000217593397588624)
--(axis cs:0,8.01619580788571e-07)
--cycle;

\addplot [very thick, color0]
table {%
0 2.00916133239269e-06
0.001 0.000122823370871258
0.003 0.000394734070599818
0.009 0.00364746351036587
0.01 0.00478957358740411
0.02 0.00625063491058781
0.03 0.00584823181073644
0.04 0.0181464267534422
0.05 0.014091266934339
0.06 0.0166735503172788
0.1 0.0202622737191333
};
\addlegendentry{$E_\mu$}
\addplot [very thick, color1]
table {%
0 1.31809614453997e-07
0.001 0.000134047850805943
0.003 0.000246954362919699
0.009 0.00189775970922289
0.01 0.00163260861397433
0.02 0.00279775810127771
0.03 0.00421578047022999
0.04 0.00720844547493582
0.05 0.00497745601434606
0.06 0.00770581631915934
0.1 0.0139647563486694
};
\addlegendentry{$E_\lambda$}
\addplot [very thick, color2]
table {%
0 8.01619580788571e-07
0.001 0.000164145618944267
0.003 0.000749042532481265
0.009 0.00393404972541729
0.01 0.003535971503955
0.02 0.00537098805113652
0.03 0.00720928223071758
0.04 0.0169049646443407
0.05 0.00621109784361869
0.06 0.016974215388119
0.1 0.0312314268396548
};
\addlegendentry{$E_\eta$}
\end{axis}

\end{tikzpicture}}
	\caption{Error versus noise levels. The left panel shows the error \Cref{equ:errorH} in the linear elasticity case, the right panel shows the error \Cref{equ:errorMu} in the viscoelasticity case. The shaded area shows the confidence interval for one standard deviation.}
	\label{fig:twodisp}
\end{figure}
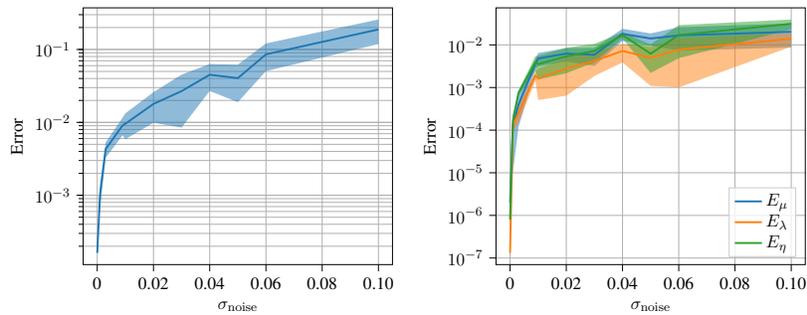

Results are reported in \Cref{fig:twodisp}, where we compare the error versus added noise levels according to \Cref{equ:noise} in the observation. We see that in general the error increases with noise levels, and the estimation is quite accurate and robust when $\sigma_{\mathrm{noise}}$ is small. 

\section{Limitations}

Despite many strengths of our proposed method, it is not without limitations, which we point out here. 

\begin{enumerate}
	\item \textit{Data sufficiency}. We demonstrated that very limited data (horizontal direction displacement data on the surface) are sufficient to quantify the viscoelastic behaviors of materials. However, there is no guarantee we can correctly estimate the constitutive relation. Many local minima may exist in the optimization problem, especially for the neural-network-based constitutive relations. This raises the question of data sufficiency: how much data do we need to quantify the viscoelasticity? More importantly, what properties can we infer from the given types of data? Those questions remain unanswered in our current work and are important for future investigations. 
	\item \textit{Numerical stability}. In general numerical stability depends on the physical parameter in the PDEs and the discretization scheme. We have used implicit schemes for numerical stability to the largest extent. However, this does not guarantee that our numerical scheme is always stable for any provided physical parameters. For example, without additional constraints, $\eta$, $\lambda$, and $\mu$ in \Cref{equ:maxwell} may become negative during the optimization process. The numerical scheme is usually developed based on certain assumptions of data and thus may break for nonphysical values. This problem is particularly salient in the context of neural networks, which may produce an ill-behaved function. Therefore, it is important that we incorporate certain constraints on the free optimization variables. For example, we restrict the linear elastic matrix to symmetric positive definite ones in \Cref{sect:5}. However, imposing proper constraints is challenging in many cases, especially for neural networks. 
	\item \textit{Predictive modeling}. One major motivation for inverse problems is predictive modeling. In predictive modeling, we first estimate the physical properties from historic/observation data (training data), and then use the estimated physical properties to predict future events with new inputs.  The implicit assumption is that the calibrated model should also work for the new inputs. However, this may not be the case if the new inputs deviate too much from the training data. For example, a learned neural-network-based constitutive relation is only valid for the range where the neural network is trained. Therefore, in principle, the calibrated model only works for new data that produces similar inputs to the neural network. To obtain a neural-network-based constitutive model that works for a large range of inputs, a large amount of training data is required. 
\end{enumerate}

\section{Conclusion}

We developed a unified framework for inverse modeling of viscoelasticity from limited displacement data. The inverse problem is formulated as a PDE-constrained optimization problem. Particularly, when we have a function inverse problem, the unknown function is approximated by a neural network and the weights and biases serve as free optimization variables. We leverage the automatic differentiation and physics constrained learning for computing the gradients of the objective function. Extensive numerical results show that our method is efficient, robust, and widely applicable. 

The developed method may be generalized to a broad range of constitutive models where only limited observation data are available. Since we use automatic differentiation to compute the gradients, we eliminate time-consuming effort in deriving and implementing gradients manually. The numerical implementation of physics constrained learning is also very flexible and can be integrated into an AD framework. We believe that the methodology of the current work is promising and can be extended to many other inverse modeling problem as well.

\appendix

\section{A Brief Introduction to ADCME.jl}

We have built our inverse modeling of viscoelasticity constitutive relations based on ADCME.jl, an automatic differentiation library for computational and mathematical engineering. ADCME.jl leverages TensorFlow for graph based optimization and automatic differentiation. It supports physics constrained optimization through custom operators. ADCME.jl is specifically designed for inverse modeling in scientific computing and is seamlessly integrated with the Julia language. Users can leverage the fast and efficient Julia programming environment instead of Python, which is the main programming languages of TensoFlow, for implementing numerical schemes. 

\section{Numerical Scheme for \Cref{sect:5}}

We use an iterative algorithm to solve the coupled equation; namely, we alternatively solve the mechanics equation and flow equation.

\paragraph{Fluid Equation} We define the fluid potential 
$$\Psi_i = P_i - \rho_i gy$$
and the capillary potential 
\begin{equation*}
  \Psi_c = \Psi_1 - \Psi_2 = P_1 -P_2 - (\rho_1-\rho_2)gy \approx - (\rho_1-\rho_2)gy
\end{equation*}
Here the capillary pressure $P_c = P_1-P_2$ is assumed to be small. We define mobilities 
\begin{equation}\label{equ:mt}
  m_i(S_i) = \frac{k_{ri}(S_i)}{\tilde\mu_i}, i=1,2\quad m_t = m_1 + m_2
\end{equation}
We have the following formula from Equations 3-4:
\begin{equation}\label{equ:Poisson}
  -\nabla\cdot (m_tK\nabla \Psi_2) = \nabla \cdot(m_1 K\nabla \Psi_c) - \frac{\partial \phi}{\partial t} + q_1 + q_2 
\end{equation}
Here $m_t=m_1+m_2$ is the \textit{total mobility}. We can solve for $\Psi_2$ using a Poisson solver. In our implementation, we use the AMGCL package \cite{Demidov2019} to solve the Poisson's equation using an algebraic multigrid method.

Next, we have from \Cref{equ:vs,equ:phi}.
\begin{equation}\label{equ:nonlinear}
  \phi\frac{\partial S_2}{\partial t} + S_2 \frac{\partial\phi}{\partial t} + \nabla \cdot (-K m_2 \nabla \Psi_2) = q_2 + q_1 \frac{m_2}{m_1} 
\end{equation}
Note we have an extra term $q_1 \frac{m_2}{m_1}$ to account for the injection of phase 2 flow. 

\Cref{equ:nonlinear} is a nonlinear equation in $S_2$ ($m_2$ is defined in terms of $S_2=1-S_1$) and requires a Newton-Raphson solver. 

\paragraph{Solid Equation} Upon solving the fluid equation, we obtain $S_1, S_2, \Psi_2$. We can use $\Psi_2$ to estimate the fluid pressure $P_1$ and $P_2$. Use \Cref{equ:mechanics}, we solve for $\mathbf{u}$ using
\begin{equation}\label{equ:solidmechanics}
  \int_\Omega \bsigma' :\delta \bepsilon \mathrm{d} x + \int_\Omega (S_1P_1+S_2P_2)\delta \epsilon_v \mathrm{d}x + \int_\Omega \mathbf{f}\delta u \mathrm{d} x = 0 \Leftrightarrow \int_\Omega \bsigma' :\delta \bepsilon \mathrm{d} x - \int_\Omega \Psi_2 \delta \epsilon_v \mathrm{d}x = 0
\end{equation}

\begin{figure}[hbt]
\centering
  \includegraphics[width=0.45\textwidth]{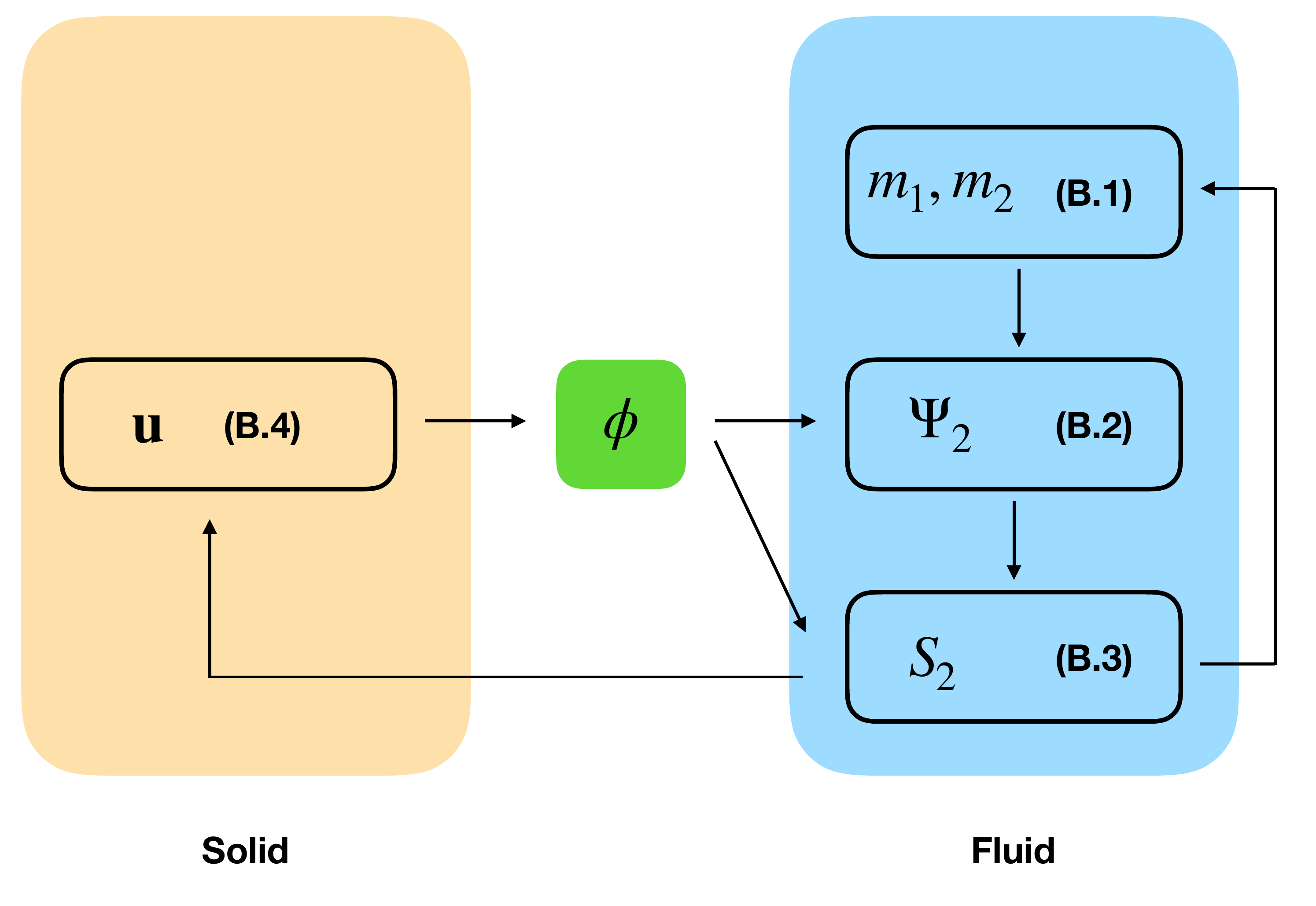}
  \caption{Flow chart of the numerical scheme for the coupled geomechanics and two-phase flow problem. }
  \label{fig:flowchart}
\end{figure}

\Cref{fig:flowchart} shows the flow chart of solving the coupled geomechanics and two-phase flow equations. 

\section{Additional Results for \Cref{sect:3}}\label{sect:additional}

We show additional results for \Cref{sect:3} in this section. 

\begin{figure}[htpb]
\begin{subfigure}[t]{0.48\textwidth}
\centering
	\includegraphics[width=1.0\textwidth]{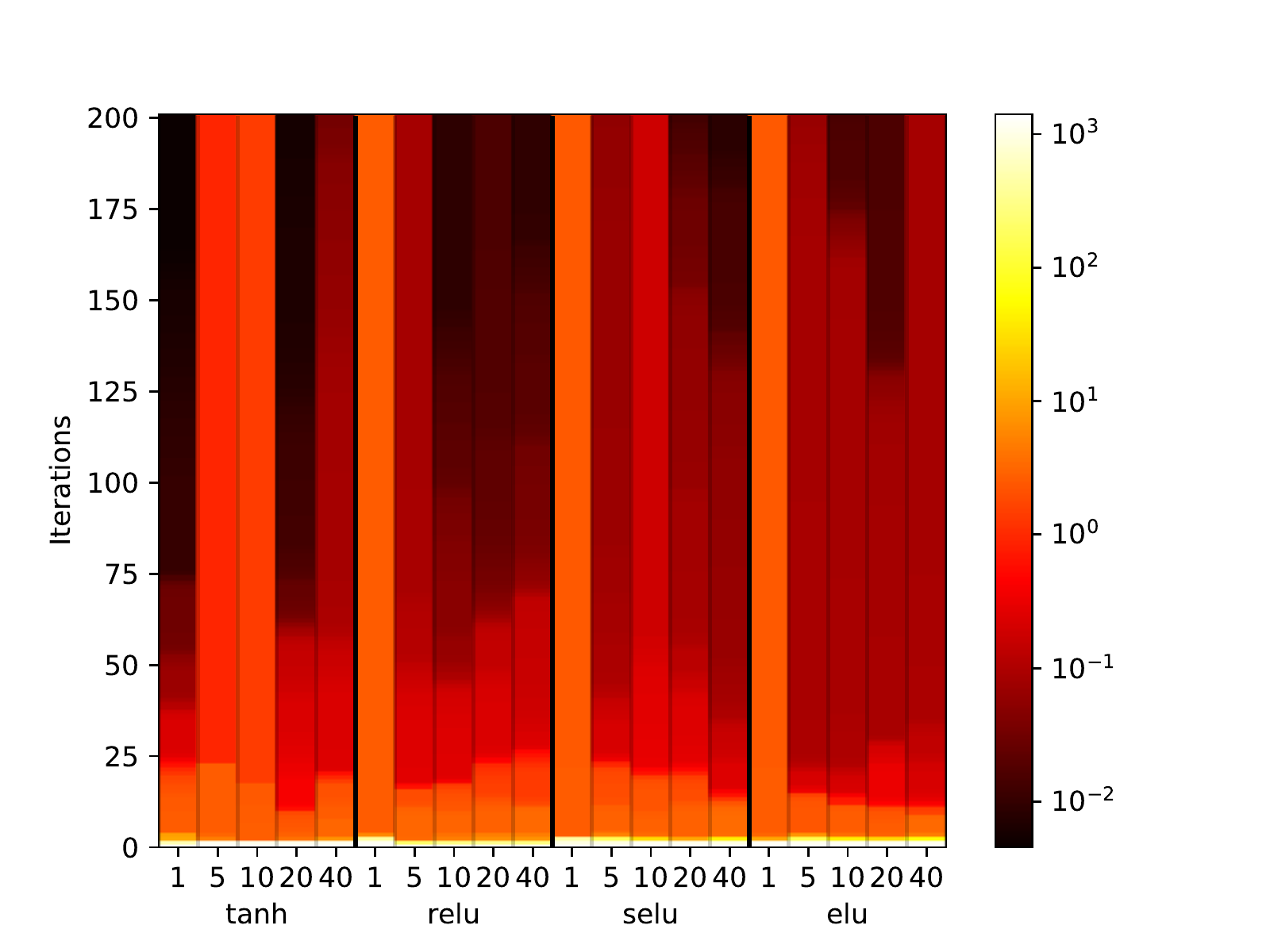}
	\caption{$\eta_1$, depth $=1$}
\end{subfigure}~
	\begin{subfigure}[t]{0.48\textwidth}
\centering
	\includegraphics[width=1.0\textwidth]{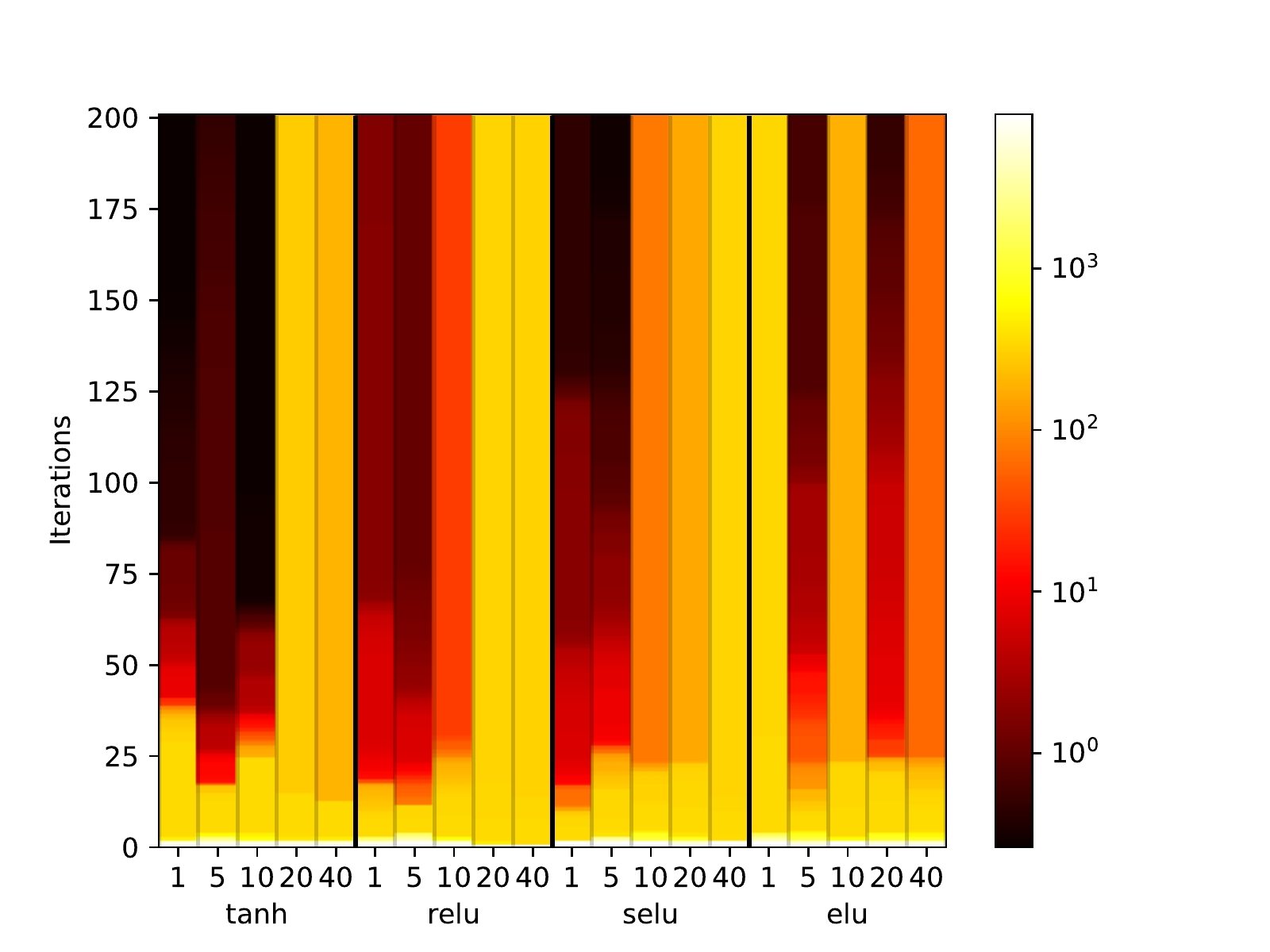}
	\caption{$\eta_2$, depth $=1$}
\end{subfigure}

\begin{subfigure}[t]{0.48\textwidth}
\centering
	\includegraphics[width=1.0\textwidth]{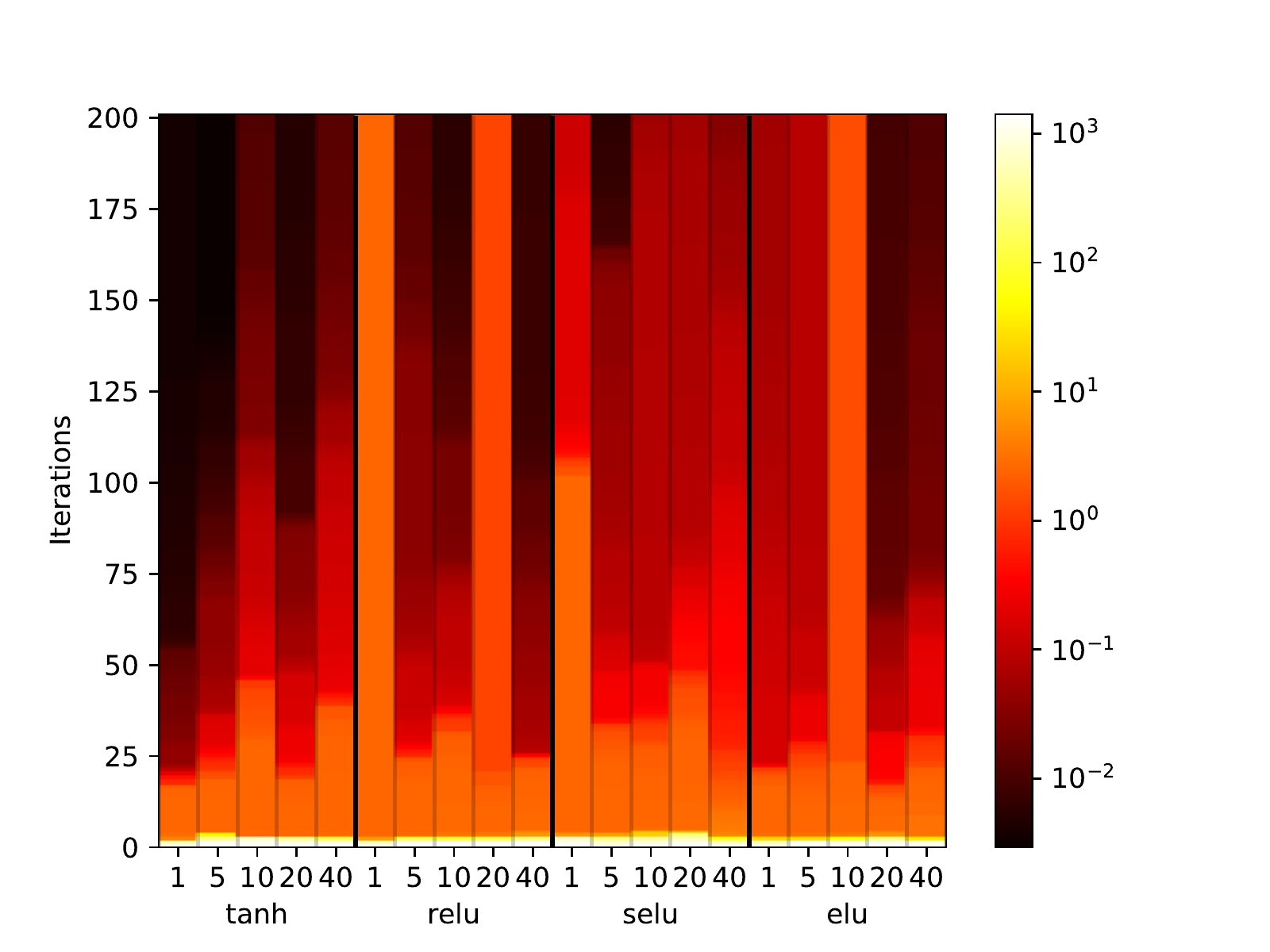}
	\caption{$\eta_1$, depth $=5$}
\end{subfigure}~
	\begin{subfigure}[t]{0.48\textwidth}
\centering
	\includegraphics[width=1.0\textwidth]{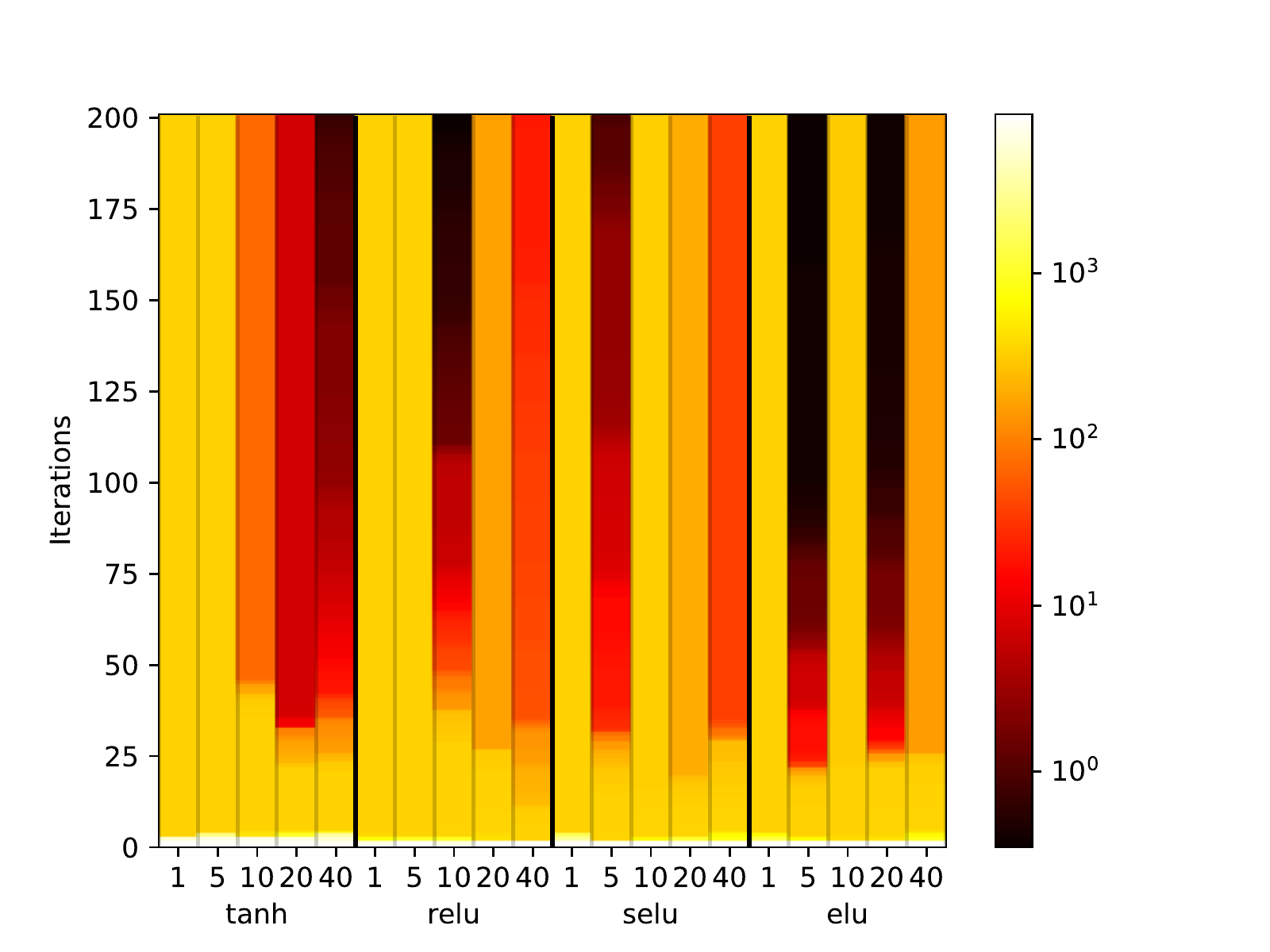}
	\caption{$\eta_2$, depth $=5$}
\end{subfigure}
\caption{Convergence of loss functions using different neural network architectures. The minor ticks 1, 5, 10, 20, and 40 represent the width of the neural network, i.e., number of neurons per layer.}
\end{figure}

\newpage

\bibliographystyle{unsrt}
\bibliography{ref}
\end{document}